\documentclass[a4paper,10pt,reqno]{amsart}

\usepackage{lineno}

\usepackage{BOONDOX-cal}
\usepackage{amssymb}
\usepackage{amsmath}
\usepackage[mathscr]{eucal}
\usepackage{mathrsfs}
\usepackage{xspace} 
\usepackage{tabu}
\usepackage{thmtools}
\usepackage[latin1]{inputenc}
\usepackage[english]{babel}
\usepackage[all,2cell]{xy}
\usepackage{xcolor}
\usepackage{graphicx}
\usepackage{wasysym} 
\usepackage{hyperref}

\usepackage{pdflscape}
\usepackage{caption}
\usepackage{multicol}
\usepackage{longtable}

\makeatletter
\newcommand{\leqnomode}{\tagsleft@true}
\newcommand{\reqnomode}{\tagsleft@false}
\makeatother

\usepackage{array}
\newcolumntype{L}[1]{>{\raggedright\let\newline\\\arraybackslash\hspace{0pt}}m{#1}}
\newcolumntype{C}[1]{>{\centering\let\newline\\\arraybackslash\hspace{0pt}}m{#1}}
\newcolumntype{R}[1]{>{\raggedleft\let\newline\\\arraybackslash\hspace{0pt}}m{#1}}

\definecolor{sg}{HTML}{df78ef}
\definecolor{sg1}{HTML}{ab47bc}
\definecolor{sg2}{HTML}{790e8b}

\usepackage{tikz}
\usetikzlibrary{mindmap,backgrounds, calc}
\usetikzlibrary{matrix,chains,scopes,positioning,arrows,fit}
\usetikzlibrary{positioning, shapes, patterns}
\usetikzlibrary{automata} 
\usetikzlibrary{shapes.geometric, arrows,arrows.meta}
\usetikzlibrary{positioning,decorations.markings}
\usetikzlibrary{cd}
\usepackage{hyperref}
\usepackage[acronym]{glossaries}

\usepackage{imakeidx}
\makeindex

\usepackage{stmaryrd} 
\UseAllTwocells

\newcommand{\at}[1]%
            {\ensuremath{\protect\underline{\mathbf{#1}}}} 

\newcommand{\op}[1]{\ensuremath{\operatorname{#1}}}        

\newcommand{\h}[1][]                                       
 {\ifthenelse{\boolean{mmode}}%
  {$\mathrm{h}$}%
  {h\nobreakdash#1\hspace{0pt}}}

\newcommand{\comp}{\circ}          
\newcommand{\adcomp}%
  {\overset{\operatorname{ad}}{\comp}} 
\newcommand{\funcomp}%
  {\overset{\operatorname{fn}}{\comp}}

\newcommand{\sccat}
{\mathbin{\kern-1pt\raisebox{6pt}{.}\kern-5pt
\downarrow\kern-5pt\raisebox{6pt}{.}\kern-1pt}}

\newcommand{\parrow}[1]
   {\underset{{\displaystyle \raisebox{5pt}%
   {$\longleftarrow$}}}{\op{#1}}{\,}}
\newcommand{\iarrow}[1]
   {\underset{{\displaystyle \raisebox{5pt}%
   {$\longrightarrow$}}}{\op{#1}}{\,}}

\newcommand{\rest}%
{\mathnormal{\restriction}}        

\newcommand{\concat}
  {\ensuremath{\text
  {\Large $\curlywedge$}}}

\newcommand{\brel}{\ensuremath{\xymatrix{{}\arity@{{*}{-}{*}}[r] & {}}}}

\newcommand{\nseq}[3]{\xymatrix@1@C=16pt{#1 \arity@{>}[r]_-{\scriptscriptstyle{#2}} & #3 }}

\NoCompileMatrices              
\xymatrixrowsep = {8ex}         
\xymatrixcolsep = {10ex}        

\newdir{<= }{:a(+25)\dir{-}*:a(-25)\dir{-}*!/:a(85)+1.8pt/:a(-25)\dir{-}}
\newdir{ = }{*@{=}}
\newdir{+>}{@{|}*@{-}!/-6.5pt/@{>}*!/-13pt/{}}
\newdir{ +}{{}*!/-5pt/@{+}}
\newdir{ >}{{}*!/-10pt/@{>}}
\newdir{.-}{{}*:a(-90)h!/2.1pt/{.}}
\newdir{.->}{:a(-90)h!/2.5pt/{.}*!/4pt/\dir{-}*!/0pt/\dir{-}*%
!/-2pt/\dir{-}*!/-8pt/\dir{>}}
\newdir{-.->}%
{:a(-90)h!/2pt/{.}*!/8pt/\dir{-}*!/4pt/\dir{-}*!/0pt/\dir{-}*%
!/-5pt/\dir{-}*!/-11pt/\dir{>}}
\newdir{~>}{!/4.6pt/\dir{~}*!/1pt/\dir{-}*!/-5pt/\dir{>}}
\newdir{~~>}{!/7pt/\dir{~}*!/0pt/\dir{~}*!/-5pt/\dir{>}}
\newdir{=~>}{*!/1pt/@2{~}*!/6pt/{}*!/-6pt/@2{>}}
\newdir{<~~>}{!/21pt/\dir{<}*!/21pt/\dir{-}!/12pt/\dir{~}!/5pt/\dir{~}*!/1pt/\dir{-}*!/-3pt/\dir{>}}
\newdir{=>}{!/5pt/\dir{=}!/2.5pt/\dir{=}*!/-5pt/\dir2{>}*!/7pt/\dir{ }}
\newdir{==>}{!/-2pt/\dir{=}!/-6pt/\dir{=}%
 !/-10pt/\dir{=}*!/-18pt/\dir2{>}}
\newdir{<==>}{!\dir2{<}!/-2pt/\dir{=}!/-6pt/\dir{=}%
 !/-10pt/\dir{=}*!/-17pt/\dir2{>}}
\newdir{< }{{}*!/12pt/\dir{<}}
\newdir{-o}{!/1.3pt/{}*!/0pt/{}@{o}*!/-5pt/{}}

\newsavebox{\xymor}  
\newsavebox{\xymon}  
\newsavebox{\xyepi}  
\newsavebox{\xytn}   
\newsavebox{\xyrel}  
\newsavebox{\xycel}  
\newsavebox{\xymdf}  
\newsavebox{\xyumor} 
\newsavebox{\xydmor} 
\newsavebox{\xyomor} 
\newsavebox{\xyemor} 

\newcommand{\xynode}{\makebox[0ex]{}}
\savebox{\xymor}{\ensuremath{%
\xymatrix@1@C=19pt{\xynode \ar@{>}[r] & \xynode }}}
\savebox{\xymon}{\ensuremath{%
\xymatrix@1@C=19pt{\xynode \ar@{{ +}{-}{>}}[r] & \xynode }}}
\savebox{\xyepi}{\ensuremath{%
\xymatrix@1@C=19pt{\xynode \ar@{{}{-}{+>}}[r] & \xynode }}}
\savebox{\xytn}{\ensuremath{%
\xymatrix@1@C=19pt{\xynode \ar[r]|(.44){\object@{.-}} & \xynode
}}}
\savebox{\xyrel}{\ensuremath{%
\xymatrix@1@C=19pt{\xynode \ar@{{}{-}{-o}}[r] & \xynode }}}
\savebox{\xycel}{\ensuremath{%
\xymatrix@1@C=19pt{\xynode \ar@{=>}[r] & \xynode }}}
\savebox{\xymdf}{\ensuremath{%
\xymatrix@1@C=16pt{\xynode \ar@{}[r]|{\dir{~>}} & \xynode}}}
\savebox{\xyumor}{\ensuremath{%
\xymatrix@1@C=19pt{\xynode \ar@{{}{-}^{>}}[r] & \xynode }}}
\savebox{\xydmor}{\ensuremath{%
\xymatrix@1@C=19pt{\xynode \ar@{{}{-}_{>}}[r] & \xynode }}}
\savebox{\xyomor}{\ensuremath{%
\xymatrix@1@C=19pt{\xynode \ar@{{}{-}^{< }}[r] & \xynode }}}
\savebox{\xyemor}{\ensuremath{%
\xymatrix@1@C=19pt{\xynode \ar@{{ >}{-}{>}}[r] & \xynode }}}

\newcommand{\mor}{\usebox{\xymor}}    

\newcommand{\functor}[9]{
 \xymatrix{
    #4 \save[]+<0ex,5ex>*+{#1}="1"  \restore
      \arity[d]_{#6}  \arity@{}[rd]|{\longmapsto}
  & #5 \save[]+<0ex,5ex>*+{#3}="3"  \restore
      \arity[d]^{#7}
  \\
   #8 & #9 \arity "1";"3"^-{#2} } }
\newcommand{\functornd}[9]{
 \xymatrix{
    #4 \save[]+<0ex,5ex>*+{#1}="1"  \restore
      \arity[d]_{#6}  \arity@{}[rd]|{\longmapsto}
  & #5 \save[]+<0ex,5ex>*+{#3}="3"  \restore
  \\
   #8 & #9 \arity[u]_{#7} \arity "1";"3"^-{#2} } }
\newcommand{\functordn}[9]{
 \xymatrix{
    #4 \save[]+<0ex,5ex>*+{#1}="1"  \restore
       \arity@{}[rd]|{\longmapsto}
  & #5 \save[]+<0ex,5ex>*+{#3}="3"  \restore
      \arity[d]^{#7}
  \\
   #8  \arity[u]^{#6}  & #9 \arity "1";"3"^-{#2} } }

\newcommand{\larr}{->}
\newcommand{\rarr}{->}

\newcommand{\xfunctor}[9]{
 \xymatrix{
    #4 \save[]+<0ex,5ex>*+{#1}="1"  \restore
      \ifthenelse{\equal{\larr}{->}}{\arity[d]_{#6}}{}
      \ifthenelse{\equal{\larr}{<-}}{\arity[d];[]^{#6}}{}
      \ifthenelse{\equal{\larr}{-<}}{\arity@{< }[d]_{#6}}{}
      \arity@{}[rd]|{\longmapsto}
  & #5 \save[]+<0ex,5ex>*+{#3}="3"  \restore
      \ifthenelse{\equal{\rarr}{->}}{\arity[d]^{#7}}{}
      \ifthenelse{\equal{\rarr}{<-}}{\arity[d];[]_{#7}}{}
      \ifthenelse{\equal{\rarr}{-<}}{\arity@{< }[d]^{#7}}{}
  \\
   #8 & #9 \arity "1";"3"^-{#2} } }

\theoremstyle{plain}
\newtheorem{theorem}{Theorem}[section]
\newtheorem{proposition}[theorem]{Proposition}

\newtheorem{corollary}[theorem]{Corollary}
\newtheorem{lemma}[theorem]{Lemma}
\newtheorem{remark}[theorem]{Remark}
\theoremstyle{definition}
\newtheorem{definition}[theorem]{Definition}

\newtheorem*{assumption}{Assumption}

\theoremstyle{remark}

\newcommand{\arity}{\mathsf{ar}}

\newtheorem{claim}[theorem]{Claim}

\numberwithin{equation}{section}

\begin{document}
\title[P{\l}onka construction has a right adjoint]{P{\l}onka adjunction}

\author[Climent]{J. Climent Vidal}
\address{Universitat de Val\`{e}ncia\\
         Departament de L\`{o}gica i Filosofia de la Ci\`{e}ncia\\
         Av. Blasco Ib\'{a}\~{n}ez, 30-$7^{\mathrm{a}}$, 46010 Val\`{e}ncia, Spain}
\email{Juan.B.Climent@uv.es}
\author[Cosme]{E. Cosme Ll\'{o}pez}
\address{Universitat de Val\`{e}ncia\\
         Departament d'\`{A}lgebra\\
         Dr. Moliner, 50, 46100 Burjassot, Val\`{e}ncia, Spain}
\email{enric.cosme@uv.es}

\subjclass[2020]{Primary: 08A62, 08B25, 08C05, 18A40. Secondary: 08A30, 18A30.} 
\keywords{P{\l}onka sum, P{\l}onka algebra, sup-semilattice inductive system of algebras, Grothendieck construction, strong Lawvere adjoint cylinder}
\date{June 19th, 2022}

\begin{abstract}


For a signature $\Sigma$ and its subsignature $\Sigma^{\neq 0}$ without $0$-ary operation symbols, we prove (1) that there are strong Lawvere adjoint cylinders between the category $\mathsf{Ssl}$, of sup-semilattices, and the categories $\int^{\mathsf{Ssl}}\mathrm{Isys}_{\Sigma}$, of sup-semilattice inductive systems of $\Sigma$-algebras, and $\int^{\mathsf{Ssl}}\mathrm{Isys}_{\Sigma^{\neq 0}}$, of sup-semilattice inductive systems of $\Sigma^{\neq 0}$-algebras; (2) that  there exists an adjunction between $\mathsf{Ssl}$ and the category $\mathsf{Alg}(\Sigma^{\neq 0})$, of $\Sigma^{\neq 0}$-algebras; (3) that there exists an adjunction between the categories 
$\mathsf{Ssl}$ and $\mathsf{Lnb}$, the category of left normal bands; (4) after defining and stating several technical results on the category $\mbox{\sffamily{\upshape{P{\l}Alg}}}(\Sigma^{\neq 0})$, of P{\l}onka $\Sigma^{\neq 0}$-algebras, and defining functors $J_{\Sigma^{\neq 0}}$ from $\mbox{\sffamily{\upshape{P{\l}Alg}}}(\Sigma^{\neq 0})$ to  $\mathsf{Alg}(\Sigma^{\neq 0})\otimes\mathsf{Lnb}$, the tensor product of $\mathsf{Alg}(\Sigma^{\neq 0})$ and $\mathsf{Lnb}$, and $P_{\Sigma^{\neq 0}}$ from $\mathsf{Alg}(\Sigma^{\neq 0})\otimes\mathsf{Lnb}$ to $\mathsf{Alg}(\Sigma^{\neq 0})$, we prove that $P_{\Sigma^{\neq 0}}\circ J_{\Sigma^{\neq 0}}$ has a left adjoint; finally, (5) after defining a functor $\mathrm{Is}_{\Sigma^{\neq 0}}$ from 
$\mbox{\sffamily{\upshape{P{\l}Alg}}}(\Sigma^{\neq 0})$ to $\int^{\mathsf{Ssl}}\mathrm{Isys}_{\Sigma^{\neq 0}}$ we prove the main result of this paper: that $\mathrm{Is}_{\Sigma^{\neq 0}}$ has a left adjoint $\mbox{\upshape{P{\l}}}_{\Sigma^{\neq 0}}$, which is the P{\l}onka sum.
\end{abstract}

\maketitle

\section{Introduction}

In the seminal paper~\cite{P67}, P{\l}onka introduced a new construction in \emph{universal algebra}, which is a very interesting particular case of the classical notion of inductive limit of an inductive system of algebras. Specifically, for a sup-semilattice inductive system of algebras, he defined its sum, later called the P{\l}onka sum, and proved, by means of the notion of partition function of the underlying set of an algebra without $0$-ary operation symbols, that the correspondence between partition functions and representations as the sum of a sup-semilattice inductive system of algebras is one-to-one. The notion of P{\l}onka sum has, as P{\l}onka himself says in~\cite{P67}, as its precursor the construction of a strong semilattice of semigroups introduced by Clfford in~\cite{C41}. P{\l}onka sums are applied to several algebraic theories, including semigroup theory and semiring theory, see, e.g.,~\cite{D03, PR82}. 
These investigations were further pursued, among others, by P{\l}onka himself in~\cite{ P67a,P68, P74}, by P{\l}onka and Romanowska in~\cite{PR92}, by  Romanowska and Smith in~\cite{RS91} and, in a massive way, by Romanowska and Smith in~\cite{RS85, RS02} (for an almost exhaustive bibliography on the subject we recommend the ones included in~\cite{BPPB22, gG08, RS02}; by the way, in~\cite{BPPB22} there are applications of P{\l}onka sums to logic).

Our main goal in this paper is to provide a category--theoretic investigation of the notion of 
P{\l}onka sum, which we also call the P{\l}onka construction. Concretely, after setting up the necessary notions and constructions specified below, we prove that the  P{\l}onka construction is part of an adjunction.

The paper is organized as follows. In Section~2 we recall those mathematical concepts and constructions from the fields of universal algebra, lattice theory and category theory which will be necessary to  understand the rest of the paper. Specifically, we begin by defining the concepts of signature, algebra, subalgebra of an algebra, congruence on an algebra, quotient of an algebra by a congruence on it. We next associate to a signature $\Sigma$ the subsignature
obtained from it by removing its $0$-ary operation symbols, which we denote by $\Sigma^{\neq 0}$, and state the corresponding adjuction between the categories of algebras associated with such a pair of signatures. We then proceed to define several key concepts of lattice theory, including residuated morphisms between ordered sets, sup-semilattices, and morphisms between sup-semilattices. Following this, for a category $\mathsf{I}$, we consider the notion of split $\mathsf{I}$-indexed category and the Grothendieck construction. Finally, we state a sufficient condition for this construction to have a left adjoint. 

In Section~3, for a signature $\Sigma$ we introduce the concept of sup-semilattice inductive system of $\Sigma$-algebras relative to a sup-semilattice and of morphism between such systems. We then define a contravariant functor $\mathrm{Isys}_{\Sigma}$ from the category $\mathsf{Ssl}$, of sup-semilattices, to the category $\mathsf{Cat}$, of categories and functors. Next, we use the Grothendieck construction to obtain the category $\int^{\mathsf{Ssl}}\mathrm{Isys}_{\Sigma}$, which has sup-semilattice inductive systems of $\Sigma$-algebras as objects. We then define functors $L_{\Sigma}$ and $K{_\Sigma}$ from $\mathsf{Ssl}$ to $\int^{\mathsf{Ssl}}\mathrm{Isys}_{\Sigma}$, such that $L_{\Sigma}\dashv \pi_{\mathrm{Isys}_{\Sigma}}\dashv K_{\Sigma}$, where $\pi_{\mathrm{Isys_{\Sigma}}}$ is the canonical split fibration  which sends each inductive system to its sup-semilattice. After that we prove that the ordered triple $(L_{\Sigma},\pi_{\mathrm{Isys}_{\Sigma}},K_{\Sigma})$ is a strong Lawvere adjoint cylinder and we prove that there exists a natural transformation $\gamma_{\Sigma}$ from $L_{\Sigma}$ to $K_{\Sigma}$. We then show that the category $\int^{\mathsf{Ssl}}\mathrm{Isys}_{\Sigma}$ is complete and, for the wide subcategory $\mathsf{Ssl}_{\mathrm{rd}}$ of $\mathsf{Ssl}$, consisting of sup-semilattices and residuated morphisms, the category $\int^{\mathsf{Ssl}_{\mathrm{rd}}}\mathrm{Isys}_{\Sigma}$ is cocomplete. Moreover, we point out that, for the subsignature $\Sigma^{\neq 0}$ of $\Sigma$, these results are also valid for the contravariant functor $\mathrm{Isys}_{\Sigma^{\neq 0}}$ (defined in the same way that $\mathrm{Isys}_{\Sigma}$) and the corresponding category $\int^{\mathsf{Ssl}}\mathrm{Isys}_{\Sigma^{\neq 0}}$.

In Section~4, we prove that there exists a functor $W_{\Sigma^{\neq 0}}$ from $\mathsf{Ssl}$ to $\mathsf{Alg}(\Sigma^{\neq 0})$ which has a left adjoint.

In Section~5, we introduce the category $\mathsf{Lnb}$ of left normal bands and, for each left normal band, we define, by recursion, the family of iterates of its structural operation. We also establish several useful properties of this family for later use. Next, we prove that there exists  an adjunction between the categories $\mathsf{Ssl}$ and $\mathsf{Lnb}$. To do it, we begin by showing that every sup-semilattice is a left-normal band. Then we prove that the inclusion functor from $\mathsf{Ssl}$ to $\mathsf{Lnb}$ has a left adjoint.

In Section~6, after defining the category $\mbox{\sffamily{\upshape{P{\l}Alg}}}(\Sigma^{\neq 0})$,  of P{\l}onka $\Sigma^{\neq 0}$-algebras, we prove several technical results on the underlying P{\l}onka operator of a P{\l}onka $\Sigma^{\neq 0}$-algebra, which will be fundamental to prove that the P{\l}onka construction is part of an adjunction. Moreover, after defining functors $J_{\Sigma^{\neq 0}}$ from $\mbox{\sffamily{\upshape{P{\l}Alg}}}(\Sigma^{\neq 0})$ to  
$\mathsf{Alg}(\Sigma^{\neq 0})\otimes\mathsf{Lnb}$, the tensor product of 
$\mathsf{Alg}(\Sigma^{\neq 0})$ and $\mathsf{Lnb}$, and $P_{\Sigma^{\neq 0}}$ from $\mathsf{Alg}(\Sigma^{\neq 0})\otimes\mathsf{Lnb}$ to $\mathsf{Alg}(\Sigma^{\neq 0})$, we prove that $P_{\Sigma^{\neq 0}}\circ J_{\Sigma^{\neq 0}}$ has a left adjoint.

Finally, in Section~7, after proving that there exists a functor $\mathrm{Is}_{\Sigma^{\neq 0}}$ from $\mbox{\sffamily{\upshape{P{\l}Alg}}}(\Sigma^{\neq 0})$ to $\int^{\mathsf{Ssl}}\mathrm{Isys}_{\Sigma^{\neq 0}}$, we prove the main result of this paper: that $\mathrm{Is}_{\Sigma^{\neq 0}}$ has a left adjoint $\mbox{\upshape{P{\l}}}_{\Sigma^{\neq 0}}$, the P{\l}onka construction, whose object part is the one that assigns to a sup-semilattice inductive system of $\Sigma^{\neq 0}$-algebras its P{\l}onka sum.  

Our underlying set theory is $\mathbf{ZFSk}$, Zermelo-Fraenkel-Skolem set theory (also known as $\mathbf{ZFC}$, i.e., Zermelo-Fraenkel set theory with the axiom of choice) plus the existence of a Grothendieck universe $\mathbf{U}$, fixed once and for all (see~\cite[pp.~21--24]{sM98}).

In all that follows we use standard concepts and constructions from set theory, see, e.g.~\cite{nB70, end77}; universal algebra, see, e.g.~\cite{gb15, bl05, bs81, gG08, w92}; lattice theory, see, e.g.,~\cite{bl05, HK71}; semigroup theory, see, e.g.,~\cite{JH80}; and category theory, see, e.g.~\cite{BW85, jb68, Gro71, hs73, L94, sM98, em76, tbg91}. Nevertheless, regarding set theory, we have adopted the following conventions. An \emph{ordinal} $\alpha$ is a transitive set that is well-ordered by $\in$; thus, $\alpha=\{\beta \mid \beta \in \alpha\}$. The first transfinite ordinal $\omega_{0}$ will be denoted by $\mathbb{N}$, which is the set of all \emph{natural numbers}, and, from what we have just said about the ordinals, for every $n\in\mathbb{N}$, $n=\{0,\dots, n-1\}$. A function from $A$ to $B$ is a subset $F$ of $A\times B$ satisfying the functional condition and a mapping from $A$ to $B$ is an ordered triple $f=(A,F,B)$, denoted by $f\colon A\mor B$, in which $F$ is a function from $A$ to $B$.

\section{Preliminaries}

In this section we collect those basic facts about universal algebra, semigroup theory, lattice theory and category theory, that we will need to carry out our research. For a full treatment of the topics of this section we refer the reader to \cite{BW85, gb15, bl05, bs81, gG08, Gro71, em76, w92}.


We begin by giving a precise definition of the concept of signature.

\begin{definition}\label{DSig}
A \emph{signature} is a mapping $\Sigma = (\Sigma_{n})_{n\in\mathbb{N}}$ from $\mathbb{N}$ to 
$\mathbf{U}$ which sends a natural number $n\in\mathbb{N}$ to the set $\Sigma_{n}$ of the \emph{formal operations}, or \emph{operation symbols}, $\sigma$ of \emph{arity} $n$, denoted by $\mathrm{ar}(\sigma) = n$.
\end{definition}

\begin{assumption}
From now on, unless otherwise stated, $\Sigma$ stands for a signature as set forth in Definition~\ref{DSig}, fixed once and for all.
\end{assumption}

We shall now give precise definitions of the concepts of algebra and homomorphism between algebras.

\begin{definition}\label{DAlg}
The set of the \emph{finitary operations on} a set $A$ is $(\mathrm{Hom}(A^{n},A))_{n\in\mathbb{N}}$, where, for every $n\in\mathbb{N}$, $A^n = \prod_{j\in n}A$ (if $n = 0$, then $A^{0}$ is a final set, i.e., a singleton set). A \emph{structure of} $\Sigma$-\emph{algebra on} a  set  $A$ is a family $(F_{n})_{n\in \mathbb{N}}$, denoted by $F$, where, for $n\in \mathbb{N}$, $F_{n}$ is a mapping from $\Sigma_{n}$ to $\mathrm{Hom}(A^{n},A)$ (if $n=0$ and $\sigma\in \Sigma_{0}$, then $F_{0}(\sigma)$ picks out an element of $A$). For a natural number $n\in \mathbb{N}$ and a formal operation $\sigma\in \Sigma_{n}$, in order to simplify the notation, the operation $F_{n}(\sigma)$ from $A^{n}$ to $A$ will be written as $F_{\sigma}$. A $\Sigma$-\emph{algebra} is a pair $(A,F)$, denoted by $\mathbf{A}$, where $A$ is a set and $F$ a structure of $\Sigma$-algebra on $A$. A $\Sigma$-\emph{homomorphism} (or, to abbreviate, \emph{homomorphism}) from $\mathbf{A}$ to $\mathbf{B}$, where $\mathbf{B} = (B,G)$, is a triple $(\mathbf{A},f,\mathbf{B})$, denoted by $f\colon \mathbf{A}\mor \mathbf{B}$, where $f$ is a  mapping from $A$ to $B$ such that, for every $n\in\mathbb{N}$, every  $\sigma\in \Sigma_{n}$ and every $(a_{j})_{j\in n}\in A^{n}$, we have 
$$
f(F_{\sigma}((a_{j})_{j\in n})) = G_{\sigma}((f(a_{j}))_{j\in n}).
$$
We will denote by $\mathsf{Alg}(\Sigma)$ the category of $\Sigma$-algebras and homomorphisms and by $\mathrm{Alg}(\Sigma)$ the set of objects  of $\mathsf{Alg}(\Sigma)$.

In some cases, to avoid mistakes, we will denote by $F^{\mathbf{A}}$ the structure of $\Sigma$-algebra on $A$, and, for $n\in \mathbb{N}$ and $\sigma\in \Sigma_{n}$, by $F^{\mathbf{A}}_{\sigma}$, or simply by $\sigma^{\mathbf{A}}$, the corresponding operation. Moreover, for $\sigma\in\Sigma_{0}$, we will, usually, denote by $\sigma^{\mathbf{A}}$ the value of the mapping $F^{\mathbf{A}}_{\sigma}\colon A^{0}\mor A$ at the unique element in $A^{0}$.

We will denote by $\mathbf{1}$ the (standard) final $\Sigma$-algebra. Moreover, if $\Sigma_{0} = \varnothing$, then we let $\pmb{\varnothing}$ stand for the initial $\Sigma$-algebra.
\end{definition}

We shall now give the definitions of the concepts of closed subset and of subalgebra of an algebra.

\begin{definition}\label{DSubAlg}
Let $\mathbf{A} = (A,F)$ be a $\Sigma$-algebra and $X\subseteq A$. Given $n\in \mathbb{N}$ and $\sigma\in\Sigma_{n}$, we will say that $X$ is \emph{closed under the operation} $F_{\sigma}\colon A^{n}\mor A$ if, for every $(a_{j})_{j\in n}\in X^{n}$, $F_{\sigma}((a_{j})_{j\in n})\in X$. We will say that $X$ is a \emph{closed subset} of $\mathbf{A}$ if $X$ is closed under the operations of $\mathbf{A}$. We will denote by $\mathrm{Cl}(\mathbf{A})$ the set of all closed subsets of $\mathbf{A}$ (which is an algebraic closure system on $A$) and by $\mathbf{Cl}(\mathbf{A})$ the algebraic lattice $(\mathrm{Cl}(\mathbf{A}),\subseteq)$. We will say that a $\Sigma$-algebra $\mathbf{B}$ is a \emph{subalgebra} of $\mathbf{A}$ if $B\subseteq A$ and 
$\mathrm{in}_{B,A}$, the canonical embedding of $B$ into $A$, determines an injective homomorphism 
$\mathrm{in}_{\mathbf{B},\mathbf{A}} = (\mathbf{B},\mathrm{in}_{B,A},\mathbf{A})$ from $\mathbf{B}$ to $\mathbf{A}$. We will denote by $\mathrm{Sub}(\mathbf{A})$ the set of all subalgebras of $\mathbf{A}$. Since $\mathrm{Cl}(\mathbf{A})$ and $\mathrm{Sub}(\mathbf{A})$ are isomorphic, we shall feel free to deal either with a closed subset of $\mathbf{A}$ or with the correlated subalgebra of $\mathbf{A}$, whichever is most convenient for the work at hand.
\end{definition}

We next give the definitions of the concepts of congruence on an algebra and of quotient of an algebra by a congruence on it.

\begin{definition}\label{DCong}
Let $\mathbf{A}$ be a $\Sigma$-algebra and $\Phi$ an equivalence on $A$. We will say that $\Phi$ is a \emph{congruence on} $\mathbf{A}$ if, for every $n\in \mathbb{N}$, every $\sigma\in \Sigma_{n}$,
and every $(a_{j})_{j\in n}, (b_{j})_{j\in n}\in A^{n}$, if, for every $j\in  n$, $(a_{j}, b_{j})\in\Phi$, then $$(F_{\sigma}((a_{j})_{j\in n}), F_{\sigma}((b_{j})_{j\in n}))\in \Phi.$$
We will denote by $\mathrm{Cgr}(\mathbf{A})$ the set of all congruences on $\mathbf{A}$ (which is an algebraic closure system on $A\times A$), by $\mathbf{Cgr}(\mathbf{A})$ the algebraic lattice $(\mathrm{Cgr}(\mathbf{A}),\subseteq)$, by $\nabla_{\mathbf{A}}$ the greatest element of $\mathbf{Cgr}(\mathbf{A})$ and by $\Delta_{\mathbf{A}}$ the least element of $\mathbf{Cgr}(\mathbf{A})$.

For a congruence $\Phi$ on $\mathbf{A}$, the \emph{quotient $\Sigma$-algebra of} $\mathbf{A}$ \emph{by} $\Phi$, denoted by $\mathbf{A}/\Phi$, is the $\Sigma$-algebra $(A/\Phi, F^{\mathbf{A}/\Phi})$, where, for every  $n\in \mathbb{N}$ and every $\sigma\in \Sigma_{n}$, the operation $F_{\sigma}^{\mathbf{A}/\Phi}$ from $(A/\Phi)^{n}$ to $A/\Phi$, sends $([a_{j}]_{\Phi})_{j\in n}$ in $(A/\Phi)^{n}$ to $[F^{\mathbf{A}}_{\sigma}((a_{j})_{j\in n})]_{\Phi}$ in $A/\Phi$,
and the \emph{canonical projection} from $\mathbf{A}$ to $\mathbf{A}/\Phi$, denoted by $\mathrm{pr}_{\Phi}\colon \mathbf{A}\mor \mathbf{A}/\Phi$, is the surjective homomorphism determined by the projection from $A$ to $A/\Phi$.
The ordered pair $(\mathbf{A}/\Phi,\mathrm{pr}_{\Phi})$ has the following universal property: $\mathrm{Ker}(\mathrm{pr}_{\Phi})$ is $\Phi$ and, for every $\Sigma$-algebra $\mathbf{B}$ and every homomorphism $f$ from $\mathbf{A}$ to $\mathbf{B}$, if $\Phi\subseteq \mathrm{Ker}(f)$, then there exists a unique homomorphism $h$ from $\mathbf{A}/\Phi$ to $\mathbf{B}$ such that $h\circ\mathrm{pr}_{\Phi} = f$. In particular, if $\Psi$ is a congruence on $A$ such that $\Phi\subseteq \Psi$, then we will denote by $\mathrm{p}_{\Phi,\Psi}$ the unique homomorphism from $\mathbf{A}/\Phi$ to $\mathbf{A}/\Psi$ such that $\mathrm{p}_{\Phi,\Psi}\circ \mathrm{pr}_{\Phi} = \mathrm{pr}_{\Psi}$.
\end{definition}

We next state that the forgetful functor $\mathrm{G}_{\Sigma}$ from $\mathsf{Alg}(\Sigma)$ to
$\mathsf{Set}$ has a left adjoint, namely $\mathbf{T}_{\Sigma}$, which assigns to a set $X$ the free $\Sigma$-algebra $\mathbf{T}_{\Sigma}(X)$ on $X$. We later state the universal property of the free $\Sigma$-algebra.

Let us note that in what follows, to construct the algebra of $\Sigma$-rows in $X$, and the free $\Sigma$-algebra on $X$, since neither the signature $\Sigma$ nor the set $X$ are subject to any constraint, coproducts must necessarily be used.

\begin{definition}
Let $X$ be a set. The \emph{algebra of} $\Sigma$-\emph{rows in} $X$, denoted by $\mathbf{W}_{\Sigma}(X)$, is defined as follows:
\begin{enumerate}
\item The underlying set of $\mathbf{W}_{\Sigma}(X)$, written as $W_{\Sigma}(X)$, is the set $(\coprod\Sigma \amalg X)^{\star}$, where $(\coprod\Sigma \amalg X)^{\star}$ is the set of all words on the set $\coprod\Sigma \amalg X$, i.e., on the set
      $$
      \textstyle
      [(\bigcup_{n\in\mathbb{N}}(\Sigma_{n}\times\{n\}))\times1]\cup
      [X\times\{1\}].
      $$

\item For every $n\in \mathbb{N}$ and every $\sigma\in\Sigma_{n}$, the structural operation $F_{\sigma}$ associated to $\sigma$ is the mapping from $W_{\Sigma}(X)^{n}$ to $W_{\Sigma}(X)$ which sends $(P_{j})_{j\in n} \in W_{\Sigma}(X)^{n}$ to $(\sigma)\curlywedge\concat_{j\in n}P_{j} \in {W_{\Sigma}(X)}$---the concatenation of the word $(\sigma)$ and the word $\concat_{j\in n}P_{j}$,  which is in its turn the concatenation of the words $P_{j}$ of the family of words $(P_{j})_{j\in n}$--- where, for every $n\in \mathbb{N}$ and every $\sigma\in\Sigma_{n}$,
      $(\sigma)$ stands for $(((\sigma,n),0))$, which is the value at $\sigma$ of the canonical mapping from $\Sigma_{n}$ to $(\coprod\Sigma \amalg X)^{\star}$.
\end{enumerate}
\end{definition}

\begin{definition}\label{DFreeAlg}
The \emph{free} $\Sigma$-\emph{algebra on} a set $X$, denoted by $\mathbf{T}_{\Sigma}(X)$, is the $\Sigma$-algebra determined by $\mathrm{Sg}_{\mathbf{W}_{\Sigma}(X)}(\{(x)\mid x\in X\})$, the subalgebra of $\mathbf{W}_{\Sigma}(X)$ generated by $\{(x)\mid x\in X\}$, where, for every  $x\in X$, $(x)$ stands for $(x,1)$, which is the value at $x$ of the canonical mapping from $X$ to $(\coprod\Sigma \amalg  X)^{\star}$. We will denote by $\mathrm{T}_{\Sigma}(X)$ the underlying set of $\mathbf{T}_{\Sigma}(X)$ and we will call the elements of $\mathrm{T}_{\Sigma}(X)$ \emph{terms} \emph{with variables in} $X$  or $X$-\emph{terms}.
\end{definition}


We have the following characterization of the elements of $\mathrm{T}_{\Sigma}(X)$.

\begin{proposition}
Let $X$ be a set. Then, for every $P\in
W_{\Sigma}(X)$, we have that $P$ is a term  with variables in $X$ if, and only if, $P = (x)$, for a unique $x\in X$, or $P = (\sigma)$, for a unique $\sigma\in\Sigma_{n}$, or $P = (\sigma)\curlywedge\concat_{j\in n}P_{j}$, for a unique $n\in \mathbb{N}$, a unique $\sigma\in\Sigma_{n}$ and a unique family $(P_{j})_{j\in n}\in\mathrm{T}_{\Sigma}(X)^{n}$.
Moreover, the three possibilities are mutually exclusive. From now on, for simplicity of notation, we will write $x$, $\sigma$ and $\sigma(P_{0},\ldots,P_{n-1})$ or $\sigma((P_{j})_{j\in n})$ instead of $(x)$, $(\sigma)$ and $(\sigma)\curlywedge\concat_{j\in n}
P_{j}$, respectively. Thus, in particular, the structural operation $F_{\sigma}$ (or, more accurately, $F_{\sigma}^{\mathbf{T}_{\Sigma}(X)}$) associated to $\sigma$ is identified with $\sigma$.
\end{proposition}

From the above proposition it follows, immediately, the universal property of the free $\Sigma$-algebra on a set $X$, as stated in the subsequent proposition.

\begin{proposition}
For every set $X$, the pair $(\mathbf{T}_{\Sigma}(X),\eta_{X})$, where $\eta_{X}$, the
\emph{insertion of (the set of generators)} $X$ \emph{into} $\mathrm{T}_{\Sigma}(X)$, is the 
corestric\-tion to $\mathrm{T}_{\Sigma}(X)$ of the canonical embedding of $X$ into $W_{\Sigma}(X)$, has the following universal property: for every $\Sigma$-algebra $\mathbf{A}$ and every  mapping $f\colon X\mor A$, there exists a unique $\Sigma$-homomorphism $f^{\sharp}\colon\mathbf{T}_{\Sigma}(X)\mor\mathbf{A}$ such that $f^{\sharp}\circ \eta_{X} = f$.
\end{proposition}


For a signature $\Sigma$, we next define the signature $\Sigma^{\neq 0}$ which is obtained from $\Sigma$ by removing the $0$-ary operation symbols from it. 

\begin{definition}\label{DSigNoCons} Let $\Sigma$ be a signature. Then we let  
$\Sigma^{\neq 0}$ stand for the signature defined as follows:
\[
\Sigma^{\neq 0}_{n}=
\begin{cases}
\Sigma_{n},&\mbox{if }n\neq 0;\\
\varnothing,&\mbox{if }n=0.
\end{cases}
\]
\end{definition}

%
%


We next define the notions of residuated morphism between ordered sets, of sup-semilattice and of morphism between sup-semilattices.

\begin{definition}
Let $\mathbf{I} = (I,\leq)$ and $\mathbf{P} = (P,\leq)$ be ordered sets and $\xi$ an isotone mapping from 
$\mathbf{I}$ to $\mathbf{P}$. We will say that $\xi$ is \emph{residuated} if there exists an isotone mapping $\zeta$ from $\mathbf{P}$ to $\mathbf{I}$ such that, for every $i\in I$, $\zeta(\xi(i))\leq i$ and , for every $p\in P$, $p\leq\xi(\zeta(p))$.
\end{definition}

\begin{remark}
To say that the isotone mapping $\xi$ from $\mathbf{I}$ to $\mathbf{P}$ is residuated is, from a 
category-theoretic point of view, equivalent to say that $\xi$ has a left adjoint. Moreover, if the isotone mapping $\xi$ from $\mathbf{I}$ to $\mathbf{P}$ is residuated, then there exists a unique isotone mapping $\zeta$ from $\mathbf{P}$ to $\mathbf{I}$ such that, for every $i\in I$, $\zeta(\xi(i))\leq i$ and, for every $p\in P$, $p\leq\xi(\zeta(p))$. Finally, 
$\zeta\circ\xi\circ\zeta = \zeta$ and $\xi\circ\zeta\circ\xi = \xi$.
\end{remark}

\begin{definition}\label{DSup}
A \emph{sup-semilattice} is an ordered set $\mathbf{I} = (I,\leq)$ in which any two elements $i$, $j$ have a least upper bound $i\vee j$ (equivalent terminology for this is a \emph{join semilattice} or a $\vee$-\emph{semilattice}). A morphism from a sup-semilattice $\mathbf{I}$ to another sup-semilattice $\mathbf{P}$ is a triple $(\mathbf{I},\xi,\mathbf{P})$, abbreviated to $\xi\colon \mathbf{I}\mor \mathbf{P}$, where $\xi$ is a mapping from $I$ to $P$ such that, for every $i$, $j\in I$, $\xi(i\vee j) = \xi(i)\vee \xi(j)$. We let $\mathsf{Ssl}$ stand for the category of \emph{sup-semilattices} and morphisms between them.
\end{definition}

\begin{remark}
The category $\mathsf{Ssl}$ is isomorphic to the category of idempotent commutative semigroups.
\end{remark}

\begin{remark}
If $\xi\colon \mathbf{I}\mor \mathbf{P}$ is a morphism between sup-semilattices, then $\xi$ is an isotone mapping between them. Moreover, in $\mathsf{Ssl}$, the monomorphisms are the injective morphisms and the epimorphisms the surjective morphisms.
\end{remark}

\begin{proposition}
Let $\mathbf{I}$ and $\mathbf{P}$ be sup-semilattices and $\xi$ a residuated morphism from $\mathbf{I}$ to $\mathbf{P}$. Then, the unique isotone mapping $\zeta$ from $\mathbf{P}$ to $\mathbf{I}$ such that, for every $i\in I$, $\zeta(\xi(i))\leq i$ and , for every $p\in P$, $p\leq\xi(\zeta(p))$ is a morphism from $\mathbf{P}$ to $\mathbf{I}$, i.e., for every $p$, $q\in P$, $\zeta(p\vee q) = \zeta(p)\vee\zeta(q)$.  
\end{proposition}


In what follows, after defining, for a category $\mathsf{I}$, the notion of split $\mathsf{I}$-indexed category, we define, for a split indexed category $(\mathsf{I},F)$, where $F$ is a split $\mathsf{I}$-indexed category, the object part of the Grothendieck construction (see~\cite{Gro71}) at $(\mathsf{I},F)$, which will be a category, denoted by $\int^{\mathsf{I}}F$. Next, after defining the notion of morphism between split indexed categories, we define, for a morphism $(G,\eta)$ from a split indexed category $(\mathsf{I},F)$ to another $(\mathsf{I}',F')$, the morphism part of the Grothendieck construction at $(G,\eta)$, which will be a functor, denoted by $\int^{G}\eta$, from $\int^{\mathsf{I}}F$ to $\int^{\mathsf{I}'}F'$. Finally, we provide a sufficient condition for the functor $\int^{G}\eta$ to have a left adjoint. 

\begin{definition}
Let $\mathsf{I}$ be a category. A \emph{split} $\mathsf{I}$-\emph{indexed category} is a contravariant functor $F$ from $\mathsf{I}$ to $\mathsf{Cat}$, the category of categories and functors. Given an object $i\in \mathrm{Ob}(\mathsf{I})$, we write $\mathsf{F}_{i}$ for the category $F(i)$, and given a morphism $\varphi\in \mathrm{Hom}_{\mathsf{I}}(i,j)$, we write $F_{\varphi}$ for the functor $F(\varphi)\colon \mathsf{F}_{j}\mor \mathsf{F}_{i}$. A \emph{split indexed category} is an ordered pair $(\mathsf{I},F)$ in which $\mathsf{I}$ is a category and $F$ a split $\mathsf{I}$-indexed category.
\end{definition}

\begin{definition}
Let $(\mathsf{I},F)$ be a split indexed category. Then $\int(\mathsf{I},F)$, the \emph{Grothen\-dieck construction} at $(\mathsf{I},F)$, denoted by $\int^{\mathsf{I}}F$, is the category defined as follows:
\begin{enumerate}
\item $\mathrm{Ob}(\int^{\mathsf{I}}F) = \bigcup_{i\in\mathrm{Ob}(\mathsf{I})}(\{i\}\times\mathrm{Ob}(\mathsf{F}_{i}))$.
\item For every $(i,x), (j,y)\in \mathrm{Ob}(\int^{\mathsf{I}}F)$, $\mathrm{Hom}_{\int^{\mathsf{I}}F}((i,x),(j,y))$ is the set of all ordered pairs $(\varphi,f)$, where $\varphi\in \mathrm{Hom}_{\mathsf{I}}(i,j)$ and $f\in \mathrm{Hom}_{\mathsf{F}_{i}}(x,F_{\varphi}(y))$.
\item For every $(i,x)\in \mathrm{Ob}(\int^{\mathsf{I}}F)$, the identity morphism at $(i,x)$ is given by $(\mathrm{id}_{i},\mathrm{id}_{x})$.
\item For every $(i,x), (j,y), (k,z)\in \mathrm{Ob}(\int^{\mathsf{I}}F)$, every $(\varphi,f)\colon (i,x)\mor (j,y)$, and every $(\psi,g)\colon (j,y)\mor (k,z)$, the composite morphism $(\psi,g)\circ (\varphi,f)$ from $(i,x)$ to $(k,z)$ is
     $$
      (\psi,g)\circ (\varphi,f) = (\psi\circ \varphi,F_{\psi}(g)\circ f)\colon (i,x)\mor (k,z).
     $$
     Notice that $f\colon x\mor F_{\varphi}(y)$, $g\colon y\mor F_{\psi}(z)$, hence, taking into account that $F_{\varphi}$ is a functor from $\mathsf{F}_{j}$ to $\mathsf{F}_{i}$, $F_{\varphi}(g)\colon F_{\varphi}(y)\mor F_{\varphi}(F_{\psi}(z))$. Therefore $F_{\psi}(g)\circ f\colon x\mor F_{\varphi}(F_{\psi}(z))$.
\end{enumerate}

We denote by $\pi_{F}$ the canonical split fibration from $\int^{\mathsf{I}}F$ to $\mathsf{I}$.
\end{definition}

\begin{proposition}\cite[Theorem~1, p.~247]{tbg91}\label{PCompl}
Let $\mathsf{I}$ be a category and $F$ a contravariant functor from $\mathsf{I}$ to $\mathsf{Cat}$. If 
$\mathsf{I}$ is complete, for every object $i$ of $\mathsf{I}$, $\mathsf{F}_{i}$ is complete and, for every objects $i$, $j$ of $\mathsf{I}$ and every morphism $\varphi\in \mathrm{Hom}_{\mathsf{I}}(i,j)$, the functor $F_{\varphi}$ from $\mathsf{F}_{j}$ to $\mathsf{F}_{i}$ is continuous, then $\int^{\mathsf{I}}F$ is complete.
\end{proposition}

\begin{proposition}\cite[Theorem~2, p.~250]{tbg91}\label{PCoCompl}
Let $\mathsf{I}$ be a category and $F$ a contravariant functor from $\mathsf{I}$ to $\mathsf{Cat}$. If 
$\mathsf{I}$ is cocomplete, for every object $i$ of $\mathsf{I}$, $\mathsf{F}_{i}$ is cocomplete and, for every objects $i$, $j$ of $\mathsf{I}$ and every morphism $\varphi\in \mathrm{Hom}_{\mathsf{I}}(i,j)$, the functor $F_{\varphi}$ from $\mathsf{F}_{j}$ to $\mathsf{F}_{i}$ has a left adjoint, then $\int^{\mathsf{I}}F$ is cocomplete.
\end{proposition}

\begin{definition}
Let $F$ be a split $\mathsf{I}$-indexed category and $F'$ a split $\mathsf{I}'$-indexed category. A morphism from $(\mathsf{I},F)$ to $(\mathsf{I}',F')$ is an ordered pair $(G,\eta)$, where $G$ is a functor from $\mathsf{I}$ to $\mathsf{I}'$ and $\eta$ a natural transformation from $F$ to $F'\circ G^{\mathrm{op}}$.
\end{definition}

\begin{proposition}
Let $(G,\eta)$ be a morphism from $(\mathsf{I},F)$ to $(\mathsf{I}',F')$. Then there exists a functor $\int^{G}\eta$ from $\int^{\mathsf{I}}F$ to $\int^{\mathsf{I}'}F'$ such that $G\circ \pi_{F} = \pi_{F'}\circ \int^{G}\eta$.
\end{proposition}


\begin{remark}
There exists a functor $\int^{\bullet}$ from the category of split indexed categories $(\mathsf{I},F)$ (and morphisms from $(\mathsf{I},F)$ to $(\mathsf{I}',F')$ the pairs $(G,\eta)$ as above) to the category of split fibrations.
\end{remark}

\begin{proposition}\label{GCLadj}
Let $(G,\eta)$ be a morphism from the split indexed category $(\mathsf{I},F)$ to the split indexed category $(\mathsf{I}',F')$. If (1) $G$ has a left adjoint $T$, (2), for every $i'\in\mathrm{Ob}(\mathsf{I}')$, $F'_{\alpha'_{i'}}\colon \mathsf{F}'_{G(T(i'))}\mor \mathsf{F}'_{i'}$, where $\alpha'_{i'}\colon i'\mor G(T(i'))$ is the value at $i'$ of the unit of the adjunction $T\dashv G$, has a left adjoint $T'_{\alpha'_{i'}}$, and (3), for every $i\in\mathrm{Ob}(\mathsf{I})$, $\eta_{i}\colon \mathsf{F}_{i}\mor \mathsf{F}'_{G(i)}$ has a left adjoint $\xi_{i}$, then the functor $\int^{G}\eta$ from $\int^{\mathsf{I}}F$ to $\int^{\mathsf{I}'}F'$ has a left adjoint.
\end{proposition}

\section{On the strong Lawvere adjoint cylinder between $\int^{\mathsf{Ssl}}\mathrm{Isys_{\Sigma}}$ and $\mathsf{Ssl}$}

In this section, for the signature $\Sigma$, we define the notions of sup-semilattice inductive system of $\Sigma$-algebras relative to a sup-semilattice and of morphism between sup-semilattice inductive systems of $\Sigma$-algebras. Moreover, we define a contravariant functor 
$\mathrm{Isys}_{\Sigma}$ from the category $\mathsf{Ssl}$ to the category $\mathsf{Cat}$ and we investigate the relationships between $\int^{\mathsf{Ssl}}\mathrm{Isys_{\Sigma}}$, the category obtained from $\mathrm{Isys}_{\Sigma}$ by means of the Grothendieck contruction, and 
$\mathsf{Ssl}$.

We begin by defining the notions of sup-semilattice inductive system of $\Sigma$-algebras relative to a sup-semilattice and of morphism between sup-semilattice inductive systems of $\Sigma$-algebras.

\begin{definition}\label{DIndSys}
Let $\mathbf{I}$ be a sup-semilattice. Then a \emph{sup-semilattice inductive system of 
$\Sigma$-algebras relative to} $\mathbf{I}$ or an $\mathbf{I}$-\emph{inductive system of 
$\Sigma$-algebras} for brevity, is an ordered pair 
$
\mathcal{A} = ((\mathbf{A}_{i})_{i\in I},(f_{i,j})_{(i,j)\in \leq}),
$ 
where $(\mathbf{A}_{i})_{i\in I}$ is an $I$-indexed family of $\Sigma$-algebras, and 
$(f_{i,j})_{(i,j)\in \leq}$ a family of homomorphisms in $\prod_{(i,j)\in \leq}\mathrm{Hom}(\mathbf{A}_{i},\mathbf{A}_{j})$ such that, for every $i\in I$, $f_{i,i} = \mathrm{id}_{\mathbf{A}_{i}}$ and, for every 
$i$, $j$, $k\in I$, if $i\leq j$ and $j\leq k$, then $f_{j,k}\circ f_{i,j} = f_{i,k}$, i.e., a covariant functor from the category canonically associated to $\mathbf{I}$ to 
$\mathsf{Alg}(\Sigma)$. Let $\mathcal{A}$ and $\mathcal{B} = ((\mathbf{B}_{i})_{i\in I},(g_{i,j})_{(i,j)\in \leq})$ be $\mathbf{I}$-inductive systems of 
$\Sigma$-algebras. Then a \emph{morphism} from $\mathcal{A}$ to $\mathcal{B}$ is an ordered triple $(\mathcal{A},u,\mathcal{B})$, abbreviated to $u\colon \mathcal{A}\mor \mathcal{B}$, where $u = (u_{i})_{i\in I}$ is a family of homomorphisms in $\prod_{i\in I}\mathrm{Hom}(\mathbf{A}_{i},\mathbf{B}_{i})$ such that, for every $(i,j)\in \leq$, $u_{j}\circ f_{i,j} = g_{i,j}\circ u_{i}$, i.e., a natural transformation from $\mathcal{A}$ to $\mathcal{B}$. We denote by $\mathsf{Alg}(\Sigma)^{\mathbf{I}}$ the corresponding category. 
\end{definition}

We next define the assignment $\mathrm{Isys_{\Sigma}}$ from $\mathsf{Ssl}$ to $\mathsf{Cat}$.

\begin{definition}\label{DIsys}
Let $\Sigma$ be a signature. Then we let $\mathrm{Isys_{\Sigma}}$ stand for the assignment from $\mathsf{Ssl}$ to $\mathsf{Cat}$ defined as follows:
\begin{enumerate}
\item for every sup-semilattice $\mathbf{I}$, $\mathrm{Isys_{\Sigma}}(\mathbf{I})$ is $\mathsf{Alg}(\Sigma)^{\mathbf{I}}$, and 
\item for every morphism $\xi\colon\mathbf{I}\mor \mathbf{P}$, $\mathrm{Isys_{\Sigma}}(\xi)$ is the functor  from $\mathsf{Alg}(\Sigma)^{\mathbf{P}}$ to $\mathsf{Alg}(\Sigma)^{\mathbf{I}}$ defined as follows:
\begin{enumerate}
\item for every $\mathbf{P}$-inductive systems of $\Sigma$-algebras $\mathcal{A} = ((\mathbf{A}_{p})_{p\in P},(f_{p,q})_{(p,q)\in \leq})$, $\mathrm{Isys}_{\Sigma}(\xi)(\mathcal{A})$, denoted by $\mathcal{A}_{\xi}$ for short,  is $((\mathbf{A}_{\xi(i)})_{i\in I},(f_{\xi(i),\xi(j)})_{(i,j)\in \leq})$, which is an $\mathbf{I}$-inductive systems of $\Sigma$-algebras, and
\item for every morphism of $\mathbf{P}$-inductive systems of $\Sigma$-algebras $u = (u_{p})_{p\in P}$ from $\mathcal{A}$ to $\mathcal{B} = ((\mathbf{B}_{p})_{p\in P},(g_{p,q})_{(p,q)\in \leq})$, 
$\mathrm{Isys}_{\Sigma}(\xi)(u)$, denoted by $u_{\xi}$ for short, is $(u_{\xi(i)})_{i\in I}$, which is a morphism of $\mathbf{I}$-inductive systems of $\Sigma$-algebras from $\mathcal{A}_{\xi}$ to $\mathcal{B}_{\xi}$.
\end{enumerate} 
\end{enumerate}
\end{definition}

\begin{proposition}
Let $\Sigma$ be a signature. Then the assignment $\mathrm{Isys_{\Sigma}}$ from $\mathsf{Ssl}$ to $\mathsf{Cat}$ is a contravariant functor.
\end{proposition}

\begin{definition}\label{DIntSsl}
From the split indexed category $(\mathsf{Ssl},\mathrm{Isys_{\Sigma}})$ we obtain, by means of the Grothendieck construction, the category $\int^{\mathsf{Ssl}}\mathrm{Isys_{\Sigma}}$ which has as objects the \emph{sup-semilattice inductive system of $\Sigma$-algebras} or, simply, the  \emph{inductive system of $\Sigma$-algebras}, i.e., the ordered pairs 
$
\mathbcal{A} = (\mathbf{I},\mathcal{A}),
$
where $\mathbf{I}$ is a sup-semilattice and $\mathcal{A}$ an $\mathbf{I}$-inductive system of 
$\Sigma$-algebras; and as morphisms from $\mathbcal{A}$ to $\mathbcal{B} = (\mathbf{P},\mathcal{B})$ the ordered triples $(\mathbcal{A},(\xi,u),\mathbcal{B})$, abbreviated to $(\xi,u)\colon \mathbcal{A}\mor \mathbcal{B}$, where $\xi$ is a morphism from $\mathbf{I}$ to $\mathbf{P}$ and $u$ a morphism from the $\mathbf{I}$-inductive system of $\Sigma$-algebras $\mathcal{A}$ to the $\mathbf{I}$-inductive system of $\Sigma$-algebras $\mathcal{B}_{\xi}$. Moreover, we let 
$\pi_{\mathrm{Isys_{\Sigma}}}$ stand for the canonical split fibration, i.e, the functor from $\int^{\mathsf{Ssl}}\mathrm{Isys_{\Sigma}}$ to $\mathsf{Ssl}$ that sends $\mathbcal{A} = (\mathbf{I},\mathcal{A})$ to $\mathbf{I}$ and $(\xi,u)\colon \mathbcal{A}\mor \mathbcal{B}$, with $\mathbcal{B} = (\mathbf{P},\mathcal{B})$, to $\xi\colon \mathbf{I}\mor \mathbf{P}$.
\end{definition}

We next show that from $\mathsf{Ssl}$ to $\int^{\mathsf{Ssl}}\mathrm{Isys}_{\Sigma}$ there are functors $L_{\Sigma}$ and $K_{\Sigma}$ such that $L_{\Sigma}\dashv \pi_{\mathrm{Isys_{\Sigma}}}\dashv K_{\Sigma}$ and that from $L_{\Sigma}$ to $K_{\Sigma}$ there exists a natural transformation, that we will denote by $\gamma_{\Sigma}$.

\begin{proposition}\label{PLLeftAdj}
There exists a full and faithful functor $L_{\Sigma}$ from $\mathsf{Ssl}$ to $\int^{\mathsf{Ssl}}\mathrm{Isys}_{\Sigma}$ which is left adjoint and right inverse of $\pi_{\mathrm{Isys_{\Sigma}}}$.
\end{proposition}

\begin{proof}
Let $L_{\Sigma}$ be the functor from $\mathsf{Ssl}$ to $\int^{\mathsf{Ssl}}\mathrm{Isys}_{\Sigma}$ defined as follows:
\begin{enumerate}
\item for every sup-semilattice $\mathbf{I}$, $L_{\Sigma}(\mathbf{I})$ is the inductive system of $\Sigma$-algebras $(\mathbf{I},((\mathbf{T}_{\Sigma}(\varnothing))_{i\in I},(\mathrm{id}_{\mathbf{T}_{\Sigma}(\varnothing)})_{(i,j)\in \leq}))$, and 
\item for every morphism $\xi$ from $\mathbf{I}$ to $\mathbf{P}$, $L_{\Sigma}(\xi)$ is $(\xi,(\mathrm{id}_{\mathbf{T}_{\Sigma}(\varnothing)})_{i\in I})$, which is a morphism from $L_{\Sigma}(\mathbf{I})$ to $L_{\Sigma}(\mathbf{P})$ (because $(\mathrm{id}_{\mathbf{T}_{\Sigma}(\varnothing)})_{i\in I}$ is a morphism from $((\mathbf{T}_{\Sigma}(\varnothing))_{i\in I},(\mathrm{id}_{\mathbf{T}_{\Sigma}(\varnothing)})_{(i,j)\in \leq})$ to $((\mathbf{T}_{\Sigma}(\varnothing))_{p\in P},(\mathrm{id}_{\mathbf{T}_{\Sigma}(\varnothing)})_{(p,q)\in \leq})^{\xi}$.
\end{enumerate} 
From the fact that $\mathbf{T}_{\Sigma}(\varnothing)$ is initial in $\mathsf{Alg}(\Sigma)$, it follows that, for every sup-semilattice $\mathbf{I}$ there exists a universal arrow from $\mathbf{I}$ to $L_{\Sigma}$. 

It is obvious that $L_{\Sigma}$ is full, faithful and right inverse of 
$\pi_{\mathrm{Isys_{\Sigma}}}$.
\end{proof}

\begin{remark}
If $\Sigma_{0} = \varnothing$, then, in Proposition~\ref{PLLeftAdj}, we have that $\mathbf{T}_{\Sigma}(\varnothing)$ is $\pmb{\varnothing}$, the initial $\Sigma$-algebra.
\end{remark}

\begin{proposition}
There exists a full and faithful functor $K_{\Sigma}$ from $\mathsf{Ssl}$ to $\int^{\mathsf{Ssl}}\mathrm{Isys}_{\Sigma}$ which is right adjoint and right inverse of $\pi_{\mathrm{Isys_{\Sigma}}}$.
\end{proposition}

\begin{proof}
Let $K_{\Sigma}$ be the functor from $\mathsf{Ssl}$ to $\int^{\mathsf{Ssl}}\mathrm{Isys}_{\Sigma}$ defined as follows:
\begin{enumerate}
\item for every sup-semilattice $\mathbf{I}$, $K_{\Sigma}(\mathbf{I})$ is the inductive system of $\Sigma$-algebras $(\mathbf{I},((\mathbf{K}_{\{i\}})_{i\in I},(t_{\mathbf{K}_{\{i\}},\mathbf{K}_{\{j\}}})_{(i,j)\in\leq}))$, where, for every $i\in I$, the underlying set of 
$\mathbf{K}_{\{i\}}$ is $\{i\}$, i.e., a single-sorted set,  (observe that the $\Sigma$-algebra $\mathbf{K}_{\{i\}}$ is isomorphic to the final 
$\Sigma$-algebra $\mathbf{1}$), and, for every $(i,j)\in \leq$, $t_{\mathbf{K}_{\{i\}},\mathbf{K}_{\{j\}}}$ is the unique homomorphism from $\mathbf{K}_{\{i\}}$ to $\mathbf{K}_{\{j\}}$, and
\item for every morphism $\xi$ from $\mathbf{I}$ to $\mathbf{P}$, $K_{\Sigma}(\xi)$ is 
$(\xi,(t_{\mathbf{K}_{\{i\}},\mathbf{K}_{\{\xi(i)\}}})_{i\in I})$, which is a morphism from $K_{\Sigma}(\mathbf{I})$ to $K_{\Sigma}(\mathbf{P})$ (because $(t_{\mathbf{K}_{\{i\}},\mathbf{K}_{\{\xi(i)\}}})_{i\in I}$ is a morphism from $((\mathbf{K}_{\{i\}})_{i\in I},(t_{\mathbf{K}_{\{i\}},\mathbf{K}_{\{j\}}})_{(i,j)\in\leq})$ to 
$((\mathbf{K}_{\{p\}})_{p\in P},(t_{\mathbf{K}_{\{p\}},\mathbf{K}_{\{q\}}})_{(p,q)\in\leq})_{\xi}$).
\end{enumerate}
From the fact that $\mathbf{K}_{\{i\}}$ is final in $\mathsf{Alg}(\Sigma)$, it follows that, for every sup-semilattice $\mathbf{I}$ there exists a universal arrow from $K_{\Sigma}$ to $\mathbf{I}$.

It is obvious that $K_{\Sigma}$ is full, faithful and right inverse of 
$\pi_{\mathrm{Isys_{\Sigma}}}$.  
\end{proof}

\begin{corollary}
The ordered triple $(L_{\Sigma},\pi_{\mathrm{Isys_{\Sigma}}},K_{\Sigma})$ is a strong Lawvere adjoint cylinder (see~\cite{L94}) (``strong'' because 
$\pi_{\mathrm{Isys}_{\Sigma}}\circ L_{\Sigma} = \mathrm{Id}_{\mathsf{Ssl}}$ and 
$\pi_{\mathrm{Isys}_{\Sigma}}\circ K_{\Sigma} = \mathrm{Id}_{\mathsf{Ssl}}$).
\end{corollary}

\begin{proposition}
There exists a natural transformation $\gamma_{\Sigma}$ from $L_{\Sigma}$ to $K_{\Sigma}$.
\end{proposition}

\begin{proof}
It suffices to take, for every sup-semilattice $\mathbf{I}$, as $\gamma_{\Sigma,\mathbf{I}}$, the component of 
$\gamma_{\Sigma}$ at $\mathbf{I}$, precisely  
$(\mathrm{id}_{\mathbf{I}},(t_{\mathbf{T}_{\Sigma}(\varnothing),\mathbf{K}_{\{i\}}})_{i\in I})$, which is a morphism from $L_{\Sigma}(\mathbf{I})$ to $K_{\Sigma}(\mathbf{I})$.
\end{proof}



In the following proposition we state that the category $\int^{\mathsf{Ssl}}\mathrm{Isys}_{\Sigma}$ is complete.

\begin{proposition}\label{PIntCompl} The category $\int^{\mathsf{Ssl}}\mathrm{Isys}_{\Sigma}$ is complete.
\end{proposition}
\begin{proof}
We make use of Proposition~\ref{PCompl}. Since $\mathsf{Ssl}$ is complete, for every sup-semilattice $\mathbf{I}$, $\mathsf{Alg}(\Sigma)^{\mathbf{I}}$ is complete, and, for every $\xi\colon \mathbf{I}\mor \mathbf{P}$, the functor $\mathrm{Isys}_{\Sigma}(\xi)$ from $\mathsf{Alg}(\Sigma)^{\mathbf{P}}$ to $\mathsf{Alg}(\Sigma)^{\mathbf{I}}$ is continuous, we have that $\int^{\mathsf{Ssl}}\mathrm{Isys}_{\Sigma}$ is complete.
\end{proof}

If instead of $\mathsf{Ssl}$ we consider the wide subcategory $\mathsf{Ssl}_{\mathrm{rd}}$ of $\mathsf{Ssl}$ of sup-semilattices and residuated morphisms, then we have that the category $\int^{\mathsf{Ssl}_{\mathrm{rd}}}\mathrm{Isys}_{\Sigma}$ is cocomplete. 

\begin{proposition}\label{PIntCoCompl} For the category $\mathsf{Ssl}_{\mathrm{rd}}$, the category $\int^{\mathsf{Ssl}_{\mathrm{rd}}}\mathrm{Isys}_{\Sigma}$ is cocomplete.
\end{proposition}
\begin{proof}
We make use of Proposition~\ref{PCoCompl}. Since $\mathsf{Ssl}$ is cocomplete, for every sup-semilattice $\mathbf{I}$, $\mathsf{Alg}(\Sigma)^{\mathbf{I}}$ is cocomplete, and, for every residuated $\xi\colon \mathbf{I}\mor \mathbf{P}$, the functor $\mathrm{Isys}_{\Sigma}(\xi)$ from $\mathsf{Alg}(\Sigma)^{\mathbf{P}}$ to $\mathsf{Alg}(\Sigma)^{\mathbf{I}}$ has an adjoint, we have that the subcategory of $\int^{\mathsf{Ssl}}\mathrm{Isys}_{\Sigma}$ whose morphisms are those $(\xi,u)$ for which $\xi$ is residuated, is cocomplete. If $\zeta$ is the left adjoint of $\xi$, then 
$\mathrm{Hom}_{\mathsf{Alg}(\Sigma)^{\mathbf{I}}}(\mathcal{B}_{\xi},\mathcal{A})$ and $\mathrm{Hom}_{\mathsf{Alg}(\Sigma)^{\mathbf{P}}}(\mathcal{B},\mathcal{A}_{\zeta})$ are naturally isomorphic, i.e., $\mathrm{Isys}_{\Sigma}(\xi)\dashv \mathrm{Isys}_{\Sigma}(\zeta)$.
Let $\mathcal{A} = ((\mathbf{A}_{i})_{i\in I},(f_{i,j})_{(i,j)\in \leq})$ be an $\mathbf{I}$-inductive system of $\Sigma$-algebras, $\mathcal{B} = ((\mathbf{B}_{p})_{p\in P},(g_{p,q})_{(p,q)\in \leq})$ 
a $\mathbf{P}$-inductive system of $\Sigma$-algebras, $u\colon \mathcal{B}_{\xi}\mor \mathcal{A}$ 
and $v\colon \mathcal{B}\mor \mathcal{A}_{\zeta}$. Then 
$\varphi(u) = (u_{\zeta(p)}\circ g_{p,\xi(\zeta(p))})_{p\in P}$, where, for every $p\in P$, we have that 
$$
u_{\zeta(p)}\circ g_{p,\xi(\zeta(p))}\colon\mathbf{B}_{p}\mor\mathbf{B}_{\xi(\zeta(p))}\mor\mathbf{A}_{\zeta(p)},
$$
is a morphism from $\mathcal{B}$ to $\mathcal{A}_{\zeta}$, and 
$\psi(v) = (f_{\zeta(\xi(i)),i}\circ v_{\xi(i)})_{i\in I}$, where, for every $i\in I$, we have that  
$$
f_{\zeta(\xi(i)),i}\circ v_{\xi(i)}\colon\mathbf{B}_{\xi(i)}\mor\mathbf{A}_{\zeta(\xi(i))}\mor\mathbf{A}_{i}
$$
is a morphism from $\mathcal{B}_{\xi}$ to $\mathcal{A}$. It remains to be shown that both applications are inverse to each other and natural. We leave it to the reader to verify this last statement, who should take into account that $\zeta\circ\xi\circ\zeta = \zeta$ and $\xi\circ\zeta\circ\xi = \xi$.
\end{proof}

\begin{remark}
For the category $\int^{\mathsf{Ssl}}\mathrm{Isys_{\Sigma^{\neq 0}}}$, obtained from the contravariant functor $\mathrm{Isys_{\Sigma^{\neq 0}}}$ (defined in the same way as 
$\mathrm{Isys_{\Sigma}}$), and suitably defined functors $L_{\Sigma^{\neq 0}}$, $\pi_{\mathrm{Isys_{\Sigma^{\neq 0}}}}$ and $K_{\Sigma^{\neq 0}}$ we also have that there exists a natural transformation $\gamma_{\Sigma^{\neq 0}}$ from $L_{\Sigma^{\neq 0}}$ to $K_{\Sigma^{\neq 0}}$ and that the ordered triple 
$(L_{\Sigma^{\neq 0}},\pi_{\mathrm{Isys_{\Sigma^{\neq 0}}}},K_{\Sigma^{\neq 0}})$ is a strong Lawvere adjoint cylinder. Moreover, Propositions~\ref{PIntCompl} and \ref{PIntCoCompl} also hold for $\int^{\mathsf{Ssl}}\mathrm{Isys_{\Sigma^{\neq 0}}}$.
\end{remark}

\section{On the relationship between $\mathsf{Ssl}$ and $\mathsf{Alg}(\Sigma^{\neq 0})$}

In this section we prove that from $\mathsf{Ssl}$ to $\mathsf{Alg}(\Sigma^{\neq 0})$ there exists a functor $W_{\Sigma^{\neq 0}}$ which has a left adjoint.

\begin{definition}\label{DSpp} Let $\Sigma$ be a signature.
Then we let $W_{\Sigma^{\neq 0}}$ stand for the assignment from $\mathsf{Ssl}$ to $\mathsf{Alg}(\Sigma^{\neq 0})$ defined as follows:
\begin{enumerate}
\item for every sup-semilattice $\mathbf{I}$, $W_{\Sigma^{\neq 0}}(\mathbf{I})$ is the 
$\Sigma^{\neq 0}$-algebra $\mathbf{W}_{\Sigma^{\neq 0}}(\mathbf{I}) = (I,F^{\mathbf{I}})$, where, for every $n\in\mathbb{N}-1$ and every $\sigma\in \Sigma_{n}$, $F^{\mathbf{I}}_{\sigma}$ is the mapping from $I^{n}$ to $I$ that sends 
$(i_{j})_{j\in n}\in I^{n}$ to $\bigvee_{j\in n}i_{j}\in I$, and
\item for every morphism $\xi$ from $\mathbf{I}$ to $\mathbf{P}$, $W_{\Sigma^{\neq 0}}(\xi)$ is 
$(\mathbf{W}_{\Sigma^{\neq 0}}(\mathbf{I}),\xi,\mathbf{W}_{\Sigma^{\neq 0}}(\mathbf{P}))$, abbreviated to $\xi\colon \mathbf{W}_{\Sigma^{\neq 0}}(\mathbf{I})\mor \mathbf{W}_{\Sigma^{\neq 0}}(\mathbf{P})$. Note that $W_{\Sigma^{\neq 0}}(\xi)$, identified to $\xi$, is a $\Sigma^{\neq 0}$-homomorphism from $\mathbf{W}_{\Sigma^{\neq 0}}(\mathbf{I})$ to $\mathbf{W}_{\Sigma^{\neq 0}}(\mathbf{P})$, because, for every $n\in\mathbb{N}-1$, every $\sigma\in \Sigma_{n}$ and every $(i_{j})_{j\in n}\in I^{n}$, we have that
$$
\textstyle
\xi(F^{\mathbf{I}}_{\sigma}((i_{j})_{j\in n})) = 
\xi(\bigvee_{j\in n}i_{j}\in I) =
\bigvee_{j\in n}\xi(i_{j}) =
F^{\mathbf{P}}_{\sigma}((\xi(i_{j}))_{j\in n}).
$$
\end{enumerate}
\end{definition}

\begin{proposition} Let $\Sigma$ be a signature. Then the assignment $W_{\Sigma^{\neq 0}}$ is a functor from $\mathsf{Ssl}$ to $\mathsf{Alg}(\Sigma^{\neq 0})$.
\end{proposition}

We next prove that the functor $W_{\Sigma^{\neq 0}}$ has a left adjoint.

\begin{proposition}\label{PUnivSsl} 
The functor $W_{\Sigma^{\neq 0}}$ from $\mathsf{Ssl}$ to $\mathsf{Alg}(\Sigma^{\neq 0})$ has a left adjoint.
\end{proposition}

\begin{proof}
Let $\mathbf{A}$ be a $\Sigma^{\neq 0}$-algebra. Then, from $A$, the underlying set of $\mathbf{A}$, we obtain $\mathbf{T}_{\vee}(A)$, the free sup-semilattice on $A$, and $\eta_{A}$, the canonical embedding of $A$ into $\mathrm{T}_{\vee}(A)$. Now, $\eta_{A}$ is not necessarily a homomorphism from $\mathbf{A}$ to $W_{\Sigma^{\neq 0}}(\mathbf{T}_{\vee}(A))$ as it may be the case that for an $n\in \mathbb{N}-1$, a $\sigma\in\Sigma_{n}$ and a family $(x_{i})_{i\in n}\in A^{n}$, $\eta_{A}(F_{\sigma}((x_{i})_{i\in n}))$ and $\bigvee_{i\in n}\eta_{A}(x_{i})$ be different. Let $\Psi^{\mathbf{A}}$ be the least congruence on the 
sup-semilattice $\mathbf{T}_{\vee}(A)$ containing all ordered pairs of the form 
$$
\textstyle
(\eta_{A}(F_{\sigma}((x_{i})_{i\in n})),\bigvee_{i\in n}\eta_{A}(x_{i})),
$$ 
for $n\in \mathbb{N}-1$, $\sigma\in\Sigma_{n}$ and $(x_{i})_{i\in n}\in A^{n}$. Then $\mathrm{pr}_{\Psi^{\mathbf{A}}}\circ\eta_{A}$, denoted by $\eta_{\mathbf{A}}$ for short, where $\mathrm{pr}_{\Psi^{\mathbf{A}}}$ is the canonical projection from $\mathbf{T}_{\vee}(A)$ to $\mathbf{T}_{\vee}(A)/\Psi^{\mathbf{A}}$, is a homomorphism from $\mathbf{A}$ to $W_{\Sigma^{\neq 0}}(\mathbf{T}_{\vee}(A)/\Psi^{\mathbf{A}})$. 

We now verify that the pair $(\mathbf{T}_{\vee}(A)/\Psi^{\mathbf{A}},\eta_{\mathbf{A}})$ is a universal morphism from $\mathbf{A}$ to $W_{\Sigma^{\neq 0}}$. Let $\mathbf{I}$ be a sup-semilattice and $f$ a homomorphism from $\mathbf{A}$ to $W_{\Sigma^{\neq 0}}(\mathbf{I})$. Then, from $f$, we obtain $f^{\sharp}$, the unique morphism from $\mathbf{T}_{\vee}(A)$ to $\mathbf{I}$ such that $f^{\sharp}\circ \eta_{A}  = f$. Moreover, since $\Psi^{\mathbf{A}}\subseteq \mathrm{Ker}(f^{\sharp})$, there exists a unique morphism $f^{\flat}$ from $\mathbf{T}_{\vee}(A)/\Psi^{\mathbf{A}}$ to $\mathbf{I}$ such that 
$f^{\flat}\circ \mathrm{pr}_{\Psi^{\mathbf{A}}} = f^{\sharp}$. Then, by virtue of the definitions of the morphisms involved, it follows that  
$$
f^{\flat}\circ \eta_{\mathbf{A}} = f^{\flat}\circ \mathrm{pr}_{\Psi^{\mathbf{A}}}\circ \eta_{A} = f^{\sharp}\circ \eta_{A} = f.
$$
It remains to be shown that $f^{\flat}$ is the unique morphism from $\mathbf{T}_{\vee}(A)/\Psi^{\mathbf{A}}$ to $\mathbf{I}$ such that the above condition is fulfilled. We leave it to the reader to verify this last assertion.

%
%
%
%
%

We denote by $M_{\Sigma^{\neq 0}}$ the left adjoint of $W_{\Sigma^{\neq 0}}$.
\end{proof}


\section{On the category $\mathsf{Lnb}$ of left normal and its relationship with the category $\mathsf{Ssl}$}

In this section we define the category $\mathsf{Lnb}$ of left normal bands, and, for a left normal band, we  define by recursion the family of the iterates of its structural operation and, for later use, establish some useful properties of such a family. Moreover, we prove that the inclusion functor $\mathrm{In}_{\mathsf{Ssl},\mathsf{Lnb}}$ from $\mathsf{Ssl}$ to $\mathsf{Lnb}$ has a left adjoint. 


We begin by defining the category $\mathsf{Lnb}$ of left normal bands.

\begin{definition}\label{DDO} 
A \emph{left normal band} is an ordered pair $(A,D)$ in which $A$ is a set and $D$ a mapping $D\colon A\times A \mor A$ such that, for every every $x$, $y$, $z\in A$, the following conditions are satisfied: 
\begin{align*}
D(x,x)&=x;
\tag{D1}\label{D1}
\\
D(x,D(y,z))&=D(D(x,y),z);
\tag{D2}\label{D2}
\\
D(x,D(y,z))&=D(x,D(z,y)).
\tag{D3}\label{D3}
\end{align*}
Let $(A,D)$ and $(B,E)$ be left normal bands. A morphism $h$ from $(A,D)$ to $(B,E)$, written $h\colon (A,D)\mor (B,E)$, is a  mapping $h\colon A \mor B$ such that 
$ E\circ (h\times h)=h\circ D$.
%
%

We let $\mathsf{Lnb}$ stand for the corresponding category. 
\end{definition}

We next define, for a left normal band $(A,D)$, the family of the iterates on the left and on the right of the structural operation $D$ of $(A,D)$.

\begin{definition}\label{DIt} Let $(A,D)$ be a left normal band. Then we define the family of the \emph{iterates of $D$ on the right}, $(D^{\mathrm{R}}_{n})_{n\geq 1}$, recursively, as follows
\[
\begin{array}{lcl}
D^{\mathrm{R}}_{1}&=&\mathrm{id}_{A};\\
D^{\mathrm{R}}_{n+1}&=&D\circ (\mathrm{id}_{A}\times D^{\mathrm{R}}_{n}),\quad n\geq 1.
\end{array}
\]
Let us note that, for every $n\in\mathbb{N}-1$, $D^{\mathrm{R}}_{n}\colon A^{n}\mor A$. In particular for $n=2$, $D^{\mathrm{R}}_{2}=D$. 

We next define the family of the \emph{iterates of $D$ on the left}, $(D^{\mathrm{L}}_{n})_{n\geq 1}$, recursively, as follows
\[
\begin{array}{lcl}
D^{\mathrm{L}}_{1}&=&\mathrm{id}_{A};\\
D^{\mathrm{L}}_{n+1}&=&D\circ ( D^{\mathrm{L}}_{n}\times \mathrm{id}_{A}),\quad n\geq 1.
\end{array}
\]
Let us note that, for every $n\in\mathbb{N}-1$, $D^{\mathrm{L}}_{n}\colon A^{n}\mor A$. In particular for $n=2$, $D^{\mathrm{L}}_{2}=D$. 
 \end{definition}

We next set out a series of technical lemmas about the families of iterates that will be useful later on.

\begin{lemma}\label{LTech} 
Let $(A,D)$ be a left normal band, $n\in\mathbb{N}-1$, $k\in n$ and $(x_{j})_{j\in n}\in A^{n}$. Then
\begin{align*}
D(D^{\mathrm{R}}_{n}((x_{j})_{j\in n}), x_{k})&=D^{\mathrm{R}}_{n}((x_{j})_{j\in n});
\\
D(x_{k}, D^{\mathrm{L}}_{n}((x_{j})_{j\in n}))&=D^{\mathrm{L}}_{n}((x_{j})_{j\in n}).
\end{align*}
\end{lemma}

\begin{lemma}\label{LTech2} 
Let $(A,D)$ be a left normal band, $n\in\mathbb{N}-1$, $(x_{j})_{j\in n}\in A^{n}$, $m\in\mathbb{N}-1$, $\varphi\colon m\mor n$ and $(x_{\varphi(k)})_{k\in m}$ the family $(x_{j})_{j\in n}\circ\varphi\in A^{m}$. Then
\begin{align*}
D(D^{\mathrm{R}}_{n}((x_{j})_{j\in n}), 
D^{\mathrm{R}}_{m}((x_{\varphi(k)})_{k\in m})
)&=D^{\mathrm{R}}_{n}((x_{j})_{j\in n});
\\
D(D^{\mathrm{R}}_{n}((x_{j})_{j\in n}), 
D^{\mathrm{L}}_{m}((x_{\varphi(k)})_{k\in m})
)&=D^{\mathrm{R}}_{n}((x_{j})_{j\in n});
\\
D(
D^{\mathrm{R}}_{m}((x_{\varphi(k)})_{k\in m}),
D^{\mathrm{L}}_{n}((x_{j})_{j\in n})
)&=D^{\mathrm{L}}_{n}((x_{j})_{j\in n});
\\
D(
D^{\mathrm{L}}_{m}((x_{\varphi(k)})_{k\in m}),
D^{\mathrm{L}}_{n}((x_{j})_{j\in n})
)&=D^{\mathrm{L}}_{n}((x_{j})_{j\in n}).
\end{align*}
\end{lemma}

\begin{remark}
Lemma~\ref{LTech} is a corollary of Lemma~\ref{LTech2}, but we have preferred to keep it as it stands because the equations occurring in it will be used several times in what follows. 
\end{remark}

\begin{corollary}\label{CTech} 
Let $(A,D)$ be a left normal band, $n\in\mathbb{N}-1$ and $(x_{j})_{j\in n}\in A^{n}$. Then
\[
D^{\mathrm{R}}_{n}((x_{j})_{j\in n})=D^{\mathrm{L}}_{n}((x_{j})_{j\in n}).
\]
\end{corollary}

%

\begin{remark} 
Corollary~\ref{CTech} makes it legitimate to use $D_{n}$ instead of $D^{\mathrm{R}}_{n}$ or $D^{\mathrm{L}}_{n}$.
\end{remark}

We next prove that there exists an adjunction between the categories $\mathsf{Ssl}$ and $\mathsf{Lnb}$. We begin by stating that every sup-semilattice is a left-normal band.

\begin{proposition}\label{PSslLNB} 
Let $\mathbf{I}$ be a sup-semilattice, then $\mathbf{I}$ is a left normal band. Moreover, if $f\colon \mathbf{I}\mor \mathbf{P}$ is a morphism of sup-semilattices, then $f\colon\mathbf{I}\mor \mathbf{P}$ is a morphism of left normal bands. Therefore, $\mathsf{Ssl}$ is a subcategory of $\mathsf{Lnb}$.
\end{proposition}

%

\begin{definition}\label{DInc} We let $\mathrm{In}_{\mathsf{Ssl},\mathsf{Lnb}}$ stand for the inclusion functor of $\mathsf{Ssl}$ into $\mathsf{Lnb}$.
\end{definition}

We next associate to the structural operation of a left normal band a binary relation on its underlying set and investigate its properties.

\begin{definition}\label{DDORel} Let $(A,D)$ be a left normal band. Then we let $\Phi^{D}$ stand for the binary relation on $A$ defined as follows: 
\[
\Phi^{D}=\mathrm{Eq}(D, \mathrm{pr}_{0})\cap \mathrm{Eq}(D\circ \mathrm{tw}_{A}, \mathrm{pr}_{1}),
\]
where $\mathrm{tw}_{A}$ is the mapping from $A\times A$ to $A\times A$ that sends $(x,y)$ to $(y,x)$, $\mathrm{pr}_{0}$ and $\mathrm{pr}_{1}$ the canonical projections from $A\times A$ to $A$, $\mathrm{Eq}(D, \mathrm{pr}_{0})$ the equalizer of $D$ and $\mathrm{pr}_{0}$ and $\mathrm{Eq}(D\circ \mathrm{tw}_{A}, \mathrm{pr}_{1})$ the equalizer of $D\circ \mathrm{tw}_{A}$ and $\mathrm{pr}_{1}$.

Let us note that, for every every $x,y\in A$ we have that $(x,y)\in \Phi^{D}$ if, and only if, $D(x,y)=x$ and $D(y,x)=y$.

We call $\Phi^{D}$ the \emph{relation induced by $D$ on $A$}.
\end{definition}

\begin{remark}
Let us note that $\Phi^{D}$ is the pull-back of the canonical inclusions of $\mathrm{Eq}(D, \mathrm{pr}_{0})$ and $\mathrm{Eq}(D\circ \mathrm{tw}_{A}, \mathrm{pr}_{1})$ into $A\times A$.
\end{remark}

%
%
%
%
%
%

We next set out a series of technical results that will be useful later on. In the first of these we state that, for a left normal band $(A,D)$, the relation $\Phi^{D}$ induced by $D$ on $A$ is a congruence on $(A,D)$.

\begin{proposition}\label{PDOA}  Let $(A,D)$ be a left normal band. Then $\Phi^{D}$ is a congruence on $(A,D)$. 
\end{proposition}

\begin{proof} The proof that $\Phi^{D}$ is an equivalence relation on $A$ is straightforward, so we leave it to the reader. Therefore, it remains to prove that, for every $x,x',y,y'\in A$, if $(x,y)$ and $(x',y')$ are pairs in $\Phi^{D}$, then 
\[
(D(x,x'),D(y,y'))\in \Phi^{D}.
\]

Let  $x,x'y,y'$ be elements in $A$ satisfying that $(x,y), (x',y')\in \Phi^{D}$. Hence, $D(x,y)=x$, $D(y,x)=y$, $D(x',y')=x'$ and $D(y',x')=y'$.

We now prove that $(D(x,x'),D(y,y'))\in \Phi^{D}$. For this note that the following chain of equalities holds
\allowdisplaybreaks
\begin{align*}
D(D(x,x'),D(y,y'))&=
D(D(x,x'),D(y',y))
\tag{by \ref{D3}}
\\&=
D(D(D(x,x'),y'),y)
\tag{by \ref{D2}}
\\&=
D(D(x, D(x',y')),y)
\tag{by \ref{D2}}
\\&=
D(D(x, x'),y)
\tag{$D(x',y')=x'$}
\\&=
D(x,D(x',y))
\tag{by \ref{D2}}
\\&=
D(x,D(y,x'))
\tag{by \ref{D3}}
\\&=
D(D(x,y),x')
\tag{by \ref{D2}}
\\&=
D(x,x').
\tag{$D(x,y)=x$}
\end{align*}

By a similar argument we obtain that $D(D(y,y'),D(x,x')) = D(y,y')$.


Thus $(D(x,x'),D(y,y'))\in \Phi^{D}$. 

We call $\Phi^{D}$ \emph{the congruence induced by $D$ on $(A,D)$}.
\end{proof}

Since, for a left normal band $(A,D)$, it is fulfilled that, for every $n\in \mathbb{N}-1$, $D_{n}$ is a derived operation of the structural operation $D$ of $(A,D)$, we have the following corollary. 

\begin{corollary}\label{PDOB} 
Let $(A,D)$ be a left normal band, $\Phi^{D}$ the congruence induced by $D$ on $(A,D)$, $n\in \mathbb{N}-1$ and $(x_{j})_{j\in n}$, $(y_{j})_{j\in n}\in A^{n}$ such that, for every $j\in n$, $(x_{j},y_{j})\in \Phi^{D}$. Then 
\begin{align*}
\left(
D_{n}\left(
(x_{j})_{j\in n}
\right),
D_{n}\left(
(y_{j})_{j\in n}
\right)
\right)&\in \Phi^{D}
.
\end{align*}
\end{corollary}

%
%
%
%
%

\begin{proposition}\label{PDOC}  
Let $(A,D)$ be a left normal band and $x$, $y\in A$. Then
\[
(D(x,y),D(y,x))\in\Phi^{D}.
\]
\end{proposition}

%
%
%


We next prove that the functor $\mathrm{In}_{\mathsf{Ssl},\mathsf{Lnb}}$ has a left adjoint.

\begin{proposition}\label{PUnivLnb} 
The functor $\mathrm{In}_{\mathsf{Ssl},\mathsf{Lnb}}$ from $\mathsf{Ssl}$ to $\mathsf{Lnb}$ has a left adjoint.
\end{proposition}

\begin{proof}
Let $(A,D)$ be a left normal band. We will prove that there exists a universal morphism from $(A,D)$ to $\mathrm{In}_{\mathsf{Ssl},\mathsf{Lnb}}$. To do so, we start by proving that there exists a binary relation $\leq^{D}$ on $A/{\Phi^{D}}$ such that the ordered pair $(A/{\Phi^{D}}, \leq^{D})$ is a sup-semilattice. Let $\leq^{D}$ be the binary relation on $A/{\Phi^{D}}$ defined as follows: for every $x$, $y\in A$, we have that $[x]_{\Phi^{D}}\leq^{D}[y]_{\Phi^{D}}$
if, and only if, there exists an element $x'\in [x]_{\Phi^{D}}$ and an element $y'\in [y]_{\Phi^{D}}$ for which $D(y',x')=y'$.
 
We begin by proving that $\leq^{D}$ is an order on $A/{\Phi^{D}}$. Let $x$ be an element of $A$. Then $[x]_{\Phi^{D}}\leq^{D}[x]_{\Phi^{D}}$ because $x\in [x]_{\Phi^{D}}$ and $D(x,x)=x$ according to \ref{D1}. This proves that $\leq^{D}$ is reflexive.

Let $x$ and $y$ be elements of $A$ and assume that $[x]_{\Phi^{D}}\leq^{D} [y]_{\Phi^{D}}$ and $[y]_{\Phi^{D}}\leq^{D}[x]_{\Phi^{D}}$. Thus we can find elements $x',x''\in [x]_{\Phi^{D}}$ and $y',y''\in [y]_{\Phi^{D}}$ such that $D(y',x')=y'$ and $D(x'',y'')=x''$. But it so happens that the following chain of equalities holds
\allowdisplaybreaks
\begin{align*}
[x]_{\Phi^{D}}&=[x'']_{\Phi^{D}}
\tag{$x''\in [x]_{\Phi^{D}}$}
\\&=
[D(x'',y'')]_{\Phi^{D}}
\tag{$D(x'',y'')=x''$}
\\&=
[D(x,y)]_{\Phi^{D}}
\tag{Prop.~\ref{PDOA}}
\\&=
[D(y,x)]_{\Phi^{D}}
\tag{Prop.~\ref{PDOC}}
\\&=
[D(y',x')]_{\Phi^{D}}
\tag{Prop.~\ref{PDOA}}
\\&=
[y']_{\Phi^{D}}
\tag{$D(y',x')=y'$}
\\&=
[y]_{\Phi^{D}}.
\tag{$y'\in [y]_{\Phi^{D}}$}
\end{align*}
Therefore $\leq^{D}$ is antisymmetric.

Let $x$, $y$ and $z$ be elements of $A$ and assume that $[x]_{\Phi^{D}}\leq^{D}[y]_{\Phi^{D}}$ and $[y]_{\Phi^{D}}\leq^{D} [z]_{\Phi^{D}}$. Thus we can find elements $x'\in [x]_{\Phi^{D}}$, $y',y''\in [y]_{\Phi^{D}}$ and $z'\in [z]_{\Phi^{D}}$ such that $D(y',x')=y'$ and $D(z',y'')=z'$. Before continuing with the proof, let us note that, from the fact that $y',y''\in [y]_{\Phi^{D}}$ we have that 
$[y']_{\Phi^{D}}=[y'']_{\Phi^{D}}$. Hence $(y',y'')\in \Phi^{D}$, i.e., $D(y',y'')=y'$ and $D(y'',y')=y''$. 
Moreover, the following claim also holds.
\begin{claim}\label{COrdI} $D(z',y')=z'$.
\end{claim}
Indeed, the following chain of equalities holds
\allowdisplaybreaks
\begin{align*}
D(z',y')&=D(D(z',y''),y')
\tag{$D(z',y'')=z'$}
\\&=
D(z',D(y'',y'))
\tag{by~\ref{D2}}
\\&=
D(z',y'')
\tag{$D(y'',y')=y''$}
\\&=
z'.
\tag{$D(z',y'')=z'$}
\end{align*}
This proves Claim~\ref{COrdI}.

To prove that $[x]_{\Phi^{D}}\leq^{D} [z]_{\Phi^{D}}$ it suffices to prove that $D(z',x')=z'$.
But it so happens that the following chain of equalities holds
\allowdisplaybreaks
\begin{align*}
D(z',x')&=D(D(z',y'),x')
\tag{Claim~\ref{COrdI}}
\\&=
D(z',D(y',x'))
\tag{by~\ref{D2}}
\\&=
D(z',y')
\tag{$D(y',x')=y'$}
\\&=
z'.
\tag{Claim~\ref{COrdI}}
\end{align*}

This proves that $\leq^{D}$ is transitive.

We conclude that $(A/{\Phi^{D}}, \leq^{D})$ is an ordered set.

We now prove that any pair of elements in $A/{\Phi^{D}}$ has a least upper bound. 

\begin{claim}\label{COrdII} Let $x$ and $y$ be elements of $A$. Then the least upper bound of $[x]_{\Phi^{D}}$ and $[y]_{\Phi^{D}}$, denoted by $[x]_{\Phi^{D}}\vee^{D} [y]_{\Phi^{D}}$, is 
$[D(x,y)]_{\Phi^{D}}$.
\end{claim}
We first prove that $[x]_{\Phi^{D}}\leq^{D} [D(x,y)]_{\Phi^{D}}$. To do this it suffices to check that $D(D(y,x),x)=D(y,x)$. But this follows from Lemma~\ref{LTech}.

We next prove that $[y]_{\Phi^{D}}\leq^{D} [D(x,y)]_{\Phi^{D}}$.  To do this it suffices to check that $D(D(x,y),y)=D(x,y)$. But this follows from Lemma~\ref{LTech}.


We next prove that $[D(x,y)]_{\Phi^{D}}$ is the least upper bound. Let $[z]_{\Phi^{D}}$ be an element of $A/{\Phi^{D}}$ such that $[x]_{\Phi^{D}}\leq^{D} [z]_{\Phi^{D}}$ and $[y]_{\Phi^{D}}\leq^{D} [z]_{\Phi^{D}}$. Then there exists elements $x'\in [x]_{\Phi^{D}}$, $y'\in [y]_{\Phi^{D}}$ and $z',z''\in [z]_{\Phi^{D}}$ such that $D(z',x')=z'$ and $D(z'',y')=z''$. Let us note that, from the fact that $z',z''\in [z]_{\Phi^{D}}$, we have that $[z']_{\Phi^{D}}=[z'']_{\Phi^{D}}$. Hence $(z',z'')\in \Phi^{D}$, i.e., $D(z',z'')=z'$ and $D(z'',z')=z''$. Before continuing with the proof, we will prove two auxiliary claims.
\begin{claim}\label{COrdU1} $D(z',x)=z'$.
\end{claim}
But it so happens that the following chain of equalities holds
\allowdisplaybreaks
\begin{align*}
D(z',x)&=
D(z',D(x,x'))
\tag{$x'\in [x]_{\Phi^{D}}$}
\\&=
D(z',D(x',x))
\tag{by~\ref{D3}}
\\&=
D(z',x')
\tag{$x'\in [x]_{\Phi^{D}}$}
\\&=
z'.
\tag{$D(z',x')=z'$}
\end{align*}
This proves Claim~\ref{COrdU1}.

\begin{claim}\label{COrdU2} $D(z'',y)=z''$.
\end{claim}
But it so happens that the following chain of equalities holds
\allowdisplaybreaks
\begin{align*}
D(z'',y)&=
D(z'',D(y,y'))
\tag{$y'\in [y]_{\Phi^{D}}$}
\\&=
D(z'',D(y',y))
\tag{by~\ref{D3}}
\\&=
D(z'',y')
\tag{$y'\in [y]_{\Phi^{D}}$}
\\&=
z''.
\tag{$D(z'',y')=z''$}
\end{align*}
This proves Claim~\ref{COrdU2}.

Now, to prove that $[D(x,y)]_{\Phi^{D}}\leq^{D}[z]_{\Phi^{D}}$ it suffices to check that $D(z'', D(x,y))=z''$. But it so happens that the following chain of equalities holds
\allowdisplaybreaks
\begin{align*}
D(z'', D(x,y))&=
D(D(z'',z'), D(x,y))
\tag{$D(z'',z')=z''$}
\\&=
D(D(D(z'',z'), x),y)
\tag{by~\ref{D2}}
\\&=
D(D(z'', D(z',x)),y)
\tag{by~\ref{D2}}
\\&=
D(D(z'', z'),y)
\tag{Claim~\ref{COrdU1}}
\\&=
D(z'',y)
\tag{$D(z'',z')=z''$}
\\&=
z''.
\tag{Claim~\ref{COrdU2}}
\end{align*}
This proves Claim~\ref{COrdII}.

This completes the proof that there exists a binary relation $\leq^{D}$ on $A/{\Phi^{D}}$ such that the ordered pair $(A/{\Phi^{D}},\leq^{D})$ or, equivalently, $(A/{\Phi^{D}},\vee^{D})$ is a sup-semilattice. We will call the relation $\leq^{D}$ the \emph{order relation associated to $\Phi^{D}$ on $A/{\Phi^{D}}$}.

On the other hand, from Claim~\ref{COrdII} it follows that the canonical projection $\mathrm{pr}_{\Phi^{D}}$ from $A$ to $A/{\Phi^{D}}$ determines a morphism from the left normal band $(A,D)$ to the left normal band $(A/{\Phi^{D}}, \vee^{D})$.

We, finally, prove that the pair $((A/{\Phi^{D}},\vee^{D}),\mathrm{pr}_{\Phi^{D}})$ is a universal morphism from 
$(A,D)$ to $\mathrm{In}_{\mathsf{Ssl},\mathsf{Lnb}}$.

Let $\mathbf{I}$ be a sup-semilattice and $h$ a morphism of left normal bands from $(A,D)$ to $\mathrm{In}_{\mathsf{Ssl},\mathsf{Lnb}}(\mathbf{I}) = (I,\vee)$, where $\vee$ is the supremum operation associated to the sup-semilattice $\mathbf{I} = (I,\leq)$. We want to show that there exists a unique sup-semilattice homomorphism
\[
h^{\flat}\colon (A/{\Phi^{D}},\leq^{D})\mor \mathbf{I}
\]
such that 
$h = \mathrm{In}_{\mathsf{Ssl},\mathsf{Lnb}}(h^{\flat})\circ \mathrm{pr}_{\Phi^{D}} = h^{\flat}\circ \mathrm{pr}_{\Phi^{D}}$. To do this we start by proving that there exists a mapping $h^{\flat}$ from $A/{\Phi^{D}}$ to $I$ such that $h = h^{\flat}\circ \mathrm{pr}_{\Phi^{D}}$. And for that, in turn, it suffices to prove that $\Phi^{D}\subseteq \mathrm{Ker}(h)$.

Since $h\colon (A,D)\mor (I,\vee)$ is a morphism of left normal bands, we have that 
\[
\vee\circ (h\times h)=h\circ D,\, i.e., \text{ for every } x,y\in A,\, h(D(x,y)) = h(x)\vee h(y).
\tag{S1}\label{S1}
\]
Now we prove the following claim.
\begin{claim}\label{CUnivSL1} The inclusion $\Phi^{D}\subseteq \mathrm{Ker}(h)$ holds.
\end{claim}
Let $x$ and $y$ be elements of $A$ and assume that $(x,y)\in \Phi^{D}$. Then, by Definition~\ref{DDORel}, $D(x,y)=x$ and $D(y,x)=y$.
But it so happens that the following chain of equalities holds
\begin{align*}
h(x)&=h(D(x,y))
\tag{$D(x,y)=x$}
\\&=
h(x)\vee h(y).
\tag{by~\ref{S1}}
\end{align*}
Thus, $h(x)=h(x)\vee h(y)$, i.e., $h(y)\leq h(x)$.

On the other hand, it so happens that the following chain of equalities holds
\begin{align*}
h(y)&=h(D(y,x))
\tag{$D(y,x)=y$}
\\&=
h(y)\vee h(x).
\tag{by~\ref{S1}}
\end{align*}
Thus, $h(y)=h(y)\vee h(x)$, i.e., $h(x)\leq h(y)$.

All in all, we conclude that $h(x)=h(y)$, i.e., that $(x,y)\in \mathrm{Ker}(h)$.

This proves Claim~\ref{CUnivSL1}.

The last claim entails, by the universal property of the quotient, that there exists a unique mapping $h^{\flat}$ from $A/{\Phi^{D}}$ to $I$ such that
\[
h = h^{\flat}\circ \mathrm{pr}_{\Phi^{D}}.
\tag{S2}\label{S2}
\]

\begin{claim}\label{CUnivSL2} The mapping $h^{\flat}$ determines a sup-semilattice homomorphism from $(A/{\Phi^{D}},\leq^{D})$ to $\mathbf{I}$.
\end{claim}
Let $[x]_{\Phi^{D}}, [y]_{\Phi^{D}}$ in $A/{\Phi^{D}}$. But it so happens that the following chain of equalities holds
\allowdisplaybreaks
\begin{align*}
h^{\flat}([x]_{\Phi^{D}}\vee [y]_{\Phi^{D}})&=
h^{\flat}([D(x,y)]_{\Phi^{D}})
\tag{Claim~\ref{COrdII}}
\\&=
h^{\flat}(\mathrm{pr}_{\Phi^{D}}(D(x,y)))
\tag{Def.~$\mathrm{pr}_{\Phi^{D}}$}
\\&=
h(D(x,y))
\tag{by~\ref{S2}}
\\&=
h(x)\vee h(y)
\tag{by~\ref{S1}}
\\&=
h^{\flat}(\mathrm{pr}_{\Phi^{D}}(x))
\vee
h^{\flat}(\mathrm{pr}_{\Phi^{D}}(y))
\tag{by~\ref{S2}}
\\&=
h^{\flat}([x]_{\Phi^{D}})
\vee
h^{\flat}([y]_{\Phi^{D}}).
\tag{Def.~$\mathrm{pr}_{\Phi^{D}}$}
\end{align*}
This proves Claim~\ref{CUnivSL2}.

It is obvious that, as morphisms, we have that $h=\mathrm{In}_{\mathsf{Ssl},\mathsf{Lnb}}(h^{\flat})\circ \mathrm{pr}_{\Phi^{D}}$.
%
%

To complete the proof it remains to prove that $h^{\flat}$ is the unique sup-semilattice homomorphism from $(A/{\Phi^{D}},\leq^{D})$ to $\mathbf{I}$ such that $h=\mathrm{In}_{\mathsf{Ssl},\mathsf{Lnb}}(h^{\flat})\circ \mathrm{pr}_{\Phi^{D}}$. We leave it to the reader to verify this last assertion.

%

We denote by $\mathrm{Sl}$ the left adjoint of $\mathrm{In}_{\mathsf{Ssl},\mathsf{Lnb}}$. Moreover, for a left normal band $(A,D)$, we let $\mathbf{Sl}(A,D)$ stand for $(A/{\Phi^{D}}, \leq^{D})$, which is the value of $\mathrm{Sl}$ at $(A,D)$, and we let $\mathrm{Sl}(A,D)$ stand for $A/{\Phi^{D}}$; and for a morphism $h$ from a left normal band $(A,D)$ to another $(B,E)$, we let 
$\mathrm{Sl}(h)$ stand for $(\mathrm{pr}_{\Phi^{E}}\circ h)^{\flat}$, the unique homomorphism of sup-semilattices from $(A/{\Phi^{D}}, \leq^{D})$ to $(B/{\Phi^{E}}, \leq^{E})$ such that 
$$
\mathrm{pr}_{\Phi^{E}}\circ h = (\mathrm{pr}_{\Phi^{E}}\circ h)^{\flat}\circ \mathrm{pr}_{\Phi^{D}}.
$$

This completes the proof of Proposition~\ref{PUnivLnb}.
\end{proof}

\section{On the category of P{\l}onka $\Sigma^{\neq 0}$-algebras}

In this section, for $\Sigma^{\neq 0}$, the subsignature of $\Sigma$ without $0$-ary operation symbols, we define the category of P{\l}onka $\Sigma^{\neq 0}$-algebras and, for later use, we establish some useful properties of the underlying structural P{\l}onka operator of a P{\l}onka $\Sigma^{\neq 0}$-algebra.


We begin by defining the category of P{\l}onka $\Sigma^{\neq 0}$-algebras.

\begin{definition}[See P{\l}onka~\cite{P74}]\label{DDOD4} 
Let $\mathbf{A}=(A,F)$ a $\Sigma^{\neq 0}$-algebra. A \emph{P{\l}onka operator for} $\mathbf{A}$ is a binary operation $D$ on $A$ such that $(A,D)$ is a left normal band and, for every $n\in\mathbb{N}-1$, every $\sigma\in \Sigma_{n}$, every $(x_{j})_{j\in n}\in A^{n}$ and every $y\in A$, the following conditions are satisfied 
\begin{align*}
D(F^{\mathbf{A}}_{\sigma}((x_{j})_{j\in n}), y)&=
F^{\mathbf{A}}_{\sigma}((D(x_{j},y))_{j\in n});
\tag{D4}\label{D4}
\\
D(y,F^{\mathbf{A}}_{\sigma}((x_{j})_{j\in n}))&=
D(y,D_{n}((x_{j})_{j\in n})).
\tag{D5}\label{D5}
\end{align*}
A \emph{P{\l}onka $\Sigma^{\neq 0}$-algebra} is an ordered pair $(\mathbf{A},D)$ in which $\mathbf{A}=(A,F)$ is a $\Sigma^{\neq 0}$-algebra and $D$ a P{\l}onka operator for 
$\mathbf{A}$.

Given two P{\l}onka $\Sigma^{\neq 0}$-algebras $(\mathbf{A},D)$ and $(\mathbf{B},E)$, a morphism $h$ from $(\mathbf{A},D)$ to $(\mathbf{B},E)$ is an ordered triple $((\mathbf{A},D),h,(\mathbf{B},E))$, abbreviated to $h\colon (\mathbf{A},D)\mor (\mathbf{B},E)$, where $h$ is a 
$\Sigma^{\neq 0}$-homomorphism from $\mathbf{A}$ to $\mathbf{B}$ which is also a morphism from the left normal band $(A,D)$ to the left normal band $(B,E)$, i.e., which is such that $ E\circ (h\times h)=h\circ D$.
%
%

We denote by $\mbox{\sffamily{\upshape{P{\l}Alg}}}(\Sigma^{\neq 0})$ the category whose objects are P{\l}onka $\Sigma^{\neq 0}$-algebras $(\mathbf{A},D)$ and whose morphisms are the morphisms $h\colon (\mathbf{A},D)\mor (\mathbf{B},E)$ of P{\l}onka $\Sigma^{\neq 0}$-algebras. 
For a P{\l}onka $\Sigma^{\neq 0}$-algebra $(\mathbf{A},D)$, we will call $D$ the \emph{underlying structural P{\l}onka operator of} $(\mathbf{A},D)$, $(A,D)$ the \emph{underlying left normal band of} $(\mathbf{A},D)$, and $\mathbf{A}$ the \emph{underlying $\Sigma^{\neq 0}$-algebra of} $(\mathbf{A},D)$.
\end{definition}

\begin{remark}
What we have called ``P{\l}onka operator for $\mathbf{A}$'' was called ``partition function for $\mathbf{A}$'' by  P{\l}onka~\cite{P67}.
\end{remark}

\begin{remark}
Let us note that the condition $\mathrm{D}4$ is a distributive law (of a monad over another monad). Moreover, $\mbox{\sffamily{\upshape{P{\l}Alg}}}(\Sigma^{\neq 0})$ is a variety.
\end{remark}

\begin{remark}\label{Sigmadelta}
We have defined a P{\l}onka $\Sigma^{\neq 0}$-algebra as an ordered pair $(\mathbf{A},D)$ in which $\mathbf{A}=(A,F)$ is a $\Sigma^{\neq 0}$-algebra and $D$ a binary operation on $A$ satisfying the laws~\ref{D1}--\ref{D5}. However, there are occasions when it would be formally more appropriate to use, instead of the category $\mbox{\sffamily{\upshape{P{\l}Alg}}}(\Sigma^{\neq 0})$, the category 
$\mbox{\sffamily{\upshape{P{\l}Alg}}}(\Sigma^{\neq 0}_{\delta})$, in which $\Sigma^{\neq 0}_{\delta}$ is the signature obtained from $\Sigma^{\neq 0}$ by adding a new binary operation symbol $\delta$ and whose objects are ordered pairs $(A,(F,F_{\delta}))$, where $F$ is a structure of $\Sigma^{\neq 0}$-algebra on $A$ corresponding to the subsignature $\Sigma^{\neq 0}$ of $\Sigma^{\neq 0}_{\delta}$ and $F_{\delta}$ a binary operation on $A$ such that $(A,(F,F_{\delta}))$ satisfies the laws~\ref{D1}--\ref{D5}. Note that $\mbox{\sffamily{\upshape{P{\l}Alg}}}(\Sigma^{\neq 0})$ and $\mbox{\sffamily{\upshape{P{\l}Alg}}}(\Sigma^{\neq 0}_{\delta})$ are isomorphic.
\end{remark}


We next establish some technical results on the P{\l}onka operator of a P{\l}onka $\Sigma$-algebra that will be useful later on.

\begin{proposition}\label{PTech2} Let $(\mathbf{A},D)$ be a P{\l}onka $\Sigma^{\neq 0}$-algebra. Then, for every $n\in \mathbb{N}-1$, every $k\in n$, every $\sigma\in \Sigma_{n}$ and every $(x_{j})_{j\in n}\in A^{n}$, the following equality holds
\[
D(F^{\mathbf{A}}_{\sigma}((x_{j})_{j\in n}),x_{k})= 
F^{\mathbf{A}}_{\sigma}((x_{j})_{j\in n}).
\]
\end{proposition}
\begin{proof}
The following chain of equalities holds
\allowdisplaybreaks
\begin{align*}
D(F^{\mathbf{A}}_{\sigma}((x_{j})_{j\in n}),x_{k})&=
D(D(F^{\mathbf{A}}_{\sigma}((x_{j})_{j\in n}),F^{\mathbf{A}}_{\sigma}((x_{j})_{j\in n})),
x_{k})
\tag{by~\ref{D1}}
\\&=
D(D(F^{\mathbf{A}}_{\sigma}((x_{j})_{j\in n}),D_{n}((x_{j})_{j\in n})),
x_{k})
\tag{by~\ref{D5}}
\\&=
D(F^{\mathbf{A}}_{\sigma}((x_{j})_{j\in n}),
D(D_{n}((x_{j})_{j\in n}),
x_{k}))
\tag{by~\ref{D2}}
\\&=
D(F^{\mathbf{A}}_{\sigma}((x_{j})_{j\in n}), 
D_{n}((x_{j})_{j\in n})
)
\tag{Lemma~\ref{LTech}}
\\&=
D(F^{\mathbf{A}}_{\sigma}((x_{j})_{j\in n}),F^{\mathbf{A}}_{\sigma}((x_{j})_{j\in n}))
\tag{by~\ref{D5}}
\\&=
F^{\mathbf{A}}_{\sigma}((x_{j})_{j\in n}).
\tag{by~\ref{D1}}
\end{align*}

This completes the proof.
\end{proof}

\begin{proposition}\label{PTech3} Let $(\mathbf{A},D)$ be a P{\l}onka $\Sigma^{\neq 0}$-algebra. Then, for every $n\in \mathbb{N}-1$, every $k\in n$, every $\sigma\in \Sigma_{n}$, every $(x_{j})_{j\in n}\in A^{n}$ and every $y\in A$, the following equality holds
\[
D(F^{\mathbf{A}}_{\sigma}((x_{j})_{j\in n}),y)= 
F^{\mathbf{A}}_{\sigma}(x_{0},\cdots, D(x_{k},y), \cdots, x_{n-1}).
\]
\end{proposition}
\begin{proof}
The following chain of equalities holds.
\begin{flushleft}
$D(
F^{\mathbf{A}}_{\sigma}((x_{j})_{j\in n}),
y
)$
\allowdisplaybreaks
\begin{align*}
&=
D(
D(F^{\mathbf{A}}_{\sigma}((x_{j})_{j\in n}),x_{k}),
y
)
\tag{Prop.~\ref{PTech2}}
\\&=
D(
F^{\mathbf{A}}_{\sigma}((x_{j})_{j\in n}),
D(x_{k},y)
)
\tag{by~\ref{D2}}
\\&=
F^{\mathbf{A}}_{\sigma}(
D(x_{0}, D(x_{k},y)),
\cdots ,
D(x_{k}, D(x_{k},y)),
\cdots,
D(x_{n-1}, D(x_{k},y))
)
\tag{by~\ref{D4}}
\\&=
F^{\mathbf{A}}_{\sigma}(
D(x_{0}, D(x_{k},y)),
\cdots ,
D(D(x_{k}, x_{k}),y),
\cdots,
D(x_{n-1}, D(x_{k},y))
)
\tag{by~\ref{D2}}
\\&=
F^{\mathbf{A}}_{\sigma}(
D(x_{0}, D(x_{k},y)),
\cdots ,
D(x_{k},y),
\cdots,
D(x_{n-1}, D(x_{k},y))
)
\tag{by~\ref{D1}}
\\&=
F^{\mathbf{A}}_{\sigma}(
D(x_{0}, D(x_{k},y)),
\cdots ,
D(D(x_{k},y),D(x_{k},y)),
\cdots,
\\&\qquad\qquad\qquad\qquad\qquad\qquad\qquad\qquad
\qquad\qquad\qquad
D(x_{n-1}, D(x_{k},y))
)
\tag{by~\ref{D1}}
\\&=
D(
F^{\mathbf{A}}_{\sigma}(
x_{0},
\cdots ,
D(x_{k},y),
\cdots,
x_{n-1}
),
D(x_{k},y)
)
\tag{by~\ref{D4}}
\\&=
F^{\mathbf{A}}_{\sigma}(
x_{0},
\cdots ,
D(x_{k},y),
\cdots,
x_{n-1}
).
\tag{Prop.~\ref{PTech2}}
\end{align*}
\end{flushleft}

This completes the proof.
\end{proof}

\begin{proposition}\label{PDHom} Let $(\mathbf{A},D)$ be a P{\l}onka $\Sigma^{\neq 0}$-algebra. Then, for every $n\in \mathbb{N}-1$, every $\sigma\in \Sigma_{n}$ and every $(x_{j})_{j\in n}, (y_{j})_{j\in n}\in A^{n}$, the following equality holds
\[
D(F^{\mathbf{A}}_{\sigma}((x_{j})_{j\in n}),F^{\mathbf{A}}_{\sigma}((y_{j})_{j\in n}))= 
F^{\mathbf{A}}_{\sigma}((D(x_{j},y_{j}))_{j\in n}).
\]
\end{proposition}
\begin{proof}
The following chain of equalities holds.
\begin{flushleft}
$D(F^{\mathbf{A}}_{\sigma}((x_{j})_{j\in n}),F^{\mathbf{A}}_{\sigma}((y_{j})_{j\in n}))$
\allowdisplaybreaks
\begin{align*}
&=D(F^{\mathbf{A}}_{\sigma}((x_{j})_{j\in n}),D_{n}((y_{j})_{j\in n}))
\tag{by~\ref{D5}}
\\&=
D(F^{\mathbf{A}}_{\sigma}((x_{j})_{j\in n}),
D(y_{0}, 
D_{n-1}((y_{j})_{j\in [1,n-1]})))
\tag{Def.~\ref{DIt}}
\\&=
D(D(F^{\mathbf{A}}_{\sigma}((x_{j})_{j\in n}),
y_{0}), 
D_{n-1}((y_{j})_{j\in [1,n-1]}))
\tag{by~\ref{D2}}
\\&=
D(
F^{\mathbf{A}}_{\sigma}(D(x_{0},y_{0}),x_{1},\cdots,x_{n-1}),
D_{n-1}((y_{j})_{j\in [1,n-1]}))
\tag{Prop.~\ref{PTech3}}
\\&\qquad\vdots
\\&=
D(
F^{\mathbf{A}}_{\sigma}(D(x_{0},y_{0}),\cdots D(x_{n-3},y_{n-3}),x_{n-2},x_{n-1}),
D(y_{n-2},y_{n-1}))
\\&=
D(D(F^{\mathbf{A}}_{\sigma}(D(x_{0},y_{0}),\cdots D(x_{n-3},y_{n-3}),x_{n-2},x_{n-1}),y_{n-2}),
y_{n-1})
\tag{by~\ref{D2}}
\\&=
D(F^{\mathbf{A}}_{\sigma}(D(x_{0},y_{0}),\cdots D(x_{n-2},y_{n-2}),x_{n-1}),
y_{n-1})
\tag{Prop.~\ref{PTech3}}
\\&=
F^{\mathbf{A}}_{\sigma}((D(x_{j},y_{j}))_{j\in n}).
\tag{Prop.~\ref{PTech3}}
\end{align*}
\end{flushleft}

This completes the proof.
\end{proof}

\begin{remark}
For a P{\l}onka $\Sigma^{\neq 0}$-algebra $(\mathbf{A},D)$, Proposition~\ref{PDHom} can be rephrased as: the underlying P{\l}onka operator $D$ of $(\mathbf{A},D)$ is a $\Sigma^{\neq 0}$-homomorphism from the $\Sigma^{\neq 0}$-algebra $\mathbf{A}^{2}$ to the $\Sigma^{\neq 0}$-algebra $\mathbf{A}$. On the other hand, Proposition~\ref{PDHom} will be fundamental to show later on that the P{\l}onka construction is part of an adjoint.
\end{remark}

Before continuing with the presentation of the technical results mentioned above and taking into account the first part of the just stated remark, the relationship between the categories $\mbox{\sffamily{\upshape{P{\l}Alg}}}(\Sigma^{\neq 0})$, $\mathsf{Alg}(\Sigma^{\neq 0})\otimes\mathsf{Lnb}$ and $\mathsf{Alg}(\Sigma^{\neq 0})$ is discussed below.

\begin{definition}[See Manes~\cite{em76}] \label{tensor}
The category $\mathsf{Alg}(\Sigma^{\neq 0})\otimes\mathsf{Lnb}$, the \emph{tensor product of} $\mathsf{Alg}(\Sigma^{\neq 0})$ \emph{and} $\mathsf{Lnb}$, has as objects the ordered triples $(A,F,D)$ such that $(A,F)$ is a $\Sigma^{\neq 0}$-algebra, $(A,D)$ a left normal band and $D$ a $\Sigma^{\neq 0}$-homomorphism from the $\Sigma^{\neq 0}$-algebra $\mathbf{A}^{2}$ to the $\Sigma^{\neq 0}$-algebra $\mathbf{A}$; and as morphisms from $(A,F,D)$ to $(B,G,E)$ the ordered triples $((A,F,D),f,(B,G,E))$, abbreviated to $f\colon (A,F,D)\mor (B,G,E)$, in which $f$ is simultaneously a $\Sigma^{\neq 0}$-homomorphism from $(A,F)$ to $(B,G)$ and a morphism of left normal bands from $(A,D)$ to $(B,E)$. Moreover, we let $P_{\Sigma^{\neq 0}}$ stand for the canonical functor from $\mathsf{Alg}(\Sigma^{\neq 0})\otimes\mathsf{Lnb}$ to $\mathsf{Alg}(\Sigma^{\neq 0})$ and by $Q_{\Sigma^{\neq 0}}$ the canonical functor from $\mathsf{Alg}(\Sigma^{\neq 0})\otimes\mathsf{Lnb}$ to $\mathsf{Lnb}$. 
\end{definition} 

\begin{proposition}\label{Pl-Sigmalad}
There exists a full embedding $J_{\Sigma^{\neq 0}}$ of the category $\mbox{\sffamily{\upshape{P{\l}Alg}}}(\Sigma^{\neq 0})$ into the category $\mathsf{Alg}(\Sigma^{\neq 0})\otimes\mathsf{Lnb}$. Moreover, the functor $P_{\Sigma^{\neq 0}}\circ J_{\Sigma^{\neq 0}}$ from $\mbox{\sffamily{\upshape{P{\l}Alg}}}(\Sigma^{\neq 0})$ to $\mathsf{Alg}(\Sigma^{\neq 0})$ has a left adjoint.
\end{proposition}

\begin{proof}
The object mapping of $J_{\Sigma^{\neq 0}}$ sends $(\mathbf{A},D)$, where $\mathbf{A} = (A,F)$, to $(A,F,D)$. That $P_{\Sigma^{\neq 0}}\circ J_{\Sigma^{\neq 0}}$ has a left adjoint follows from Theorem~7.3 (b), on p. 116, of Barr and Wells~\cite{BW85}.
\end{proof}


\begin{proposition}\label{PTech4} Let $(\mathbf{A},D)$ be a P{\l}onka $\Sigma^{\neq 0}$-algebra. Then, for every $m,n\in \mathbb{N}-1$, every $\sigma\in \Sigma_{n}$, every $(x_{j})_{j\in n}\in A^{n}$, every $\varphi\colon m\mor n$ and the sequence $(x_{\varphi(k)})_{k\in m}\in A^{m}$, the following equality holds
\[
D(F^{\mathbf{A}}_{\sigma}((x_{j})_{j\in n}),
D_{m}((x_{\varphi(k)})_{k\in m})
)= 
F^{\mathbf{A}}_{\sigma}((x_{j})_{j\in n}).
\]
\end{proposition}
\begin{proof}
The following chain of equalities holds.
\begin{flushleft}
$D(F^{\mathbf{A}}_{\sigma}((x_{j})_{j\in n}),
D_{m}((x_{\varphi(k)})_{k\in m})
)$
\allowdisplaybreaks
\begin{align*}
&=D(D(F^{\mathbf{A}}_{\sigma}((x_{j})_{j\in n}),
F^{\mathbf{A}}_{\sigma}((x_{j})_{j\in n}))
,
D_{m}((x_{\varphi(k)})_{k\in m})
)
\tag{by~\ref{D1}}
\\&=
D(D(F^{\mathbf{A}}_{\sigma}((x_{j})_{j\in n}),
D_{n}((x_{j})_{j\in n}))
,
D_{m}((x_{\varphi(k)})_{k\in m})
)
\tag{by~\ref{D1}}
\\&=
D(F^{\mathbf{A}}_{\sigma}((x_{j})_{j\in n}),
D(
D_{n}((x_{j})_{j\in n})
,
D_{m}((x_{\varphi(k)})_{k\in m}))
)
\tag{by~\ref{D2}}
\\&=
D(F^{\mathbf{A}}_{\sigma}((x_{j})_{j\in n}),
D_{n}((x_{j})_{j\in n})
)
\tag{Lemma~\ref{LTech2}}
\\&=
D(F^{\mathbf{A}}_{\sigma}((x_{j})_{j\in n}),
F^{\mathbf{A}}_{\sigma}((x_{j})_{j\in n}))
\tag{by~\ref{D5}}
\\&=
F^{\mathbf{A}}_{\sigma}((x_{j})_{j\in n}).
\tag{by~\ref{D1}}
\end{align*}
\end{flushleft}

This completes the proof.
\end{proof}

\begin{proposition}\label{CTech4} Let $(\mathbf{A},D)$ be a P{\l}onka $\Sigma^{\neq 0}$-algebra. Then, for every $n\in \mathbb{N}-1$, every $\sigma\in \Sigma_{n}$ and every $(x_{j})_{j\in n}\in A^{n}$, 
the following equality holds
\[
D(F^{\mathbf{A}}_{\sigma}((x_{j})_{j\in n}),
F^{\mathbf{A}}_{\sigma}(
(D_{n}((x_{j})_{j\in n}))_{j\in n}
)
)= 
F^{\mathbf{A}}_{\sigma}((x_{j})_{j\in n}).
\]
Note that $(D_{n}((x_{j})_{j\in n}))_{j\in n}$, the argument of $F^{\mathbf{A}}_{\sigma}$ in the left-hand side of the above equation, is the family in $A^{n}$ which is constantly $D_{n}((x_{j})_{j\in n})$.
\end{proposition}

\begin{proof}
The following chain of equalities holds.
\begin{flushleft}
$D(F^{\mathbf{A}}_{\sigma}((x_{j})_{j\in n}),
F^{\mathbf{A}}_{\sigma}(
(D_{n}((x_{j})_{j\in n}))_{j\in n}
)
)$
\allowdisplaybreaks
\begin{align*}
&=
D(F^{\mathbf{A}}_{\sigma}((x_{j})_{j\in n}),
D_{n}(
(D_{n}((x_{j})_{j\in n}))_{j\in n}
)
)
\tag{by~\ref{D5}}
\\&=
D(F^{\mathbf{A}}_{\sigma}((x_{j})_{j\in n}),
D(D_{n}((x_{j})_{j\in n}),
D_{n-1}(
(D_{n}((x_{j})_{j\in n}))_{j\in [1,n-1]}
)
)
\tag{Def.~\ref{DIt}}
\\&=
D(D(
F^{\mathbf{A}}_{\sigma}((x_{j})_{j\in n}),
D_{n}((x_{j})_{j\in n})
),
D_{n-1}(
(D_{n}((x_{j})_{j\in n}))_{j\in [1,n-1]}
))
\tag{by~\ref{D2}}
\\&=
D(
F^{\mathbf{A}}_{\sigma}((x_{j})_{j\in n}),
D_{n-1}(
(D_{n}((x_{j})_{j\in n}))_{j\in [1,n-1]}
))
\tag{Prop.~\ref{PTech4}}
\\&\qquad\vdots
\\&=
D(F^{\mathbf{A}}_{\sigma}((x_{j})_{j\in n}),
D(D_{n}((x_{j})_{j\in n}),D_{n}((x_{j})_{j\in n})))
\\&=
D(D(F^{\mathbf{A}}_{\sigma}((x_{j})_{j\in n}),D_{n}((x_{j})_{j\in n})),
D_{n}((x_{j})_{j\in n}))
\tag{by~\ref{D2}}
\\&=
D(F^{\mathbf{A}}_{\sigma}((x_{j})_{j\in n}),D_{n}((x_{j})_{j\in n}))
\tag{Prop.~\ref{PTech4}}
\\&=
F^{\mathbf{A}}_{\sigma}((x_{j})_{j\in n}).
\tag{Prop.~\ref{PTech4}}
\end{align*}
\end{flushleft}

This completes the proof.
\end{proof}

\begin{remark} Proposition~\ref{CTech4}, together with Proposition~\ref{PDHom}, will be fundamental to show later on that the P{\l}onka construction is part of an adjunction.
\end{remark}

In the following proposition we state that, for a P{\l}onka $\Sigma^{\neq 0}$-algebra $(\mathbf{A},D)$, the congruence induced by $D$ on the underlying left normal band of $(\mathbf{A},D)$ (this by virtue of Propositions~\ref{PDOA}) is a congruence on $(\mathbf{A},D)$.

\begin{proposition}\label{PDOD4Cong} Let $(\mathbf{A},D)$ be a P{\l}onka $\Sigma^{\neq 0}$-algebra. Then the congruence $\Phi^{D}$ induced by $D$ on $(A,D)$, the underlying left normal band of $(\mathbf{A},D)$, is also a congruence on $\mathbf{A}$, the underlying $\Sigma^{\neq 0}$-algebra of $(\mathbf{A},D)$. Therefore $\Phi^{D}$ is a congruence on $(\mathbf{A},D)$.
\end{proposition}

\begin{proof} 
Let $n$ be an element of $\mathbb{N}-1$, $\sigma\in \Sigma_{n}$, $(x_{j})_{j\in n}$ and $(y_{j})_{j\in n}\in A^{n}$ such that, for every $j\in n$, $(x_{j},y_{j})\in\Phi^{D}$, i.e., $D(x_{j},y_{j})=x_{j}$ and $D(y_{j},x_{j})=y_{j}$.
We now prove that 
$
(F^{\mathbf{A}}_{\sigma}((x_{j})_{j\in n}), F^{\mathbf{A}}_{\sigma}((y_{j})_{j\in n}))\in \Phi^{D}.
$

Note that the following chain of equalities holds
\allowdisplaybreaks
\begin{align*}
D(F^{\mathbf{A}}_{\sigma}((x_{j})_{j\in n}), F^{\mathbf{A}}_{\sigma}((y_{j})_{j\in n}))
&=
F^{\mathbf{A}}_{\sigma}((D(x_{j},y_{j}))_{j\in n})
\tag{Prop.~\ref{PDHom}}
\\&=
F^{\mathbf{A}}_{\sigma}((x_{j})_{j\in n}).
\tag{$D(x_{j},y_{j})=x_{j}$}
\end{align*}

By a similar argument we obtain that 
$$
D(F^{\mathbf{A}}_{\sigma}((y_{j})_{j\in n}), F^{\mathbf{A}}_{\sigma}((x_{j})_{j\in n})) = 
F^{\mathbf{A}}_{\sigma}((y_{j})_{j\in n}).
$$


We will call $\Phi^{D}$ the \emph{congruence induced by $D$ on $(\mathbf{A},D)$}.
\end{proof}

\section{On the relationship between $\mbox{\sffamily{\upshape{P{\l}Alg}}}(\Sigma^{\neq 0})$ and $\int^{\mathsf{Ssl}}\mathrm{Isys}_{\Sigma^{\neq 0}}$: P{\l}onka adjunction}

In this section we prove that from $\mbox{\sffamily{\upshape{P{\l}Alg}}}(\Sigma^{\neq 0})$ to $\int^{\mathsf{Ssl}}\mathrm{Isys}_{\Sigma^{\neq 0}}$ there exists a functor $\mathrm{Is}_{\Sigma^{\neq 0}}$ which has a left adjoint: the P{\l}onka construction, whose object part is the one that assigns to an inductive system of $\Sigma^{\neq 0}$-algebras its P{\l}onka sum.  


We begin by proving that, for a P{\l}onka $\Sigma^{\neq 0}$-algebra $(\mathbf{A},D)$, every equivalence class with respect to $\Phi^{D}$ is a closed subset of $\mathbf{A}$.

\begin{proposition}\label{PClos} Let $(\mathbf{A},D)$ be a P{\l}onka $\Sigma^{\neq 0}$-algebra, $\Phi^{D}$ the congruence induced by $D$ on $(\mathbf{A},D)$, $x\in A$ and $\mathfrak{X}=[x]_{\Phi^{D}}$ its equivalence class. Then $\mathfrak{X}$ is a closed subset of $\mathbf{A}$.
\end{proposition}

\begin{proof}
Let $n$ be an element of $\mathbb{N}-1$, $\sigma\in \Sigma_{n}$ and $(x_{j})_{j\in n}\in \mathfrak{X}^{n}$. Hence, for every $j\in n$, $(x_{j},x)\in \Phi^{D}$. Therefore, by Proposition~\ref{PDOD4Cong}, we have that
\[
(F^{\mathbf{A}}_{\sigma}((x_{j})_{j\in n}), F^{\mathbf{A}}_{\sigma}((x)_{j\in n}))\in \Phi^{D},
\tag{C1}\label{C1}
\]
where $(x)_{j\in n}$ is the family in $A^{n}$ which is constantly $x$.

To prove that $F^{\mathbf{A}}_{\sigma}((x_{j})_{j\in n})\in \mathfrak{X}$, it suffices to verify that $(x,F^{\mathbf{A}}_{\sigma}((x_{j})_{j\in n}))\in \Phi^{D}$, i.e., that $D(x,F^{\mathbf{A}}_{\sigma}((x_{j})_{j\in n}))=x$ and $D(F^{\mathbf{A}}_{\sigma}((x_{j})_{j\in n}), x)=F^{\mathbf{A}}_{\sigma}((x_{j})_{j\in n})$.

The following chain of equalities holds.
\allowdisplaybreaks
\begin{align*}
D(x,F^{\mathbf{A}}_{\sigma}((x_{j})_{j\in n}))&=
D(x,
D(F^{\mathbf{A}}_{\sigma}((x_{j})_{j\in n}),
F^{\mathbf{A}}_{\sigma}((x)_{j\in n}))
)
\tag{by~\ref{C1}}
\\&=
D(x,
D(
F^{\mathbf{A}}_{\sigma}((x)_{j\in n}),
F^{\mathbf{A}}_{\sigma}((x_{j})_{j\in n})
)
)
\tag{by~\ref{D3}}
\\&=
D(x,
F^{\mathbf{A}}_{\sigma}((x)_{j\in n})
)
\tag{by~\ref{C1}}
\\&=
D(x,
D_{n}((x)_{j\in n})
)
\tag{by~\ref{D5}}
\\&=
x.
\tag{Lemma~\ref{LTech2}}
\end{align*}

The following chain of equalities holds.
\allowdisplaybreaks
\begin{align*}
D(F^{\mathbf{A}}_{\sigma}((x_{j})_{j\in n}), x)
&=
D(
D(F^{\mathbf{A}}_{\sigma}((x_{j})_{j\in n}), 
F^{\mathbf{A}}_{\sigma}((x)_{j\in n})
),
x)
\tag{by~\ref{C1}}
\\&=
D(
D(F^{\mathbf{A}}_{\sigma}((x_{j})_{j\in n}), 
D_{n}((x)_{j\in n})
),
x)
\tag{by~\ref{D5}}
\\&=
D(F^{\mathbf{A}}_{\sigma}((x_{j})_{j\in n}), 
D(
D_{n}((x)_{j\in n}),
x))
\tag{by~\ref{D2}}
\\&=
D(F^{\mathbf{A}}_{\sigma}((x_{j})_{j\in n}), 
D_{n}((x)_{j\in n})
)
\tag{Lemma~\ref{LTech}}
\\&=
D(F^{\mathbf{A}}_{\sigma}((x_{j})_{j\in n}), 
F^{\mathbf{A}}_{\sigma}((x)_{j\in n})
)
\tag{by~\ref{D5}}
\\&=
F^{\mathbf{A}}_{\sigma}((x_{j})_{j\in n}).
\tag{by~\ref{C1}}
\end{align*}

This completes the proof.
\end{proof}

\begin{definition}\label{DConsSub} Let $(\mathbf{A},D)$ be a P{\l}onka $\Sigma^{\neq 0}$-algebra, $\Phi^{D}$ the congruence induced by $D$ on $(\mathbf{A},D)$, $x\in A$ and $\mathfrak{X}=[x]_{\Phi^{D}}$ its equivalence class. Then we denote by $\boldsymbol{\mathfrak{X}}$ the $\Sigma^{\neq 0}$-algebra determined by $\mathfrak{X}$.  
\end{definition}

In the following proposition we state that if two equivalence classes are related by means of the order introduced in Proposition~\ref{PUnivLnb}, then there exist a $\Sigma^{\neq 0}$-homomorphism from the subalgebra associated to the smaller one to the subalgebra associated to the larger one.

\begin{proposition}\label{PConsHom} Let $(\mathbf{A},D)$ be a P{\l}onka $\Sigma^{\neq 0}$-algebra, $\Phi^{D}$ the congruence induced by $D$ on $(\mathbf{A},D)$ and $\mathfrak{X}=[x]_{\Phi^{D}}$, $\mathfrak{Y}=[y]_{\Phi^{D}}$  two equivalence classes of $A/{\Phi^{D}}$ such that $(\mathfrak{X}, \mathfrak{Y})\in\leq^{D}$. Then the mapping defined as follows
\[
f_{\mathfrak{X},\mathfrak{Y}}
\left\lbrace
\begin{array}{ccl}
[x]_{\Phi^{D}}&\mor&[y]_{\Phi^{D}}\\
z&\longmapsto& D(z, y)
\end{array}.
\right.
\]
determines a $\Sigma^{\neq 0}$-homomorphism from $\boldsymbol{\mathfrak{X}}$ to $\boldsymbol{\mathfrak{Y}}$.
\end{proposition}

\begin{proof}

We need to prove several statements to conclude that $f^{\mathfrak{X},\mathfrak{Y}}$ is a $\Sigma^{\neq 0}$-homomor\-phism.

\begin{claim}\label{CConsHomI} For every $z\in [x]_{\Phi^{D}}$, $f_{\mathfrak{X},\mathfrak{Y}}(z)$ is well-defined, i.e., for every $y'\in [y]_{\Phi^{D}}$, the equality $D(z,y)=D(z,y')$ holds.
\end{claim}

The following chain of equalities holds
\allowdisplaybreaks
\begin{align*}
D(z,y)&=D(z, D(y,y'))
\tag{$(y,y')\in \Phi^{D}$}
\\&=
D(z, D(y',y))
\tag{by~\ref{D3}}
\\&=
D(z,y').
\tag{$(y,y')\in \Phi^{D}$}
\end{align*}

This proves Claim~\ref{CConsHomI}.

\begin{claim}\label{CConsHomII} For every $z\in [x]_{\Phi^{D}}$,  $f_{\mathfrak{X},\mathfrak{Y}}(z)\in [y]_{\Phi^{D}}$.
\end{claim}

By hypothesis we have that $(\mathfrak{X},\mathfrak{Y})\in \leq^{D}$. Therefore, by unpacking the definition of $\leq^{D}$, as stated in Proposition~\ref{PUnivLnb}, we have that there exists $x'\in [x]_{\Phi^{D}}$ and $y'\in[y]_{\Phi^{D}}$ such that $D(y', x')=y'$.

Now, taking into account that $y'\in[y]_{\Phi^{D}}$ and Claim~\ref{CConsHomI}, to prove that $f_{\mathfrak{X},\mathfrak{Y}}(z)=D(z,y)$ is an element of $[y]_{\Phi^{D}}$ it suffices to prove that $[D(z,y')]_{\Phi^{D}}= [y']_{\Phi^{D}}$, i.e., that 
$D(D(z,y'), y')=D(z,y')$ and  $D(y', D(z,y'))=y'$. But these two equalities follow from Lemma~\ref{LTech}.

%

This proves Claim~\ref{CConsHomII}.

\begin{claim}\label{CConsHomIII} $f_{\mathfrak{X},\mathfrak{Y}}$ is a $\Sigma^{\neq 0}$-homomorphism.
\end{claim}
Let $n$ be an element of $\mathbb{N}-1$, $\sigma \in \Sigma_{n}$ and $(z_{j})_{j\in n}\in \mathfrak{X}^{n}$. We need to prove that 
\[
f_{\mathfrak{X},\mathfrak{Y}}(
F^{\boldsymbol{\mathfrak{X}}}_{\sigma}((
z_{j}
)_{j\in n})
)
=
F^{\boldsymbol{\mathfrak{Y}}}_{\sigma}((
f_{\mathfrak{X},\mathfrak{Y}}(z_{j})
)_{j\in n}).
\]

The following chain of equalities holds
\allowdisplaybreaks
\begin{align*}
f_{\mathfrak{X},\mathfrak{Y}}(
F^{\boldsymbol{\mathfrak{X}}}_{\sigma}((
z_{j}
)_{j\in n})
)&=
f_{\mathfrak{X},\mathfrak{Y}}(
F^{\mathbf{A}}_{\sigma}((
z_{j}
)_{j\in n})
)
\tag{Def.~\ref{DConsSub}}
\\&=
D(
F^{\mathbf{A}}_{\sigma}((
z_{j}
)_{j\in n}), y
)
\tag{Def. $f^{\mathfrak{X},\mathfrak{Y}}$}
\\&=
D(
F^{\mathbf{A}}_{\sigma}((
z_{j}
)_{j\in n}), 
D(
y, D_{n}((y)_{j\in n})
)
)
\tag{Lemma~\ref{LTech}}
\\&=
D(
F^{\mathbf{A}}_{\sigma}((
z_{j}
)_{j\in n}), 
D(
D_{n}((y)_{j\in n}),y
)
)
\tag{by~\ref{D2}}
\\&=
D(
F^{\mathbf{A}}_{\sigma}((
z_{j}
)_{j\in n}), 
D_{n}((y)_{j\in n})
)
\tag{Lemma~\ref{LTech}}
\\&=
D(
F^{\mathbf{A}}_{\sigma}((
z_{j}
)_{j\in n}),
F^{\mathbf{A}}_{\sigma}((
y
)_{j\in n})
)
\tag{by~\ref{D5}}
\\&=
F^{\mathbf{A}}_{\sigma}((
D(
z_{j},y
)
)_{j\in n}
)
\tag{Prop.~\ref{PConsHom}}
\\&=
F^{\mathbf{A}}_{\sigma}((
f_{\mathfrak{X},\mathfrak{Y}}(
z_{j})
)_{j\in n}
)
\tag{Def. $f_{\mathfrak{X},\mathfrak{Y}}$}
\\&=
F^{\boldsymbol{\mathfrak{Y}}}_{\sigma}((
f_{\mathfrak{X},\mathfrak{Y}}(
z_{j})
)_{j\in n}
).
\tag{Def.~\ref{DConsSub}}
\end{align*}

This proves Claim~\ref{CConsHomIII}.

This completes the proof.
\end{proof}

In the following proposition we prove that every P{\l}onka $\Sigma^{\neq 0}$-algebra induces an inductive system of $\Sigma^{\neq 0}$-algebras.

\begin{proposition}\label{PD4IndSys} Let $(\mathbf{A},D)$ be a P{\l}onka $\Sigma^{\neq 0}$-algebra, $\Phi^{D}$ the congruence induced by $D$ on $(\mathbf{A},D)$ and $\mathbf{Sl}(A,D) = (A/{\Phi^{D}}, \leq^{D})$ the sup-semilattice defined in Proposition~\ref{PUnivLnb}. Then the ordered pair 
\[
\mathcal{C}(\mathbf{A},D)=
(
(
\boldsymbol{\mathfrak{X}}
)_{\mathfrak{X}\in \mathrm{Sl}(A,D)},
(
f_{\mathfrak{X},\mathfrak{Y}}
)_{(\mathfrak{X},\mathfrak{Y})\in \leq^{D}}
)
\]
is an $\mathbf{Sl}(A,D)$-inductive system of $\Sigma^{\neq 0}$-algebras.
\end{proposition}
\begin{proof}
Le us note that, by Proposition~\ref{PClos} and Definition~\ref{DConsSub},  
$(\boldsymbol{\mathfrak{X}})_{\mathfrak{X}\in A/{\Phi^{D}}}$ is an $A/{\Phi^{D}}$-indexed family of 
$\Sigma^{\neq 0}$-algebras and that, by Proposition~\ref{PConsHom},  
$(f_{\mathfrak{X},\mathfrak{Y}})_{(\mathfrak{X},\mathfrak{Y})\in \leq^{D}}$ is a family of $\Sigma^{\neq 0}$-homomorphisms in $\prod_{(\mathfrak{X},\mathfrak{Y})\in \leq^{D}}\mathrm{Hom}(\boldsymbol{\mathfrak{X}},\boldsymbol{\mathfrak{Y}})$.

Therefore, by Definition~\ref{DIndSys}, it only remains to prove the following two claims.

\begin{claim}\label{CD4IndSysI} For every $\mathfrak{X}\in A/{\Phi^{D}}$, the equality $f_{\mathfrak{X},\mathfrak{X}}=\mathrm{id}_{\boldsymbol{\mathfrak{X}}}$ holds.
\end{claim}
Assume that $\mathfrak{X}=[x]_{\Phi^{D}}$. Then, by Proposition~\ref{PConsHom}, the mapping $f_{\mathfrak{X},\mathfrak{X}}$ is an endomapping of $[x]_{\Phi^{D}}$. Let $z\in [x]_{\Phi^{D}}$. Then the following chain of equalities holds
\allowdisplaybreaks
\begin{align*}
f_{\mathfrak{X},\mathfrak{X}}(z)&=D(z,x)
\tag{Prop.~\ref{PConsHom}}
\\&=
D(z,z)
\tag{Claim~\ref{CConsHomI}}
\\&=
z
\tag{by~\ref{D1}}
\\&=
\mathrm{id}_{[x]_{\Phi^{D}}}(z)
\tag{Def.~$\mathrm{id}_{[x]_{\Phi^{D}}}$}
\\&=
\mathrm{id}_{\boldsymbol{\mathfrak{X}}}(z).
\tag{Def.~$\mathfrak{X}$}
\end{align*}

This proves Claim~\ref{CD4IndSysI}.

\begin{claim}\label{CD4IndSysII} For every $\mathfrak{X}, \mathfrak{Y}$ and $\mathfrak{Z}$ in $A/{\Phi^{D}}$ if $(\mathfrak{X}, \mathfrak{Y})\in \leq^{D}$ and $(\mathfrak{Y}, \mathfrak{Z})\in \leq^{D}$, then $f_{\mathfrak{Y},\mathfrak{Z}}\circ f_{\mathfrak{X},\mathfrak{Y}}= f_{\mathfrak{X},\mathfrak{Z}}$.
\end{claim}
Assume that $\mathfrak{X}=[x]_{\Phi^{D}}$, $\mathfrak{Y}=[y]_{\Phi^{D}}$ and $\mathfrak{Z}=[z]_{\Phi^{D}}$. Since $(\mathfrak{Y},\mathfrak{Z})\in \leq^{D}$ we have that, by Proposition~\ref{PUnivLnb}, there exists elements $y'\in [y]_{\Phi^{D}}$ and $z'\in [z]_{\Phi^{D}}$ such that $D(z',y')=z'$.

Let $t\in [x]_{\Phi^{D}}$. Then the following chain of equalities holds
\allowdisplaybreaks
\begin{align*}
f_{\mathfrak{Y},\mathfrak{Z}}(
f_{\mathfrak{X},\mathfrak{Y}}(
t
))
&=
D(D(t,y),z)
\tag{Prop.~\ref{PConsHom}}
\\&=
D\left(D(t,y'),z'\right)
\tag{Claim.~\ref{CConsHomI}}
\\&=
D(t, D(y',z'))
\tag{by~\ref{D2}}
\\&=
D(t, D(z',y'))
\tag{by~\ref{D3}}
\\&=
D(t, z')
\tag{$D(z',y')=z'$}
\\&=
D(t, z)
\tag{Claim.~\ref{CConsHomI}}
\\&=
f_{\mathfrak{X},\mathfrak{Z}}(t).
\tag{Prop.~\ref{PConsHom}}
\end{align*}

This proves Claim~\ref{CD4IndSysII}.

This completes the proof.
\end{proof}

We next prove that every morphism of P{\l}onka $\Sigma^{\neq 0}$-algebras induces a morphism in $\int^{\mathsf{Ssl}}\mathrm{Isys}_{\Sigma^{\neq 0}}$. In the following proposition we will make use of the functor $\mathrm{Sl}$ defined in Proposition~\ref{PUnivLnb}.

\begin{proposition}\label{PConsMorII} 
To each morphism of $\mbox{\sffamily{\upshape{P{\l}Alg}}}(\Sigma^{\neq 0})$ we can uniquely assign a morphism of $\int^{\mathsf{Ssl}}\mathrm{Isys}_{\Sigma^{\neq 0}}$.
\end{proposition}

\begin{proof} 
Let $h$ be a morphism in $\mbox{\sffamily{\upshape{P{\l}Alg}}}(\Sigma^{\neq 0})$ from $(\mathbf{A},D)$ to $(\mathbf{B},E)$. By Definition~\ref{DIntSsl}, a morphism in $\int^{\mathsf{Ssl}}\mathrm{Isys}_{\Sigma^{\neq 0}}$ from $(\mathbf{Sl}(A,D),\mathcal{C}(\mathbf{A},D))$ to $(\mathbf{Sl}(B,E),\mathcal{C}(\mathbf{B},E))$ is given by a morphism from $\mathbf{Sl}(A,D)$ to $\mathbf{Sl}(B,E)$, which in this case is $\mathrm{Sl}(h)$, and a morphism $t(h)$ from the $\mathbf{Sl}(A,D)$-inductive system of $\Sigma^{\neq 0}$-algebras $\mathcal{C}(\mathbf{A},D)$ to the $\mathbf{Sl}(A,D)$-inductive system of $\Sigma^{\neq 0}$-algebras $\mathcal{C}(\mathbf{B},E)_{\mathrm{Sl}(h)}$.

That $\mathrm{Sl}(h)$ is a morphism in $\mathsf{Ssl}$ from $\mathbf{Sl}(A,D)$ to $\mathbf{Sl}(B,E)$ has already been proven in Proposition~\ref{PUnivLnb}. 

Let $\mathcal{C}(\mathbf{A},D)$ be the $\mathbf{Sl}(A,D)$-inductive system of $\Sigma^{\neq 0}$-algebras 
\[
\mathcal{C}(\mathbf{A},D)=
(
(
\boldsymbol{\mathfrak{X}}
)_{\mathfrak{X}\in\mathrm{Sl}(A,D)},
(
f_{\mathfrak{X},\mathfrak{Y}}
)_{(\mathfrak{X},\mathfrak{Y})\in \leq^{D}}
),
\]
and $\mathcal{C}(\mathbf{B},E)$ the $\mathbf{Sl}(B,E)$-inductive system of $\Sigma^{\neq 0}$-algebras
\[
\mathcal{C}(\mathbf{B},E)=
(
(
\boldsymbol{\mathfrak{U}}
)_{\mathfrak{U}\in\mathrm{Sl}(B,E)},
(
g_{\mathfrak{U},\mathfrak{V}}
)_{(\mathfrak{U},\mathfrak{V})\in \leq^{E}}
).
\]
Then, by Definition~\ref{DIsys}, we have that $\mathcal{C}(\mathbf{B},E)_{\mathrm{Sl}(h)}$ is the 
$\mathbf{Sl}(A,D)$-inductive system of $\Sigma^{\neq 0}$-algebras 
\[
\mathcal{C}(\mathbf{B},E)_{\mathrm{Sl}(h)}=
(
(
\mathrm{Sl}(h)(
\boldsymbol{\mathfrak{X}}
)
)_{\mathfrak{X}\in\mathrm{Sl}(\mathbf{A},D)},
(
g_{\mathrm{Sl}(h)(\mathfrak{X}),\mathrm{Sl}(h)(\mathfrak{Y})}
)_{(\mathfrak{X},\mathfrak{Y})\in \leq^{D}}
).
\]
Now, by Definition~\ref{DIndSys}, to give a morphism of $\mathbf{Sl}(A,D)$-inductive systems from
$\mathcal{C}(\mathbf{A},D)$ to $\mathcal{C}(\mathbf{B},E)_{\mathrm{Sl}(h)}$ is to give a family $t(h)=(t(h)_{\mathfrak{X}})_{\mathfrak{X}\in\mathrm{Sl}(\mathbf{A},D)}$ of $\Sigma^{\neq 0}$-homomor\-phisms in 
\[
\textstyle
\prod_{\mathfrak{X}\in\mathrm{Sl}(\mathbf{A},D)}\mathrm{Hom}(
\boldsymbol{\mathfrak{X}}, 
\mathrm{Sl}(h)(\boldsymbol{\mathfrak{X}}
))\]
such that, for every $(\mathfrak{X},\mathfrak{Y})\in\leq^{D}$, we have that 
\begin{align*}
t(h)_{\mathfrak{Y}}\circ f_{\mathfrak{X},\mathfrak{Y}}
&=
g_{\mathrm{Sl}(h)(\mathfrak{X}),\mathrm{Sl}(h)(\mathfrak{Y})}\circ t(h)_{\mathfrak{X}}.
\tag{A}\label{EqA}
\end{align*}

Let $\mathfrak{X}=[x]_{\Phi^{D}}$ be an equivalence class in $A/{\Phi^{D}}$. Then, by Proposition~\ref{PUnivLnb}, we have that 
$
\mathrm{Sl}(h)(\boldsymbol{\mathfrak{X}})
=[h(x)]_{\Phi^{E}}.
$
Let $t(h)_{\mathfrak{X}}$ be the correspondence from $[x]_{\Phi^{D}}$ to $[h(x)]_{\Phi^{E}}$ defined  as follows:
 \[
 t(h)_{\mathfrak{X}}
 \left\lbrace
\begin{array}{ccc}
[x]_{\Phi^{D}}&\mor &[h(x)]_{\Phi^{E}}
\\
z&\longmapsto &h(z)
\end{array}
\right.
\]
We need to check that $t(h)_{\mathfrak{X}}$ is a well-defined $\Sigma^{\neq 0}$-homomorphism from 
$\mathfrak{X}$ to $\mathrm{Sl}(h)(\boldsymbol{\mathfrak{X}})$.

\begin{claim}\label{CConsMorIIa} For $\mathfrak{X}=[x]_{\Phi^{D}}$, if $z\in [x]_{\Phi^{D}}$, then $h(z)\in [h(x)]_{\Phi^{E}}$.
\end{claim}

By Definition~\ref{DDORel}, to prove that $h(z)\in [h(x)]_{\Phi^{E}}$, under the hypothesis that $z\in [x]_{\Phi^{D}}$, is equivalent to prove that 
$$
E(h(x), h(z))= h(x) \quad\text{and}\quad E(h(z), h(x))= h(z).
$$

The following chain of equalities holds
\begin{align*}
E(h(x), h(z))&=
h(D(x,z))
\tag{$E\circ (h\times h)=h\circ D$}
\\&=
h(x).
\tag{$(x,z)\in \Phi^{D}$}
\end{align*}

By a similar argument we obtain that 
$$
E(h(z), h(x)) = h(z).
$$

This proves Claim~\ref{CConsMorIIa}.

\begin{claim}\label{CConsMorIIb} For $\mathfrak{X}=[x]_{\Phi^{D}}$ we have that $t(h)_{\mathfrak{X}}$ is a $\Sigma^{\neq 0}$-homomorphism from $\boldsymbol{\mathfrak{X}}$ to $ \mathrm{Sl}(h)(
\boldsymbol{\mathfrak{X}})$.
\end{claim}

Let $n$ be an element of $\mathbb{N}-1$, $\sigma\in \Sigma_{n}$ and 
$(z_{j})_{j\in n} \in \mathfrak{X}^{n}$. We need to check that 
\[t(h)_{\mathfrak{X}}(
F^{\boldsymbol{\mathfrak{X}}}_{\sigma}((z_{j})_{j\in n})
)
=
F^{\mathrm{Sl}(h)(
\boldsymbol{\mathfrak{X}}
)}_{\sigma}((
t(h)_{\mathfrak{X}}(
z_{j}
)
)_{j\in n}).
\]
The following chain of equalities holds
\allowdisplaybreaks
\begin{align*}
t(h)_{\mathfrak{X}}(
F^{\boldsymbol{\mathfrak{X}}}_{\sigma}((z_{j})_{j\in n})
)&=
h(
F^{\boldsymbol{\mathfrak{X}}}_{\sigma}((z_{j})_{j\in n})
)
\tag{Def.~$t(h)_{\mathfrak{X}}$}
\\&=
h(
F^{\mathbf{A}}_{\sigma}((z_{j})_{j\in n})
)
\tag{Def.~\ref{DConsSub}}
\\&=
F^{\mathbf{B}}_{\sigma}((
h(z_{j})
)_{j\in n})
\tag{$h$ is $\Sigma$-hom.}
\\&=
F^{\mathrm{Sl}(h)(
\boldsymbol{\mathfrak{X}}
)}_{\sigma}((
h(z_{j})
)_{j\in n})
\tag{Def.~\ref{DConsSub}}
\\&=
F^{\mathrm{Sl}(h)(
\boldsymbol{\mathfrak{X}}
)}_{\sigma}((
t(h)_{\mathfrak{X}}(
z_{j}
)
)_{j\in n}).
\tag{Def.~$t(h)_{\mathfrak{X}}$}
\end{align*}

This proves Claim~\ref{CConsMorIIb}.

It remains to prove Equation~\ref{EqA}. Let $(\mathfrak{X},\mathfrak{Y})$ be an element of $\leq^{D}$ and $z$ an element of $[x]_{\Phi^{D}}$. Then the following chain of equalities holds
\allowdisplaybreaks
\begin{align*}
t(h)_{\mathfrak{Y}}(f_{\mathfrak{X},\mathfrak{Y}}(
z))&=
t(h)_{\mathfrak{Y}}(D(z,y))
\tag{Prop.~\ref{PConsHom}}
\\&=
h(
D(z,y)
)
\tag{Def.~$t(h)_{\mathfrak{Y}}$}
\\&=
E(h(z),h(y))
\tag{$E\circ (h\times h)=h\circ D$}
\\&=
g_{\mathrm{Sl}(h)(\mathfrak{X}),\mathrm{Sl}(h)(\mathfrak{Y})}(
h(z)
)
\tag{Prop.~\ref{PConsHom}}
\\&=
g_{\mathrm{Sl}(h)(\mathfrak{X}),\mathrm{Sl}(h)(\mathfrak{Y})}(
t(h)_{\mathfrak{X}}(z)
).
\tag{Def.~$t(h)_{\mathfrak{X}}$}
\end{align*}

All in all, we conclude that $t(h)_{\mathfrak{Y}}\circ f_{\mathfrak{X},\mathfrak{Y}}
= g_{\mathrm{Sl}(h)(\mathfrak{X}),\mathrm{Sl}(h)(\mathfrak{Y})}\circ t(h)_{\mathfrak{X}}.$

This completes the proof.
\end{proof}

We next show that there exists a functor from $\mbox{\sffamily{\upshape{P{\l}Alg}}}(\Sigma^{\neq 0})$ to $\int^{\mathsf{Ssl}}\mathrm{Isys}_{\Sigma^{\neq 0}}$.

\begin{definition}\label{DIS}
We denote by $\mathrm{Is}_{\Sigma^{\neq 0}}$ the assignment from $\mbox{\sffamily{\upshape{P{\l}Alg}}}(\Sigma^{\neq 0})$ to $\int^{\mathsf{Ssl}}\mathrm{Isys}_{\Sigma^{\neq 0}}$ defined as follows:
\begin{enumerate}
\item for every P{\l}onka $\Sigma$-algebra $(\mathbf{A},D)$ in $\mbox{\sffamily{\upshape{P{\l}Alg}}}(\Sigma^{\neq 0})$, $\mathrm{Is}_{\Sigma^{\neq 0}}(\mathbf{A},D)$ is the ordered pair 
$
(\mathbf{Sl}(A,D),\mathcal{C}(\mathbf{A},D))
$, 
where, by Proposition~\ref{PUnivLnb}, $\mathbf{Sl}(A,D)$ is the sup-semilattice
$
(A/\Phi^{D},\leq^{D}),
$ 
and, by Proposition~\ref{PD4IndSys}, $\mathcal{C}(\mathbf{A},D)$ is the $\mathbf{Sl}(A,D)$-inductive system of $\Sigma^{\neq 0}$-algebras 
$$
((\boldsymbol{\mathfrak{X}})_{\mathfrak{X}\in\mathrm{Sl}(A,D)},(f_{\mathfrak{X},\mathfrak{Y}})_{(\mathfrak{X},\mathfrak{Y})\in \leq^{D}});
$$
\item for every morphism $h\colon (\mathbf{A},D)\mor (\mathbf{B},E)$ in  
$\mbox{\sffamily{\upshape{P{\l}Alg}}}(\Sigma^{\neq 0})$ between P{\l}onka $\Sigma$-algebras, $\mathrm{Is}_{\Sigma^{\neq 0}}(h)$ is the morphism in $\int^{\mathsf{Ssl}}\mathrm{Isys}_{\Sigma^{\neq 0}}$
\[
\textstyle
(\mathrm{Sl}(h),t(h))\colon 
(
\mathbf{Sl}(A,D),
\mathcal{C}(\mathbf{A},D)
)
\mor 
(
\mathbf{Sl}(B,E),
\mathcal{C}(\mathbf{B},E)
)
,
\]
introduced in Proposition~\ref{PConsMorII}.
\end{enumerate}
\end{definition}

\begin{proposition}\label{PIS}
The assignment $\mathrm{Is}_{\Sigma^{\neq 0}}$ is a functor from 
$\mbox{\sffamily{\upshape{P{\l}Alg}}}(\Sigma^{\neq 0})$ to $\int^{\mathsf{Ssl}}\mathrm{Isys}_{\Sigma^{\neq 0}}$.
\end{proposition}

\begin{proof}
That $\mathrm{Is}_{\Sigma^{\neq 0}}$ maps objects and morphisms of $\mbox{\sffamily{\upshape{P{\l}Alg}}}(\Sigma^{\neq 0})$ to objects and morphisms of $\int^{\mathsf{Ssl}}\mathrm{Isys}_{\Sigma^{\neq 0}}$ has been proved in Propositions~\ref{PD4IndSys} and~\ref{PConsMorII}, respectively. Therefore, all that remains is to check is that $\mathrm{Is}_{\Sigma}$ preserves identities and compositions.

\textsf{Preservation of identities.}

Let $(\mathbf{A},D)$ be an object in $\mbox{\sffamily{\upshape{P{\l}Alg}}}(\Sigma^{\neq 0})$. We need to prove that 
\[
\mathrm{Is}_{\Sigma^{\neq 0}}(\mathrm{id}_{(\mathbf{A},D)})=\mathrm{id}_{\mathrm{Is}_{\Sigma}(\mathbf{A},D)}.
\]
Let us recall, from Proposition~\ref{PConsMorII}, that 
\[
\mathrm{Is}_{\Sigma^{\neq 0}}(\mathrm{id}_{(\mathbf{A},D)})=(
\mathrm{Sl}(
\mathrm{id}_{(\mathbf{A},D)}
),
t(\mathrm{id}_{(\mathbf{A},D)})
)
\]

By Proposition~\ref{PUnivLnb}, $\mathrm{Sl}$ is a functor from $\mathsf{Lnb}$ to $\mathsf{Ssl}$. Thus, the following equation holds
\[
\mathrm{Sl}(
\mathrm{id}_{(\mathbf{A},D)}
)
=\mathrm{id}_{\mathrm{Sl}(\mathbf{A},D)}.
\]

Recall, from Proposition~\ref{PConsMorII}, that $t(\mathrm{id}_{(\mathbf{A},D)})=(t(\mathrm{id}_{(\mathbf{A},D)})_{\mathfrak{X}})_{\mathfrak{X}\in\mathrm{Sl}(\mathbf{A},D)}$. Let $\mathfrak{X}=[x]_{\Phi^{D}}$ be an equivalence class in $\mathrm{Sl}(\mathbf{A},D)=A/{\Phi^{D}}$. For an element $z\in [x]_{\Phi^{D}}$, the following chain of equalities holds
\begin{align*}
t(\mathrm{id}_{(\mathbf{A},D)})_{\mathfrak{X}}(z)
&=\mathrm{id}_{(\mathbf{A},D)}(z)
\tag{Prop.~\ref{PConsMorII}}
\\&=
\mathrm{id}_{\mathbf{A}}(z)
\tag{Prop.~\ref{DDOD4}}
\\&=
z
\tag{Def.~$\mathrm{id}_{\mathbf{A}}$}
\\&=
\mathrm{id}_{\boldsymbol{\mathfrak{X}}}(z).
\tag{Def.~$\mathrm{id}_{\boldsymbol{\mathfrak{X}}}$}
\end{align*}
Therefore, we conclude that $(t(\mathrm{id}_{(\mathbf{A},D)})_{\mathfrak{X}})_{\mathfrak{X}\in\mathrm{Sl}(\mathbf{A},D)}=(\mathrm{id}_{\boldsymbol{\mathfrak{X}}})_{\mathfrak{X}\in\mathrm{Sl}(\mathbf{A},D)}$.

All in all, we have that 
\begin{align*}
\mathrm{Is}_{\Sigma^{\neq 0}}(\mathrm{id}_{(\mathbf{A},D)})&=
(
\mathrm{Sl}(
\mathrm{id}_{(\mathbf{A},D)}
),
t(\mathrm{id}_{(\mathbf{A},D)})
)
\\&=
(
\mathrm{id}_{\mathbf{Sl}(A,D)},
(\mathrm{id}_{\boldsymbol{\mathfrak{X}}})_{\mathfrak{X}\in\mathrm{Sl}(\mathbf{A},D)}
)
\\&=\mathrm{id}_{\mathrm{Is}_{\Sigma^{\neq 0}}(\mathbf{A},D)}.
\end{align*}

This proves that $\mathrm{Is}_{\Sigma^{\neq 0}}$ preserves identities.

\textsf{Preservation of compositions.}

Let $g\colon (\mathbf{A},D)\mor (\mathbf{B},E)$ and $h\colon (\mathbf{B},E)\mor (\mathbf{C},F)$ be  morphisms in $\mbox{\sffamily{\upshape{P{\l}Alg}}}(\Sigma^{\neq 0})$. We need to prove that 
\[
\mathrm{Is}_{\Sigma^{\neq 0}}(h\circ g)=\mathrm{Is}_{\Sigma^{\neq 0}}(h)\circ \mathrm{Is}_{\Sigma^{\neq 0}}(g).
\]
Let us recall, from Proposition~\ref{PConsMorII}, that 
\[
\mathrm{Is}_{\Sigma^{\neq 0}}\left(h\circ g\right)=\left(
\mathrm{Sl}\left(
h\circ g
\right),
t\left(h\circ g\right)
\right)
\]

By Proposition~\ref{PUnivLnb}, $\mathrm{Sl}$ is a functor from $\mathsf{Lnb}$ to $\mathsf{Ssl}$. Thus, the following equation holds
\[
\mathrm{Sl}\left(
h\circ g
\right)
=\mathrm{Sl}\left(
h
\right)
\circ 
\mathrm{Sl}\left(
g
\right).
\]
Let us recall, from Proposition~\ref{PConsMorII}, that 
$t(h\circ g)=(t(h\circ g)_{\mathfrak{X}})_{\mathfrak{X}\in\mathrm{Sl}(A,D)}$. Fix an equivalence class $\mathfrak{X}=[x]_{\Phi^{D}}$ in $\mathrm{Sl}(A,D)=A/{\Phi^{D}}$. For an element $z\in [x]_{\Phi^{D}}$, the following chain of equalities holds
\allowdisplaybreaks
\begin{align*}
t(h\circ g)_{\mathfrak{X}}(z)
&=h(g(z))
\tag{Prop.~\ref{PConsMorII}}
\\&=
h(
t(g)_{\mathfrak{X}}(z)
)
\tag{Prop.~\ref{PConsMorII}}
\\&=
t(h)_{\mathrm{Sl}(g)(\mathfrak{X})}(
t(g)_{\mathfrak{X}}(z)
).
\tag{Prop.~\ref{PConsMorII}}
\end{align*}
Therefore $(t(h\circ g)_{\mathfrak{X}})_{\mathfrak{X}\in\mathrm{Sl}(A,D)}=(
t(h)_{\mathrm{Sl}(g)(\mathfrak{X})}\circ
t(g)_{\mathfrak{X}}
)_{\mathfrak{X}\in\mathrm{Sl}(A,D)}$.

All in all, we have that 
\begin{align*}
\mathrm{Is}_{\Sigma^{\neq 0}}(h\circ g)&=
(
\mathrm{Sl}(
h\circ g
),
t(h\circ g)
)
\\&=
(
\mathrm{Sl}(h)\circ \mathrm{Sl}(g),
(
t(h)_{\mathrm{Sl}(g)(\mathfrak{X})}\circ
t(g)_{\mathfrak{X}}
)_{\mathfrak{X}\in\mathrm{Sl}(A,D)}
)
\\&=
\mathrm{Is}_{\Sigma^{\neq 0}}(h)\circ \mathrm{Is}_{\Sigma^{\neq 0}}(g).
\end{align*}

This proves that $\mathrm{Is}_{\Sigma^{\neq 0}}$ preserves compositions.

This completes the proof.
\end{proof}

%
%
%


We next prove that the functor $\mathrm{Is}_{\Sigma^{\neq 0}}$ has a left adjoint.

\begin{proposition}\label{LAdj}
The functor $\mathrm{Is}_{\Sigma^{\neq 0}}$ from $\mbox{\sffamily{\upshape{P{\l}Alg}}}(\Sigma^{\neq 0})$ to $\int^{\mathsf{Ssl}}\mathrm{Isys}_{\Sigma^{\neq 0}}$ has a left adjoint.
\end{proposition}

\begin{proof}
Let $\mathbcal{A} = (\mathbf{I},\mathcal{A})$ be an inductive system of $\Sigma^{\neq 0}$-algebras, where the $\mathbf{I}$-inductive system of $\Sigma^{\neq 0}$-algebras $\mathcal{A}$ is $((\mathbf{A}_{i})_{i\in I},(f_{i,j})_{(i,j)\in \leq})$. We will prove that there exists a universal morphism from $\mathbcal{A}$ to $\mathrm{Is}_{\Sigma^{\neq 0}}$, i.e., that there exists a P{\l}onka $\Sigma^{\neq 0}$-algebra 
$(\mbox{\bfseries{\upshape{P{\l}}}}_{\Sigma^{\neq 0}}(\mathbcal{A}), D^{\mathbcal{A}})$ and a morphism $\eta_{\mathbcal{A}}$ from $\mathbcal{A}$ to 
$\mathrm{Is}_{\Sigma^{\neq 0}}(\mbox{\bfseries{\upshape{P{\l}}}}_{\Sigma^{\neq 0}}(\mathbcal{A}), D^{\mathbcal{A}})$ such that, for every P{\l}onka $\Sigma^{\neq 0}$-algebra $(\mathbf{B},E)$ and every morphism $(\xi, u)$ from $\mathbcal{A}$ to $\mathrm{Is}_{\Sigma^{\neq 0}}(\mathbf{B},E)$, there exists a unique morphism $(\xi, u)^{\sharp}$ in $\mbox{\sffamily{\upshape{P{\l}Alg}}}(\Sigma^{\neq 0})$ 
\[
(\xi, u)^{\sharp}\colon 
\left(
\mbox{\bfseries{\upshape{P{\l}}}}_{\Sigma^{\neq 0}}(\mathbcal{A}), D^{\mathbcal{A}}
\right)
\mor 
(\mathbf{B},E)
\]
such that $(\xi, u)=\mathrm{Is}_{\Sigma^{\neq 0}}((\xi, u)^{\sharp})\circ \eta_{\mathbcal{A}}$.

\textsf{Definition of} $(\mbox{\bfseries{\upshape{P{\l}}}}_{\Sigma^{\neq 0}}(\mathbcal{A}), D^{\mathbcal{A}})$.

We begin by defining the $\Sigma^{\neq 0}$-algebra 
$\mbox{\bfseries{\upshape{P{\l}}}}_{\Sigma^{\neq 0}}(\mathbcal{A})$ as follows:
\begin{enumerate}
\item the underlying set of $\mbox{\bfseries{\upshape{P{\l}}}}_{\Sigma^{\neq 0}}(\mathbcal{A})$ is 
\[
\textstyle
\coprod_{i\in I}A_{i} = \bigcup_{i\in I}(A_{i}\times \{i\}),\] 
\item for every $n\in\mathbb{N}-1$ and every $\sigma\in \Sigma_{n}$, the structural operation $F^{\scriptsize{\mbox{\bfseries{\upshape{P{\l}}}}_{\Sigma^{\neq 0}}(\mathbcal{A})}}_{\sigma}$ is defined as
\[
\textstyle
F^{\scriptsize{\mbox{\bfseries{\upshape{P{\l}}}}_{\Sigma^{\neq 0}}(\mathbcal{A})}}_{\sigma}
\left\lbrace
\begin{array}{ccc}
(\bigcup_{i\in I}(A_{i}\times\{i\}))^{n}
&\mor&
\bigcup_{i\in I} (A_{i}\times\{i\})
\\
(x_{j},i_{j})_{j\in n}
&\longmapsto&
(
F^{\mathbf{A}_{\bigvee_{j\in n}i_{j}}}_{\sigma}(
(f_{i_{j},\bigvee_{j\in n}i_{j}}(
x_{j}
))_{j\in n}),
\bigvee_{j\in n}i_{j}
).
\end{array}
\right.
\]
\end{enumerate} 
We will call the $\Sigma^{\neq 0}$-algebra $\mbox{\bfseries{\upshape{P{\l}}}}_{\Sigma^{\neq 0}}(\mathbcal{A})$ the \emph{P{\l}onka sum} of $\mathbcal{A}$ and we will let $\mbox{\upshape{P{\l}}}_{\Sigma^{\neq 0}}(\mathbcal{A})$ stand for its underlying set.

We next define a binary operation on the underlying set of $\mbox{\bfseries{\upshape{P{\l}}}}_{\Sigma^{\neq 0}}(\mathbcal{A})$ and prove that it is a P{\l}onka operator for $\mbox{\bfseries{\upshape{P{\l}}}}_{\Sigma^{\neq 0}}(\mathbcal{A})$.

Let $D^{\mathbcal{A}}$ be the mapping defined as follows:
\[
D^{\mathbcal{A}}
\left\lbrace
\begin{array}{ccc}
 \bigcup_{i\in I}\left(A_{i}\times\{i\}\right)\times \bigcup_{i\in I}\left(A_{i}\times\{i\}\right)&\mor&\bigcup_{i\in I}\left(A_{i}\times\{i\}\right)
 \\
 ((x,j),(y,k))&\longmapsto&
 \left(
f_{j,j\vee k}(x),j\vee k
\right)
\end{array}
\right.
\]
Note that, for a pair $((x,j),(y,k))$, with $j,k\in I$, $x\in A_{j}$ and $y\in A_{k}$, the pair 
$D^{\mathbcal{A}}((x,j),(y,k))$ is an element of $A_{j\vee k}\times \{j\vee k\}\subseteq \bigcup_{i\in I}(A_{i}\times\{i\})$. 

We next prove that $D^{\mathbcal{A}}$ is a P{\l}onka operator for $\mbox{\bfseries{\upshape{P{\l}}}}_{\Sigma^{\neq 0}}(\mathbcal{A})$.

We start by showing that $D^{\mathbcal{A}}$  satisfies all the conditions stated in Definition~\ref{DDO}.

\begin{claim}\label{PD1} Let $(x,j)$ be an element of $\bigcup_{i\in I}(A_{i}\times \{i\})$. Then the following equality holds
\[D^{\mathbcal{A}}((x,j),(x,j))=(x,j).\]
\end{claim}
The following chain of equalities holds
\allowdisplaybreaks
\begin{align*}
D^{\mathbcal{A}}((x,j),(x,j))&=
(
f_{j,j\vee j}(x),j\vee j
)
\tag{Def.~$D^{\mathbcal{A}}$}
\\&=
(
f_{j,j}(x),j
)
\tag{Def.~\ref{DSup}}
\\&=
(
\mathrm{id}_{\mathbf{A}_{j}}(x),j
)
\tag{Def.~\ref{DIndSys}}
\\&=
(
x,j
).
\tag{Def.~$\mathrm{id}_{\mathbf{A}_{j}}$}
\end{align*}
This proves Claim~\ref{PD1}.

\begin{claim}\label{PD2} Let $(x,j)$ and $(y,k)$ be elements of $\bigcup_{i\in I}(A_{i}\times \{i\})$. Then  the following equality holds
\[D^{\mathbcal{A}}((x,j),D^{\mathbcal{A}}((y,k),(z,l)))=
D^{\mathbcal{A}}(D^{\mathbcal{A}}((x,j),(y,k)),(z,l))
.\]
\end{claim}
The following chain of equalities holds
\allowdisplaybreaks
\begin{align*}
D^{\mathbcal{A}}((x,j),D^{\mathbcal{A}}((y,k),(z,l)))
&=
D^{\mathbcal{A}}((x,j),
(f_{k,k\vee l}(y),k\vee l)
)
\tag{Def.~$D^{\mathbcal{A}}$}
\\&=
(
f_{j,j\vee (k\vee l)}(x), j\vee (k\vee l)
)
\tag{Def.~$D^{\mathbcal{A}}$}
\\&=
(
f_{j,j\vee k\vee l}(x), j\vee k\vee l
)
\tag{Def.~\ref{DSup}}
\\&=
(
f_{j,(j\vee k)\vee l}(x), (j\vee k)\vee l
)
\tag{Def.~\ref{DSup}}
\\&=
(
f_{j\vee k, (j\vee k)\vee l}(
f_{j,j\vee k}(x)
),
(j\vee k)\vee l
)
\tag{Def.~\ref{DIndSys}}
\\&=
D^{\mathbcal{A}}(
(f_{j,j\vee k}(x),j\vee k),(z,l)
)
\tag{Def.~$D^{\mathbcal{A}}$}
\\&=
D^{\mathbcal{A}}(
D^{\mathbcal{A}}((x,j),(y,k))
,(z,l)
).
\tag{Def.~$D^{\mathbcal{A}}$}
\end{align*}
This proves Claim~\ref{PD2}.

\begin{claim}\label{PD3} Let $(x,j)$, $(y,k)$ and $(z,l)$ elements of $\bigcup_{i\in I}(A_{i}\times \{i\})$. Then the following equality holds 
\[D^{\mathbcal{A}}((x,j),D^{\mathbcal{A}}((y,k),(z,l)))=
D^{\mathbcal{A}}((x,j),D^{\mathbcal{A}}((z,l),(y,k)))
.\]
\end{claim}
The following chain of equalities holds
\allowdisplaybreaks
\begin{align*}
D^{\mathbcal{A}}((x,j),D^{\mathbcal{A}}((y,k),(z,l)))
&=
D^{\mathbcal{A}}((x,j),
(f_{k,k\vee l}(y),k\vee l
)
)
\tag{Def.~$D^{\mathbcal{A}}$}
\\&=
(
f_{j,j\vee(k\vee l)}(x), j\vee (k\vee l)
)
\tag{Def.~$D^{\mathbcal{A}}$}
\\&=
(
f_{j,j\vee k\vee l}(x), j\vee k\vee l
)
\tag{Def.~\ref{DSup}}
\\&=
(
f_{j,j\vee (l\vee k)}(x), j\vee (l\vee k)
)
\tag{Def.~\ref{DSup}}
\\&=
D^{\mathbcal{A}}((x,j),
(
f_{l,l\vee k}(z), l\vee k
)
)
\tag{Def.~$D^{\mathbcal{A}}$}
\\&=
D^{\mathbcal{A}}((x,j),D^{\mathbcal{A}}((z,l),(y,k))).
\tag{Def.~$D^{\mathbcal{A}}$}
\end{align*}
This proves Claim~\ref{PD3}.

This proves that $(\coprod_{i\in I}A_{i}, D^{\mathbcal{A}})$ is a left normal band. Now, we check that it  also satisfies the two conditions set forth in Definition~\ref{DDOD4}.

\begin{claim}\label{PD4} Let $n$ be an element of $\mathbb{N}-1$, $\sigma\in \Sigma_{n}$, 
$(x_{j},i_{j})_{j\in n}\in(\bigcup_{i\in I}(A_{i}\times \{i\}))^{n}$ and $(y,k)\in\bigcup_{i\in I}(A_{i}\times \{i\})$. Then the following equality holds
\[
D^{\mathbcal{A}}(
F^{\scriptsize{\mbox{\bfseries{\upshape{P{\l}}}}_{\Sigma^{\neq 0}}(\mathbcal{A})}}_{\sigma}((x_{j},i_{j})_{j\in n}),
(y,k))=
F^{\scriptsize{\mbox{\bfseries{\upshape{P{\l}}}}_{\Sigma^{\neq 0}}(\mathbcal{A})}}_{\sigma}(
(
D^{\mathbcal{A}}(
(x_{j},i_{j}),(y,k)
)
)_{j\in n}
).
\]
\end{claim}
The following chain of equalities holds
\begin{flushleft}
$D^{\mathbcal{A}}(
F^{\scriptsize{\mbox{\bfseries{\upshape{P{\l}}}}_{\Sigma^{\neq 0}}(\mathbcal{A})}}_{\sigma}((x_{j},i_{j})_{j\in n}),
(y,k))$
\allowdisplaybreaks
\begin{align*}
&=
\textstyle
D^{\mathbcal{A}}(
(
F^{\mathbf{A}_{\bigvee_{j\in n}i_{j}}}_{\sigma}((
f_{i_{j}, \bigvee_{j\in n}i_{j}}(x_{j})
)_{j\in n}
),
\bigvee_{j\in n}i_{j}
),
(y,k)
)
\tag{Def.~$\mbox{\bfseries{\upshape{P{\l}}}}_{\Sigma^{\neq 0}}(\mathbcal{A})$}
\\&=
\textstyle
(
f_{\bigvee_{j\in n} i_{j},(\bigvee_{j\in n} i_{j})\vee k}
(
F^{\mathbf{A}_{\bigvee_{j\in n}i_{j}}}_{\sigma}((
f_{i_{j}, \bigvee_{j\in n}i_{j}}(x_{j})
)_{j\in n}
)),
(\bigvee_{j\in n} i_{j})\vee k
)
\tag{Def.~$D^{\mathbcal{A}}$}
\\&=
\textstyle
(
F^{\mathbf{A}_{(\bigvee_{j\in n}i_{j})\vee k}}_{\sigma}(
(
f_{\bigvee_{j\in n} i_{j}, (\bigvee_{j\in n} i_{j})\vee k}(
f_{i_{j},\bigvee_{j\in n} i_{j}}(x_{j})
)
)_{j\in n}
),
\\&\qquad\qquad\qquad\qquad\qquad\qquad\qquad\qquad
\qquad\qquad\qquad
\textstyle
(\bigvee_{j\in n} i_{j})\vee k
)
\tag{$\Sigma^{\neq 0}$-\text{hom}.}
\\&=
\textstyle
(
F^{\mathbf{A}_{(\bigvee_{j\in n}i_{j})\vee k}}_{\sigma}(
(
f_{i_{j},(\bigvee_{j\in n} i_{j})\vee k}(x_{j})
)_{j\in n}
)
,
(\bigvee_{j\in n} i_{j})\vee k
)
\tag{Def.~\ref{DIndSys}}
\\&=
\textstyle
(
F^{\mathbf{A}_{\bigvee_{j\in n}(i_{j}\vee k)}}_{\sigma}(
(
f_{i_{j},\bigvee_{j\in n} (i_{j}\vee k)}(x_{j})
)_{j\in n}
)
,
\bigvee_{j\in n} (i_{j}\vee k)
)
\tag{Def.~\ref{DSup}}
\\&=
\textstyle
(
F^{\mathbf{A}_{\bigvee_{j\in n}(i_{j}\vee k)}}_{\sigma}(
(
f_{i_{j}\vee k,\bigvee_{j\in n} (i_{j}\vee k)}(
f_{i_{j},i_{j}\vee k}(x_{j}))
)_{j\in n}
)
,
\bigvee_{j\in n} (i_{j}\vee k)
)
\tag{Def.~\ref{DIndSys}}
\\&=
F^{\scriptsize{\mbox{\bfseries{\upshape{P{\l}}}}_{\Sigma^{\neq 0}}(\mathbcal{A})}}_{\sigma}(
(
(f_{i_{j},i_{j}\vee k}(x_{j}), i_{j}\vee k)
)_{j\in n}
)
\tag{Def.~$\mbox{\bfseries{\upshape{P{\l}}}}_{\Sigma^{\neq 0}}(\mathbcal{A})$}
\\&=
F^{\scriptsize{\mbox{\bfseries{\upshape{P{\l}}}}_{\Sigma^{\neq 0}}(\mathbcal{A})}}_{\sigma}(
(
D^{\mathbcal{A}}(
(x_{j},i_{j}),(y,k)
)
)_{j\in n}
).
\tag{Def.~$D^{\mathbcal{A}}$}
\end{align*}
\end{flushleft}
This proves Claim~\ref{PD4}.

It remains to prove that the mapping $D^{\mathbcal{A}}$ satisfies condition~\ref{D5}. But before proving it, we prove the following claim.

\begin{claim}\label{CIt} Let $n$ be an element of $\mathbb{N}-1$ and $(x_{j},i_{j})_{j\in n}\in (\bigcup_{i\in I} (A_{i}\times \{i\}))^{n}$. Then
\[
\textstyle
D^{\mathbcal{A}}_{n}(
(x_{j},i_{j})_{j\in n}
)=
(
f_{i_{0},\bigvee_{j\in n}i_{j}}(x_{0}), \bigvee_{j\in n}i_{j}
).
\]
\end{claim}
We prove it by induction on $n$. 

\textsf{Base case}. For $n=1$, the following equality holds
\begin{align*}
D^{\mathbcal{A}}_{1}(x_{0},i_{0})
&=\mathrm{id}_{\scriptsize{\mbox{\bfseries{\upshape{P{\l}}}}_{\Sigma^{\neq 0}}(\mathbcal{A})}}(x_{0},i_{0})
\tag{Def~\ref{DIt}}
\\&=
(x_{0},i_{0})
\tag{Def.~$\mathrm{id}_{\scriptsize{\mbox{\bfseries{\upshape{P{\l}}}}_{\Sigma^{\neq 0}}(\mathbcal{A})}}$}
\\&=
(f_{i_{0},i_{0}}(x_{0}), i_{0}).
\tag{Def.~\ref{DIndSys}}
\end{align*}

This proves the base case.

\textsf{Inductive step}. Assume that the statement holds for sequences of length $n-1$. Let us prove it for sequences of length $n$. The following chain of equalities holds
\allowdisplaybreaks
\begin{align*}
D^{\mathbcal{A}}_{n}(
(x_{j},i_{j})_{j\in n}
)&=
D^{\mathbcal{A}}_{n}((x_{0},i_{0}),
D^{\mathbcal{A}}_{n-1}(
(x_{j},i_{j})_{j\in [1,n-1]}
))
\tag{Def.~\ref{DIt}}
\\&=
\textstyle
D^{\mathbcal{A}}((x_{0},i_{0}),
(
f_{i_{1},\bigvee_{j\in [1,n-1]}i_{j}}(x_{1}), \bigvee_{j\in [1,n-1]}i_{j}
)
)
\tag{Induction}
\\&=
\textstyle
(
f_{i_{0},i_{0}\vee \bigvee_{j\in [1,n-1]}i_{j}}(x_{0}),
i_{0}\vee \bigvee_{j\in [1,n-1]}i_{j}
)
\tag{Def.~$D^{\mathbcal{A}}$}
\\&=
\textstyle
(
f_{i_{0},\bigvee_{j\in n}i_{j}}(x_{0}),
\bigvee_{j\in n}i_{j}
).
\tag{Def.~\ref{DSup}}
\end{align*}

This proves Claim~\ref{CIt}.

We are now in position to prove that Condition~\ref{D5} holds.

\begin{claim}\label{PD5} Let $n$ be an element of $\mathbb{N}-1$, $\sigma\in \Sigma_{n}$, $(x_{j},i_{j})_{j\in n} \in (\bigcup_{i\in I}(A_{i}\times \{i\}))^{n}$ and $(y,k)\in \bigcup_{i\in I}(A_{i}\times \{i\})$. Then the following equality holds
\[
D^{\mathbcal{A}}(
(y,k),
F^{\scriptsize{\mbox{\bfseries{\upshape{P{\l}}}}_{\Sigma^{\neq 0}}(\mathbcal{A})}}_{\sigma}((x_{j},i_{j})_{j\in n}))=
D^{\mathbcal{A}}(
(y,k),
D^{\mathbcal{A}}_{n}((x_{j},i_{j})_{j\in n})).
\]
\end{claim}

The following chain of equalities holds
\begin{flushleft}
$D^{\mathbcal{A}}(
(y,k),
F^{\scriptsize{\mbox{\bfseries{\upshape{P{\l}}}}_{\Sigma^{\neq 0}}(\mathbcal{A})}}_{\sigma}((x_{j},i_{j})_{j\in n}))$
\allowdisplaybreaks
\begin{align*}
&=
\textstyle
D^{\mathbcal{A}}(
(y,k),
(
F^{\mathbf{A}_{\bigvee_{j\in n} i_{j}}}_{\sigma}((
f_{i_{j}, \bigvee_{j\in n} i_{j}}(x_{j})
)_{j\in n})
,
\bigvee_{j\in n}i_{j}
)
)
\tag{Def.~$\mbox{\bfseries{\upshape{P{\l}}}}_{\Sigma^{\neq 0}}(\mathbcal{A})$}
\\&=
\textstyle
(f_{k, k\vee (\bigvee_{j\in n} i_{j})}(y),
k\vee (\bigvee_{j\in n} i_{j})
)
\tag{Def.~$D^{\mathbcal{A}}$}
\\&=
\textstyle
D^{\mathbcal{A}}(
(y,k),
(
f_{i_{0},\bigvee_{j\in n}i_{j}}(x_{0}), \bigvee_{j\in n}i_{j}
)
)
\tag{Def.~$D^{\mathbcal{A}}$}
\\&=
D^{\mathbcal{A}}(
(y,k),
D^{\mathbcal{A}}_{n}((x_{j},i_{j})_{j\in n})).
\tag{Claim~\ref{CIt}}
\end{align*}
\end{flushleft}

This proves Claim~\ref{PD5}.

This completes the proof that $D^{\mathbcal{A}}$ is a P{\l}onka operator for 
$\mbox{\bfseries{\upshape{P{\l}}}}_{\Sigma^{\neq 0}}(\mathbcal{A})$. Therefore 
$(\mbox{\bfseries{\upshape{P{\l}}}}_{\Sigma^{\neq 0}}(\mathbcal{A}),D^{\mathbcal{A}})$ is a P{\l}onka $\Sigma^{\neq 0}$-algebra.

For the  P{\l}onka $\Sigma^{\neq 0}$-algebra $(\mbox{\bfseries{\upshape{P{\l}}}}_{\Sigma^{\neq 0}}(\mathbcal{A}),D^{\mathbcal{A}})$ we next describe the congruence induced by 
$D^{\mathbcal{A}}$ $\Phi^{D^{\mathbcal{A}}}$ on 
$(\mbox{\bfseries{\upshape{P{\l}}}}_{\Sigma^{\neq 0}}(\mathbcal{A}),D^{\mathbcal{A}})$ and the order $\leq^{D^{\mathbcal{A}}}$ associated to $D^{\mathbcal{A}}$ on 
$\coprod_{i\in I}A_{i}/\Phi^{D^{\mathbcal{A}}}$.

\begin{claim}\label{cgr}
For every pair $(x,j)$, $(y,k)$ in $\bigcup_{i\in I}(A_{i}\times \{i\})$, it holds that 
\[
((x,j),(y,k))\in \Phi^{D^{\mathbcal{A}}}\quad\mbox{if, and only if,}\quad j=k.
\]
\end{claim}

Let $(x,j)$ and $(y,k)$ be elements of $\bigcup_{i\in I}(A_{i}\times \{i\})$.  Let us assume that $((x,j),(y,k))\in \Phi^{D^{\mathbcal{A}}}$. Then, by Definition~\ref{DDORel}, this is equivalent to 
$$
D^{\mathbcal{A}}((x,j),(y,k))= (x,j) \quad\text{and}\quad D^{\mathbcal{A}}((y,k),(x,j))= (y,k),
$$
which, in turn, unpacking the definition of $D^{\mathbcal{A}}$, is equivalent to  
$$
(f_{j,j\vee k}(x),j\vee k)= (x,j)  \quad\text{and}\quad (f_{k,k\vee j}(y),k\vee j)= (y,k).
$$
From these two equations, we conclude that $j\vee k=j$ and $k\vee j=k$, i.e., that $k\leq j$ and $j\leq k$. By the antisymmetry of the order $\leq$ in $\mathbf{I}$, we conclude that $j=k$.

Reciprocally, let us assume that $j=k$. Then, to prove that $((x,j),(y,k))\in \Phi^{D^{\mathbcal{A}}}$, it suffices to check that $D^{\mathbcal{A}}((x,j),(y,k))=(x,j)$ and 
$D^{\mathbcal{A}}((y,k),(x,j))=(y,k)$.

Note that the following chain of equalities holds
\allowdisplaybreaks
\begin{align*}
D^{\mathbcal{A}}((x,j),(y,k))&=
(f_{j,j\vee k}(x),j\vee k)
\tag{Def.~$D^{\mathbcal{A}}$}
\\&=
(f_{j,j\vee j}(x),j\vee j)
\tag{$j=k$}
\\&=
(f_{j,j}(x),j)
\tag{Def.~\ref{DSup}}
\\&=
(
\mathrm{id}_{\mathbf{A}_{j}}(x),j
)
\tag{Def.~\ref{DIndSys}}
\\&=
(x,j).
\tag{Def.~$\mathrm{id}_{\mathbf{A}_{j}}$}
\end{align*}

By a similar argument we obtain that
$$
D^{\mathbcal{A}}((y,k), (x,j)) = (y,k).
$$
We conclude that $((x,j),(y,k))\in \Phi^{D^{\mathbcal{A}}}$

This proves Claim~\ref{cgr}.

\begin{claim}\label{ord}
For every pair $(x,j)$, $(y,k)$ in $\bigcup_{i\in I}(A_{i}\times \{i\})$, it holds that 
\[
[(x,j)]_{\Phi^{D^{\mathbcal{A}}}}
\leq^{D^{\mathbcal{A}}}
[(y,k)]_{\Phi^{D^{\mathbcal{A}}}}
\quad\mbox{if, and only if,}\quad j\leq k.
\]
\end{claim}

Let $(x,j)$ and $(y,k)$ be elements of $\bigcup_{i\in I}(A_{i}\times \{i\})$ such that 
\[
[(x,j)]_{\Phi^{D^{\mathbcal{A}}}}
\leq^{D^{\mathbcal{A}}}
[(y,k)]_{\Phi^{D^{\mathbcal{A}}}}
\]
Then, by Proposition~\ref{PUnivLnb}, there exists elements $(x',j)\in [(x,j)]_{\Phi^{D^{\mathbcal{A}}}}$ and $(y',k)\in [(y,k)]_{\Phi^{D^{\mathbcal{A}}}}$ such that 
$D^{\mathbcal{A}}((y',k),(x',j))=(y',k)$. Unpacking the definition of $D^{\mathbcal{A}}$, we have that 
$(
f_{k,k\vee j}(y'), k\vee j)
=(y',k)$.
From this equation it follows that $k\vee j=k$, i.e., that $j\leq k$.

Reciprocally, let $(x,j)$ and $(y,k)$ be of elements of $\bigcup_{i\in I}(A_{i}\times \{i\})$ such that 
$j\leq k$. Note that the following chain of equalities holds
\allowdisplaybreaks
\begin{align*}
D^{\mathbcal{A}}((y,k),(x,j))&=
(
f_{k,k\vee j}(y), k\vee j
)
\tag{Def.~$D^{\mathbcal{A}}$}
\\&=
(
f_{k,k}(y), k
)
\tag{$j\leq k$}
\\&=
(
\mathrm{id}_{\mathbf{A}_{k}}(y),k
)
\tag{Def.~\ref{DIndSys}}
\\&=
(
y,k
).
\tag{Def.~$\mathrm{id}_{\mathbf{A}_{k}}$}
\end{align*}
Then, by Proposition~\ref{PUnivLnb}, we conclude that 
$
[(x,j)]_{\Phi^{D^{\mathbcal{A}}}}
\leq^{D^{\mathbcal{A}}}
[(y,k)]_{\Phi^{D^{\mathbcal{A}}}}
$.

This proves Claim~\ref{ord}.

\textsf{Definition of the morphism} $\eta_{\mathbcal{A}}$.

We next prove that there exists a morphism $\eta_{\mathbcal{A}}$ from $\mathbcal{A}$ to 
$\mathrm{Is}_{\Sigma^{\neq 0}}(\mbox{\bfseries{\upshape{P{\l}}}}_{\Sigma^{\neq 0}}(\mathbcal{A}), D^{\mathbcal{A}})$. But before doing so, we recall that, by Definition~\ref{DIS}, $\mathrm{Is}_{\Sigma^{\neq 0}}(\mbox{\bfseries{\upshape{P{\l}}}}_{\Sigma^{\neq 0}}(\mathbcal{A}), D^{\mathbcal{A}})$ is the ordered pair
\[
(\mathbf{Sl}(\mbox{\upshape{P{\l}}}_{\Sigma^{\neq 0}}(\mathbcal{A}), D^{\mathbcal{A}}),\mathcal{C}(\mbox{\bfseries{\upshape{P{\l}}}}_{\Sigma^{\neq 0}}(\mathbcal{A}), D^{\boldsymbol{\mathcal{A}}})),
\]
where, by Proposition~\ref{PUnivLnb}, 
$\mathbf{Sl}(\mbox{\upshape{P{\l}}}_{\Sigma^{\neq 0}}(\mathbcal{A}), D^{\boldsymbol{\mathcal{A}}})$ is the sup-semilattice 
\[
\textstyle
((\coprod_{i\in I}A_{i})/{\Phi^{D^{\mathbcal{A}}}}, \leq^{D^{\boldsymbol{\mathcal{A}}}})
\]
associated to the left normal band $(\mbox{\upshape{P{\l}}}_{\Sigma^{\neq 0}}(\mathbcal{A}), D^{\boldsymbol{\mathcal{A}}})$ and, by Proposition~\ref{PD4IndSys}, $\mathcal{C}(\mbox{\bfseries{\upshape{P{\l}}}}_{\Sigma^{\neq 0}}(\mathbcal{A}), D^{\boldsymbol{\mathcal{A}}})$ is the $\mathbf{Sl}(\mbox{\upshape{P{\l}}}_{\Sigma^{\neq 0}}(\mathbcal{A}), D^{\boldsymbol{\mathcal{A}}})$-inductive system of $\Sigma^{\neq 0}$-algebras
$$
((\boldsymbol{\mathfrak{X}})_{\mathfrak{X}\in\mathrm{Sl}(\mbox{\upshape{{\scriptsize P{\l}}}}_{\Sigma^{\neq 0}}(\boldsymbol{\mathcal{A}}), D^{\boldsymbol{\mathcal{A}}})},(g_{\mathfrak{X},\mathfrak{Y}
})_{(\mathfrak{X},\mathfrak{Y})\in \leq^{D^{\boldsymbol{\mathcal{A}}}}}).
$$

The morphism 
$
\eta_{\mathbcal{A}}\colon \mathbcal{A} \mor \mathcal{C}\left(
\mbox{\bfseries{\upshape{P{\l}}}}_{\Sigma^{\neq 0}}(\mathbcal{A}), D^{\mathbcal{A}}
\right)
$
is given by a pair $\eta_{\mathbcal{A}}=(\alpha_{\mathbcal{A}}, v_{\mathbcal{A}})$ in which $\alpha_{\mathbcal{A}}$ is a morphism  of sup-semilattices from $\mathbf{I}$ to $\mathbf{Sl}(\mbox{\upshape{P{\l}}}_{\Sigma^{\neq 0}}(\mathbcal{A}), D^{\boldsymbol{\mathcal{A}}})$ and $v_{\mathbcal{A}}$ a morphism from  $\mathbcal{A}$ to $\mathcal{C}
(\mbox{\bfseries{\upshape{P{\l}}}}_{\Sigma^{\neq 0}}(\boldsymbol{\mathcal{A}}), D^{\boldsymbol{\mathcal{A}}})_{\alpha_{\mathbcal{A}}}$



To obtain $\alpha_{\mathbcal{A}}$ we proceed as follows. Let $i\in I$ and let $z_{i}$ be any element of $A_{i}$. Then we define the mapping $\alpha_{\mathbcal{A}}$ from $I$ to $\mathrm{Sl}\left(
\mbox{\upshape{P{\l}}}_{\Sigma^{\neq 0}}(\mathbcal{A}), D^{\boldsymbol{\mathcal{A}}}
\right)$ as follows
\[
\alpha_{\mathbcal{A}}
\left\lbrace
\begin{array}{ccc}
I&\mor&
\mathrm{Sl}\left(
\mbox{\upshape{P{\l}}}_{\Sigma^{\neq 0}}(\mathbcal{A}), D^{\boldsymbol{\mathcal{A}}}
\right)
\\
i&\longmapsto&
[(z_{i},i)]_{\Phi^{D^{\mathbcal{A}}}}
\end{array}
\right.
\]

\begin{claim}\label{CNat1} $\alpha_{\mathbcal{A}}$ is a well-defined mapping that does not depend on the choice of the element $z_{i}$, for every $i\in I$.
\end{claim}

Let $i\in I$ and let $z'_{i}$ be any element in $A_{i}$. According to Claim~\ref{cgr}, we have that 
$
((z_{i},i),
(z'_{i},i)
)
\in \Phi^{D^{\mathbcal{A}}}
$, then the assignment $\alpha_{\mathbcal{A}}(i)$ does not depend on the choice of  the element $z_{i}$, for every $i\in I$.

This proves Claim~\ref{CNat1}.

\begin{claim}\label{CNat2} $\alpha_{\mathbcal{A}}$ is a morphism of sup-semilattices of the form
\[
\alpha_{\mathbcal{A}}\colon
\mathbf{I}
\mor
\mathbf{Sl}\left(
\mbox{\upshape{P{\l}}}_{\Sigma^{\neq 0}}(\mathbcal{A}), D^{\boldsymbol{\mathcal{A}}}
\right).
\]
\end{claim}

Let $i$ and $j$ be two elements in $I$. The following chain of equalities holds
\allowdisplaybreaks
\begin{align*}
\alpha_{\mathbcal{A}}(i\vee j)&=
[(z_{i\vee j},i\vee j)]_{\Phi^{D^{\mathbcal{A}}}}
\tag{Def. $\alpha_{\mathbcal{A}}$}
\\&=
[(f_{i,i\vee j}(z_{i}),i\vee j)]_{\Phi^{D^{\mathbcal{A}}}}
\tag{Claim~\ref{CNat1}}
\\&=
[
D^{\mathbcal{A}}(
(z_{i},i),
(z_{j},j)
)
]_{\Phi^{D^{\mathbcal{A}}}}
\tag{Def.~$D^{\mathbcal{A}}$}
\\&=
[(z_{i},i)]_{\Phi^{D^{\mathbcal{A}}}}
\vee^{D^{\mathbcal{A}}}
[(z_{j},j)]_{\Phi^{D^{\mathbcal{A}}}}
\tag{Claim~\ref{COrdII}}
\\&=
\alpha_{\mathbcal{A}}(i)
\vee ^{D^{\mathbcal{A}}}
\alpha_{\mathbcal{A}}(j).
\tag{Def. $\alpha_{\mathbcal{A}}$}
\end{align*}

This proves Claim~\ref{CNat2}.

To obtain $v_{\mathbcal{A}}$ we proceed as follows. For every $i\in I$, we will denote by 
$\boldsymbol{\alpha}_{\mathbcal{A}}(i)$ the $\Sigma^{\neq 0}$-algebra whose underlying set is 
$\alpha_{\mathbcal{A}}(i)$, which, by Definition~\ref{DConsSub}, is a subalgebra of 
$\mbox{\bfseries{\upshape{P{\l}}}}_{\Sigma^{\neq 0}}(\boldsymbol{\mathcal{A}})$. 

Let $i\in I$ and let us define the mapping $v_{\mathbcal{A},i}$ as follows
\[
v_{\mathbcal{A},i}
\left\lbrace
\begin{array}{ccc}
A_{i}&\mor&{\alpha}_{\mathbcal{A}}(i)
\\
z&\longmapsto&(z,i)
\end{array}
\right.
\] 

\begin{claim}\label{CNat3} The mapping $v_{\mathbcal{A},i}$ is well-defined. 
\end{claim}

Note that, for every $i\in I$ and every $z\in A_{i}$, $v_{\mathbcal{A},i}(z)=(z,i)$ is an element in $[(z_{i},i)]_{\Phi^{D^{\mathbcal{A}}}}$. This last statement follows directly from Claim~\ref{cgr}.

This proves Claim~\ref{CNat3}.

We next prove that the just defined mapping $v_{\mathbcal{A},i}$ is a $\Sigma^{\neq 0}$-homomorphism.

\begin{claim}\label{CNat4} For every $i\in I$, $v_{\mathbcal{A},i}$ is a $\Sigma^{\neq 0}$-homomorphism of the form
\[
v_{\mathbcal{A},i}\colon 
\mathbf{A}_{i}
\mor
\boldsymbol{\alpha}_{\mathbcal{A}}(i).
\]
\end{claim}

Let $n$ be an element of $\mathbb{N}-1$, $\sigma$ an operation symbol in $\Sigma_{n}$ and 
$(z_{j})_{j\in n} \in A_{i}^{n}$. The following chain of equalities holds
\allowdisplaybreaks
\begin{align*}
v_{\mathbcal{A},i}(
F^{\mathbf{A}_{i}}_{\sigma}(
(
z_{j}
)_{j\in n}
)
)
&=
(
F^{\mathbf{A}_{i}}_{\sigma}(
(
z_{j}
)_{j\in n}
),
i
)
\tag{Def.~$v_{\mathbcal{A},i}$}
\\&=
(
F^{\mathbf{A}_{i}}_{\sigma}(
(
f_{i,i}(
z_{j}
)
)_{j\in n}
),
i
)
\tag{Def~\ref{DIndSys}}
\\&=
\textstyle
(
F^{\mathbf{A}_{\bigvee_{j\in n}i}}_{\sigma}(
(
f_{i,\bigvee_{j\in n}i}(
z_{j}
)
)_{j\in n}
),
\bigvee_{j\in n}i
)
\tag{Def.~\ref{DSup}}
\\&=
F^{\mbox{\bfseries{\upshape{\scriptsize P{\l}}}}_{\Sigma^{\neq 0}}(\mathbcal{A})}_{\sigma}(
(z_{j},i)_{j\in n}
)
\tag{Def.~$\mbox{\bfseries{\upshape{P{\l}}}}_{\Sigma^{\neq 0}}(\mathbcal{A})$}
\\&=
F^{\boldsymbol{\alpha}_{\mathbcal{A}}(i)}
_{\sigma}(
(z_{j},i)_{j\in n}
)
\tag{Def.~\ref{DConsSub}}
\\&=
F^{\boldsymbol{\alpha}_{\mathbcal{A}}(i)}
_{\sigma}(
v_{\mathbcal{A},i}(
(z_{j}))_{j\in n}
),
\tag{Def.~$v_{\mathbcal{A},i}$}
\end{align*}
where $\bigvee_{j\in n}i$ is the supremum of the family $(i)_{j\in n}\in A_{i}^{n}$ which is constantly $i$.

This proves Claim~\ref{CNat4}.

%

Let $v_{\mathbcal{A}}$ be $(v_{\mathbcal{A},i})_{i\in I}$. In the next claim, we prove that $v_{\mathbcal{A}}$ is a morphism from 
$\mathcal{A}$ to 
$\mathcal{C}(\mbox{\bfseries{\upshape{P{\l}}}}_{\Sigma^{\neq 0}}(\mathbcal{A}), D^{\boldsymbol{\mathcal{A}}})_{\alpha_{\mathbcal{A}}}$.

%
%
%
%
%

\begin{claim}\label{CNat5} For every $i,j\in I$, if $(i,j)\in\leq$, then the following equality holds
\[
v_{\mathbcal{A},j}\circ f_{i,j}=
g_{{\alpha}_{\mathbcal{A}}(i),{\alpha}_{\mathbcal{A}}(j)}\circ v_{\mathbcal{A},i}.
\]
\end{claim}

Let $z$ be an element in $A_{i}$. The following chain of equalities holds.
\allowdisplaybreaks
\begin{align*}
v_{\mathbcal{A},j}(
f_{i,j}(z
)
)
&=
(
f_{i,j}(z
),j
)
\tag{Def.~$v_{\mathbcal{A},j}$}
\\&=
(
f_{i,i\vee j}(z
),
i\vee j
)
\tag{$(i,j)\in \leq$}
\\&=
D^{\mathbcal{A}}(
(z,i),
(z_{j},j)
)
\tag{Def.~$D^{\mathbcal{A}}$}
\\&=
g_{{\alpha}_{\mathbcal{A}}(i),{\alpha}_{\mathbcal{A}}(j)}(
z,i
)
\tag{Prop.~\ref{PConsHom} $\And$ Def.~${\alpha}_{\mathbcal{A}}(j)$}
\\&=
g_{{\alpha}_{\mathbcal{A}}(i),{\alpha}_{\mathbcal{A}}(j)}(
v_{\mathbcal{A},i}(
z
)
).
\tag{Def.~$v_{\mathbcal{A},i}$}
\end{align*}

This proves Claim~\ref{CNat5}.

Let $\eta_{\mathbcal{A}}$ be $(\alpha_{\mathbcal{A}}, v_{\mathbcal{A}})$.

\begin{claim}\label{CNat6} $\eta^{\mathbcal{A}}$ is a morphism from $\mathbcal{A}$ to 
$\mathrm{Is}_{\Sigma^{\neq 0}}(\mbox{\bfseries{\upshape{P{\l}}}}_{\Sigma^{\neq 0}}(\mathbcal{A}), D^{\mathbcal{A}})$.
\end{claim}

This follows from Claims~\ref{CNat1},~\ref{CNat2},~\ref{CNat3},~\ref{CNat4} and~\ref{CNat5}.

This proves Claim~\ref{CNat6}.

\begin{figure}
\begin{center}
\begin{tikzpicture}
[ACliment/.style={-{To [angle'=45, length=5.75pt, width=4pt, round]}},scale=1,
AClimentD/.style={double equal sign distance,
-implies
}]

\node[] (P) at (0,0) [] {$\left(
\mbox{\bfseries{\upshape{P{\l}}}}_{\Sigma^{\neq 0}}(\mathbcal{A}), D^{\mathbcal{A}}
\right)$};
\node[] (B) at (0,-2) [] {$(\mathbf{B},E)$};

\node[] (A) at (3,0) [] {$\mathbcal{A}$};
\node[] (IPA) at (7,0) [] {$\mathrm{Is}_{\Sigma^{\neq 0}}(\mbox{\bfseries{\upshape{P{\l}}}}_{\Sigma^{\neq 0}}(\mathbcal{A}), D^{\mathbcal{A}})$};
\node[] (IB) at (7,-2) [] {$\mathrm{Is}_{\Sigma^{\neq 0}}(\mathbf{B},E)$};

\draw[ACliment] (P) to node [left] {$(\xi, u)^{\sharp}$} (B);
\draw[ACliment] (A) to node [above] {$\eta_{\mathbcal{A}}$} (IPA);
\draw[ACliment] (A) to node [left, pos=0.6] {$(\xi, u)\quad$} (IB);
\draw[ACliment] (IPA) to node [right] {$\mathrm{Is}_{\Sigma^{\neq 0}}((\xi, u)^{\sharp})$} (IB);
\end{tikzpicture}
\end{center}
\caption{The Universal Property in Proposition~\ref{LAdj}.}
\label{FUnivEta}
\end{figure}

\textsf{Universal property.}

Finally, we prove the universal property. The reader is advised to consult the diagram in Figure~\ref{FUnivEta}. Let $(\mathbf{B},E)$ be a P{\l}onka $\Sigma^{\neq 0}$-algebra  and $(\xi, u)$ a morphism from $\mathbcal{A}$ to $\mathrm{Is}_{\Sigma^{\neq 0}}(\mathbf{B},E)$. Then we want to prove that there exists a unique morphism $(\xi, u)^{\sharp}$ in $\mbox{\sffamily{\upshape{P{\l}Alg}}}(\Sigma^{\neq 0})$  of the form
\[
(\xi, u)^{\sharp}\colon 
\left(
\mbox{\bfseries{\upshape{P{\l}}}}_{\Sigma^{\neq 0}}(\mathbcal{A}), D^{\mathbcal{A}}
\right)
\mor 
(\mathbf{B},E)
\]
such that $(\xi, u)=\mathrm{Is}_{\Sigma^{\neq 0}}((\xi, u)^{\sharp})\circ \eta_{\mathbcal{A}}$.

To this end, we first recall that, by Definition~\ref{DIS}, $\mathrm{Is}_{\Sigma^{\neq 0}}(\mathbf{B},E)$ is the ordered pair
$
(\mathbf{Sl}(B,E),\mathcal{C}(\mathbf{B},D)),
$ 
where, by Proposition~\ref{PUnivLnb}, $\mathbf{Sl}(B,E)$ is the sup-semilattice
$
(B/{\Phi^{E}},\leq^{E}),
$
and $\mathcal{C}(\mathbf{B},E)$ the $\mathbf{Sl}(B,E)$-inductive system of $\Sigma^{\neq 0}$-algebras given by
\[
(
(
\boldsymbol{\mathfrak{A}}
)_{\mathfrak{A}\in\mathrm{Sl}(B,E)},
(
g_{\mathfrak{A},\mathfrak{B}}
)_{(\mathfrak{A},\mathfrak{B})\in \leq^{E}}
).
\]

By Definition~\ref{DIsys}, the morphism $(\xi,u)$ from $\mathbcal{A}$ to $\mathrm{Is}_{\Sigma^{\neq 0}}(\mathbf{B},E)$ consists of a morphism of sup-semilattices
$
\xi\colon \mathbf{I}\mor \mathbf{Sl}(B,E)
$ 
and a family $u=(u_{i})_{i\in I}$ where, for every $i\in I$, $u_{i}$ is  $\Sigma^{\neq 0}$-homomorphism of the form
$
u_{i}\colon \mathbf{A}_{i}\mor 
\boldsymbol{\xi}(i)$, 
where, by Definition~\ref{DConsSub}, $\boldsymbol{\xi}(i)$ stands for the $\Sigma^{\neq 0}$-algebra whose underlying set is $\xi(i)$ which, we recall, is a subalgebra of $\mathbf{B}$. Moreover, the family $u=(u_{i})_{i\in I}$ is such that, for every $i,j\in I$ with $(i,j)\in\leq$, the following equality holds
\begin{align*}
u_{j}\circ f_{i,j}=g_{\xi(i),\xi(j)}\circ u_{i}.
\tag{E1}\label{E1}
\end{align*}

For every $i \in I$, consider the mapping $w_{i}=\mathrm{in}_{\boldsymbol{\xi}(i), \mathbf{B}}\circ u_{i}$, which is a  $\Sigma^{\neq 0}$-homomorphism of the form $w_{i}\colon \mathbf{A}_{i}\mor \mathbf{B}$. By the universal property of the coproduct there exists a unique mapping
\[
\textstyle
[(w_{i})_{i\in I}]\colon \coprod_{i\in I}A_{i}\mor B
\]
such that, for every $i\in I$, $[(w_{i})_{i\in I}]\circ \mathrm{in}_{i}=w_{i}$, where $\mathrm{in}_{i}$ stands for the canonical inclusion of $A_{i}$ in $\coprod_{i\in I}A_{i}$. Let 
$(\xi,u)^{\sharp}$ be $[(w_{i})_{i\in I}]$. We will present a series of claims to conclude that this is the desired universal morphism.

The following claim is introduced to avoid further calculations.
\begin{claim}\label{CUniv0} For every $(x,j)$ in $\coprod_{i\in I}A_{i}$, the following equality holds
\[
(\xi,u)^{\sharp}(x,j)=u_{j}(x).
\]
\end{claim}
The following chain of equality holds.
\allowdisplaybreaks
\begin{align*}
(\xi,u)^{\sharp}(x,j)&=
[(w_{i})_{i\in I}](x,j)
\tag{Def.~$(\xi,u)^{\sharp}$}
\\&=
w_{j}(x)
\tag{Def.~$[(w_{i})_{i\in I}]$}
\\&=
\mathrm{in}_{\boldsymbol{\xi}(j),\mathbf{B}}\left(u_{j}(x)\right)
\tag{Def.~$w_{j}$}
\\&=
u_{j}(x).
\tag{Def.~$\mathrm{in}_{\boldsymbol{\xi}(j),\mathbf{B}}$}
\end{align*}

This proves Claim~\ref{CUniv0}.

In order to be able to state that $(\xi,u)^{\sharp}$ is the desired morphism, we have to show that it satisfies certain properties that we set out below. We start by proving that $(\xi,u)^{\sharp}$ is a $\Sigma^{\neq 0}$-homomorphism from the P{\l}onka sum determined by $\mathbcal{A}$ to $\mathbf{B}$.

\begin{claim}\label{CUniv1} The mapping $(\xi,u)^{\sharp}$ determines a $\Sigma^{\neq 0}$-homomorphism from $\mbox{\bfseries{\upshape{P{\l}}}}_{\Sigma^{\neq 0}}(\mathbcal{A})$ to 
$\mathbf{B}$
\end{claim}

Let $n$ be an element of $\mathbb{N}-1$, $\sigma\in \Sigma_{n}$ and $(x_{j},i_{j})_{j\in n}\in (\coprod_{i\in I}A_{i})^{n}$. Note that, for every $j\in n$,  $u_{i_{j}}(x_{j})\in \xi(i_{j})$. Thus, for every $j\in n$, $\xi(i_{j})=[u^{i_{j}}(x_{j})]_{\Phi^{E}}$.  In virtue of Claim~\ref{COrdII} and taking into account the notation introduced in Definition~\ref{DIt}, we have that
\begin{align*}
\textstyle
\bigvee_{j\in n}\xi(i_{j})=[
E_{n}(
(u_{i_{j}}(x_{j}))_{j\in n}
)
]_{\Phi^{E}}.
\tag{E2}\label{E2}
\end{align*}

Then the following chain of equalities holds
\begin{flushleft}
$(\xi,u)^{\sharp}(
F^{\scriptsize{\mbox{\bfseries{\upshape{P{\l}}}}_{\Sigma^{\neq 0}}(\mathbcal{A})}}_{\sigma}(
(x_{j},i_{j})_{j\in n}
)
)$
\allowdisplaybreaks
\begin{align*}
&=
\textstyle
(\xi,u)^{\sharp}(
F^{\mathbf{A}_{\bigvee_{j\in n}i_{j}}}_{\sigma}(
(
f_{i_{j}, \bigvee_{j\in n}i_{j}}(x_{j})
)_{j\in n}
),
\bigvee_{j\in n}i_{j}
)
\tag{Def.~$\mbox{\bfseries{\upshape{P{\l}}}}_{\Sigma^{\neq 0}}(\mathbcal{A})$}
\\&=
u_{\bigvee_{j\in n}i_{j}}(
F^{\mathbf{A}_{\bigvee_{j\in n}i_{j}}}_{\sigma}(
(
f_{i_{j}, \bigvee_{j\in n}i_{j}}(x_{j})
)_{j\in n}
)
)
\tag{Claim~\ref{CUniv0}}
\\&=
F^{\mathbf{B}}_{\sigma}(
(
u_{\bigvee_{j\in n}i_{j}}(
f_{i_{j}, \bigvee_{j\in n}i_{j}}(x_{j})
)
)_{j\in n}
)
\tag{$\Sigma^{\neq 0}$-hom.}
\\&=
F^{\mathbf{B}}_{\sigma}(
(
g_{\xi(i_{j}),\xi(\bigvee_{j\in n}i_{j})}(
u_{i_{j}}(x_{j})
)
)_{j\in n}
)
\tag{by~\ref{E1}}
\\&=
F^{\mathbf{B}}_{\sigma}(
(
g_{\xi(i_{j}),\bigvee_{j\in n}\xi(i_{j})}(
u_{i_{j}}(x_{j})
)
)_{j\in n}
)
\tag{Def.~\ref{DSup}}
\\&=
F^{\mathbf{B}}_{\sigma}(
(
E(
u_{i_{j}}(x_{j}),
E_{n}((u_{i_{j}}(x_{j}))_{j\in n})
)
)_{j\in n}
)
\tag{Prop.~\ref{PConsHom} $\And$ by~\ref{E2}}
\\&=
E(
F^{\mathbf{B}}_{\sigma}(
(
u_{i_{j}}(x_{j})
)_{j\in n}
)
,
F^{\mathbf{B}}_{\sigma}(
(
E(
(u_{i_{j}}(x_{j}))_{j\in n}
)
)_{j\in n}
)
)
\tag{Prop.~\ref{PDHom}}
\\&=
F^{\mathbf{B}}_{\sigma}(
(
u_{i_{j}}(x_{j})
)_{j\in n}
)
\tag{Prop.~\ref{CTech4}}
\\&=
F^{\mathbf{B}}_{\sigma}(
(
(\xi,u)^{\sharp}(
x_{j},i_{j}
)
)_{j\in n}
).
\tag{Claim~\ref{CUniv0}}
\end{align*}
\end{flushleft}

This proves Claim~\ref{CUniv1}.

%
%
%

\begin{claim}\label{CUniv2} The equality
$E\circ ((\xi,u)^{\sharp}\times (\xi,u)^{\sharp})=(\xi,u)^{\sharp}\circ D^{\mathbcal{A}}$
holds.
\end{claim}

Let $(x,j)$ and $(y,k)$ be two elements of $\coprod_{i\in I}A_{i}$.
Then we have that $u_{j}(x)\in \xi(j)$ and $u_{k}(y)\in \xi(k)$. From this, we infer that $\xi(j)=[u_{j}(x)]_{\Phi^{E}}$ and $\xi(k)=[u_{k}(y)]_{\Phi^{E}}$. Then, by Claim~\ref{COrdII}, we have that 
\begin{align*}
\xi(j)\vee^{E} \xi(k)&=
\left[
E\left(
u_{j}(x),
u_{k}(y)
\right)
\right]_{\Phi^{E}}.
\tag{E3}\label{E3}
\end{align*}

The following chain of equalities holds.
\begin{flushleft}
$E(((\xi,u)^{\sharp}\times (\xi,u)^{\sharp})(
(x,j),(y,k)
))
$
\allowdisplaybreaks
\begin{align*}
\qquad&=
E(
(\xi,u)^{\sharp}(x,j),
(\xi,u)^{\sharp}(y,k)
)
\tag{Def.~$((\xi,u)^{\sharp})^{2}$}
\\&=
E(
u_{j}(x),
u_{k}(y)
)
\tag{Claim~\ref{CUniv0}}
\\&=
E\left(
E\left(u_{j}(x),u_{j}(x)\right),
u_{k}(y)
\right)
\tag{by~\ref{D1}}
\\&=
E(
u_{j}(x),
E(u_{j}(x),
u_{k}(y))
)
\tag{by~\ref{D2}}
\\&=
g_{\xi(j),\xi(j\vee k)}(
u_{j}(x)
)
\tag{Prop.~\ref{PConsHom} $\And$ by~\ref{E3}}
\\&=
u_{j\vee k}(
f_{j,j\vee k}(
x
)
)
\tag{by~\ref{E1}}
\\&=
(\xi,u)^{\sharp}(
f_{j,j\vee k}(
x
),
j\vee k
)
\tag{Claim~\ref{CUniv0}}
\\&=
(\xi,u)^{\sharp}(
D^{\mathbcal{A}}(
(x,j),(y,k)
)
)
\tag{Def.~$D^{\mathbcal{A}}$}
\end{align*}
\end{flushleft}

This proves Claim~\ref{CUniv2}.

The last two claims entail the next one.

\begin{claim}\label{CUniv3} The mapping
$(\xi,u)^{\sharp}$ is a morphism in  $\mbox{\sffamily{\upshape{P{\l}Alg}}}(\Sigma^{\neq 0})$ of the form
\[
(\xi, u)^{\sharp}\colon 
\left(
\mbox{\bfseries{\upshape{P{\l}}}}_{\Sigma^{\neq 0}}(\mathbcal{A}), D^{\mathbcal{A}}
\right)
\mor 
(\mathbf{B},E)
\]
\end{claim}

It follows from Claims~\ref{CUniv1} and~\ref{CUniv2}. 

This proves Claim~\ref{CUniv3}.

Now, by Proposition~\ref{PConsMorII}, $\mathrm{Is}_{\Sigma^{\neq 0}}((\xi,u)^{\sharp})$ is the morphism in  $\int^{\mathsf{Ssl}}\mathrm{Isys}_{\Sigma^{\neq 0}}$ from
$
(\mathbf{Sl}(\mbox{\upshape{P{\l}}}_{\Sigma^{\neq 0}}(\mathbcal{A}), D^{\mathbcal{A}}),\mathcal{C}(\mbox{\bfseries{\upshape{P{\l}}}}_{\Sigma^{\neq 0}}(\mathbcal{A}), D^{\boldsymbol{\mathcal{A}}}))
$
to
$
(\mathbf{Sl}(B,E),\mathcal{C}(\mathbf{B},E))
$
given by the pair 
$
(
\mathrm{Sl}((\xi,u)^{\sharp}),
t(
(\xi,u)^{\sharp}
)
)$, 
where, we recall, by Proposition~\ref{PUnivLnb}, $\mathrm{Sl}((\xi,u)^{\sharp})$ is the morphism of sup-semilattices
\[
\mathrm{Sl}((\xi,u)^{\sharp})
\left\lbrace
\begin{array}{ccl}
\mathbf{Sl}\left(
\mbox{\upshape{P{\l}}}_{\Sigma^{\neq 0}}(\mathbcal{A}), D^{\mathbcal{A}}
\right)
&\mor&
\mathbf{Sl}(B,E)
\\
{[(x,j)]_{\Phi^{D^{\mathbcal{A}}}}}
&\longmapsto&
[(\xi,u)^{\sharp}\left(x,j\right)]_{\Phi^{E}}
\end{array}
\right.
\]

\begin{claim}\label{CUniv4} For every equivalence class $\mathfrak{X}=[(x,j)]_{\Phi^{D^{\mathbcal{A}}}}$ in $\mathrm{Sl}(
\mbox{\upshape{P{\l}}}_{\Sigma^{\neq 0}}(\mathbcal{A}), D^{\mathbcal{A}}
)$, the following equality holds
\[
\mathrm{Sl}((\xi,u)^{\sharp})(\mathfrak{X})=\xi(j).
\]
\end{claim}
The following chain of equalities hold
\allowdisplaybreaks
\begin{align*}
\mathrm{Sl}((\xi,u)^{\sharp})(\mathfrak{X})&=
\mathrm{Sl}((\xi,u)^{\sharp})([(x,j)]_{\Phi^{D^{\mathbcal{A}}}})
\tag{Def.~$\mathfrak{X}$}
\\&=
[(\xi,u)^{\sharp}(x,j)]_{\Phi^{E}}
\tag{Def.~$\mathrm{Sl}((\xi,u)^{\sharp})$}
\\&=
[u_{j}(x)]_{\Phi^{E}}
\tag{Claim~\ref{CUniv0}}
\\&=
\xi(j).
\tag{$u_{j}(x)\in \xi(j)$}
\end{align*}

This proves Claim~\ref{CUniv4}.

All in all, we conclude that $\mathrm{Sl}((\xi,u)^{\sharp})$ is the morphism of sup-semilattices
\[
\mathrm{Sl}((\xi,u)^{\sharp})
\left\lbrace
\begin{array}{ccl}
\mathbf{Sl}(
\mbox{\upshape{P{\l}}}_{\Sigma^{\neq 0}}(\mathbcal{A}), D^{\mathbcal{A}}
)
&\mor&
\mathbf{Sl}(B,E)
\\
{[(x,j)]_{\Phi^{D^{\mathbcal{A}}}}}
&\longmapsto&
\xi(j)
\end{array}
\right.
\tag{E4}\label{E4}
\]

Let us recall that, by Definition~\ref{DIsys}, $\mathcal{C}(\mathbf{B},E)_{\mathrm{Sl}((\xi,u)^{\sharp})}$ is the $\mathbf{Sl}\left(
\mbox{\upshape{P{\l}}}_{\Sigma^{\neq 0}}(\mathbcal{A}), D^{\mathbcal{A}}
\right)$-inductive system given by
\[
(
(
(\xi,u)^{\sharp}[\boldsymbol{\mathfrak{X}}]
)_{\mathfrak{X}\in \mathrm{Sl}(
\scriptsize{\mbox{\upshape{P{\l}}}_{\Sigma^{\neq 0}}(\mathbcal{A})}, D^{\mathbcal{A}}
)},
(
g_{(\xi,u)^{\sharp}(\mathfrak{X}),(\xi,u)^{\sharp}(\mathfrak{Y})}
)_{(\mathfrak{X},\mathfrak{Y})\in\leq^{D^{\mathbcal{A}}}}
).
\]

By Proposition~\ref{PConsMorII}, $t((\xi,u)^{\sharp})$ is the morphism of $\mathbf{Sl}(
\mbox{\upshape{P{\l}}}_{\Sigma^{\neq 0}}(\mathbcal{A}), D^{\mathbcal{A}}
)$-inductive systems
\[
t((\xi,u)^{\sharp})\colon 
\mathcal{C}(
\mbox{\bfseries{\upshape{P{\l}}}}_{\Sigma^{\neq 0}}(\mathbcal{A}), D^{\mathbcal{A}}
)
\mor
\mathcal{C}(\mathbf{B},E)_{\mathrm{Sl}((\xi,u)^{\sharp})}
\]
given by the family 
$t((\xi,u)^{\sharp})=(t((\xi,u)^{\sharp})_{\mathfrak{X}})_{\mathfrak{X}\in \mathrm{Sl}(
\mbox{\bfseries{\upshape{\scriptsize P{\l}}}}_{\boldsymbol{\Sigma^{\neq 0}}}(\mathbcal{A}), D^{\mathbcal{A}}
)}$
where, for an equivalence class $\mathfrak{X}=[(x,j)]_{\Phi^{D^{\mathbcal{A}}}}$ in $\mathrm{Sl}(
\mbox{\upshape{P{\l}}}_{\Sigma^{\neq 0}}(\mathbcal{A}), D^{\mathbcal{A}}
)$, the $\Sigma^{\neq 0}$-homomorphism $t((\xi,u)^{\sharp})_{\mathfrak{X}}$ is given by
\[
t((\xi,u)^{\sharp})_{\mathfrak{X}}
\left\lbrace
\begin{array}{ccc}
 [(x,j)]_{\Phi^{D^{\mathbcal{A}}}}&\mor &
[
(\xi,u)^{\sharp}(x,j)
]_{\Phi^{E}}
\\
(z,j)&\mor &(\xi,u)^{\sharp}(z,j)
\end{array}
\right.
\]

Let us point out that in the above description we have taken into account Claim~\ref{cgr}. We can further simplify the presentation of $t\left((\xi,u)^{\sharp}\right)_{\mathfrak{X}}$ by considering Claims~\ref{CUniv0} and~\ref{CUniv4}. All in all the $\Sigma^{\neq 0}$-homomorphism $t((\xi,u)^{\sharp})_{\mathfrak{X}}$ is given by
\[
t((\xi,u)^{\sharp})_{\mathfrak{X}}
\left\lbrace
\begin{array}{ccl}
 \boldsymbol{\mathfrak{X}}&\mor &
\boldsymbol{\xi}(j)
\\
(z,j)&\mor &u_{j}(z)
\end{array}
\right.
\tag{E5}\label{E5}
\]

%

The following two claims will finally enable us to prove that
$$(\xi, u) = 
\mathrm{Is}_{\Sigma^{\neq 0}}((\xi, u)^{\sharp})\circ \eta_{\mathbcal{A}}.
$$ 

\begin{claim}\label{CUniv5} The equality $\xi=\mathrm{Sl}\left((\xi,u)^{\sharp}\right)\circ \alpha_{\mathbcal{A}}$ holds.
\end{claim}
Let $i$ be an element of $I$ and $x\in A_{i}$. Then the following chain of equalities holds.
\allowdisplaybreaks
\begin{align*}
\xi(i)&=\mathrm{Sl}((\xi,u)^{\sharp})(
[(x,i)]_{\Phi^{D^{\mathbcal{A}}}}
)
\tag{Claim~\ref{CUniv4}}
\\&=
\mathrm{Sl}((\xi,u)^{\sharp})( \alpha_{\mathbcal{A}}(i)).
\tag{Def.~$\alpha_{\mathbf{A}}$ $\And$ Claim~\ref{CNat1}}
\end{align*}

This proves Claim~\ref{CUniv5}.

\begin{claim}\label{CUniv6} For every $i\in I$, the equality $u_{i}=t((\xi,u)^{\sharp}
)_{\alpha_{\mathbcal{A}}(i)}\circ 
v_{\mathbcal{A},i}$ holds.
\end{claim}
Let $z$ be an element of $A_{i}$. Then the following chain of equalities holds.
\allowdisplaybreaks
\begin{align*}
t((\xi,u)^{\sharp}
)_{\alpha^{\mathbcal{A}}(i)}( 
v_{\mathbcal{A},i}(z))&=
t((\xi,u)^{\sharp}
)_{\alpha^{\mathbcal{A}}(i)}(
z,i
)
\tag{Def.~$v_{\mathbcal{A},i}$}
\\&=
u_{i}(z).
\tag{by~\ref{E5}}
\end{align*}

This proves Claim~\ref{CUniv6}.

We are now in position to prove the following equality.

\begin{claim}\label{CUniv7} The equality $(\xi, u)=\mathrm{Is}_{\Sigma^{\neq 0}}((\xi, u)^{\sharp})\circ \eta_{\mathbcal{A}}$ holds.
\end{claim}
The following chain of equalities holds.
\allowdisplaybreaks
\begin{align*}
(\xi, u)&=
(
\xi, (u_{i})_{i\in I}
)
\tag{$u=(u_{i})_{i\in I}$}
\\&=
(
\mathrm{Sl}((\xi,u)^{\sharp})\circ \alpha_{\mathbcal{A}},
(u_{i})_{i\in I}
)
\tag{Claim~\ref{CUniv5}}
\\&=
(
\mathrm{Sl}((\xi,u)^{\sharp})\circ \alpha_{\mathbcal{A}},
(
t(
(\xi,u)^{\sharp}
)_{\alpha^{\mathbcal{A}}(i)}\circ 
v_{\mathbcal{A},i}
)_{i\in I}
)
\tag{Claim~\ref{CUniv6}}
\\&=
(
\mathrm{Sl}((\xi,u)^{\sharp}),
t(
(\xi,u)^{\sharp}
)
)
\circ 
(
\alpha_{\mathbcal{A}}, v_{\mathbcal{A}}
)
\tag{Composition}
\\&=
(
\mathrm{Sl}((\xi,u)^{\sharp}),
t(
(\xi,u)^{\sharp}
)
)
\circ 
\eta_{\mathbcal{A}}
\tag{Def.~$\eta_{\mathbcal{A}}$}
\\&=
\mathrm{Is}_{\Sigma^{\neq 0}}((\xi, u)^{\sharp})\circ \eta_{\mathbcal{A}}
\tag{Def.~$\mathrm{Is}_{\Sigma^{\neq 0}}((\xi, u)^{\sharp})$}
\end{align*}

This proves Claim~\ref{CUniv7}.

To complete the proof it remains to prove the uniqueness of  $(\xi,u)^{\sharp}$.
\begin{claim}\label{CUniv8} Let $
h\colon 
\left(
\mbox{\bfseries{\upshape{P{\l}}}}_{\Sigma^{\neq 0}}(\mathbcal{A}), D^{\mathbcal{A}}
\right)
\mor 
(\mathbf{B},E)
$ be a morphism in $\mbox{\sffamily{\upshape{P{\l}Alg}}}(\Sigma^{\neq 0})$ 
such that $(\xi,u)=\mathrm{Is}_{\Sigma^{\neq 0}}(h)\circ \eta_{\mathbcal{A}}$, then $h=(\xi,u)^{\sharp}$.
\end{claim}

By Proposition~\ref{PIS}, $\mathrm{Is}_{\Sigma^{\neq 0}}(h)$ is the morphism in $\int^{\mathsf{Ssl}}\mathrm{Isys}_{\Sigma^{\neq 0}}$
$$
\mathrm{Is}_{\Sigma^{\neq 0}}(h)\colon 
(\mathbf{Sl}(\mbox{\upshape{P{\l}}}_{\Sigma^{\neq 0}}(\mathbcal{A}), D^{\mathbcal{A}}),\mathcal{C}(\mbox{\bfseries{\upshape{P{\l}}}}_{\Sigma^{\neq 0}}(\mathbcal{A}), D^{\boldsymbol{\mathcal{A}}}))
\mor
(\mathbf{Sl}(B,E),\mathcal{C}(\mathbf{B},E))
$$
given by
$
\mathrm{Is}_{\Sigma^{\neq 0}}(h)=
(
\mathrm{Sl}(h),
t(h)
)$ where, we recall, by Proposition~\ref{PConsMorII}, $t(h)$ is the morphism of $\mathbf{Sl}(
\mbox{\upshape{P{\l}}}_{\Sigma^{\neq 0}}(\mathbcal{A}), D^{\mathbcal{A}}
)$-inductive systems
\[
t(h)\colon 
\mathcal{C}(
\mbox{\bfseries{\upshape{P{\l}}}}_{\Sigma^{\neq 0}}(\mathbcal{A}), D^{\mathbcal{A}}
)
\mor
\mathcal{C}(\mathbf{B},E)_{\mathrm{Sl}((\xi,u)^{\sharp})}
\]
given by the family 
$t(h)=(t(h)_{\mathfrak{X}})_{\mathfrak{X}\in \mathrm{Sl}(\mbox{\upshape{{\scriptsize P{\l}}}}_{\Sigma^{\neq 0}}(\boldsymbol{\mathcal{A}}), D^{\boldsymbol{\mathcal{A}}}), D^{\mathbcal{A}}
)}$
where, for an equivalence class $\mathfrak{X}=[(x,j)]_{\Phi^{D^{\mathbcal{A}}}}$ in $\mathrm{Sl}(
\mbox{\upshape{P{\l}}}_{\Sigma^{\neq 0}}(\mathbcal{A}), D^{\mathbcal{A}}
)$, the $\Sigma^{\neq 0}$-homomorphism $t(h)_{\mathfrak{X}}$ is given by
\[
t(h)_{\mathfrak{X}}
\left\lbrace
\begin{array}{ccc}
 [(x,j)]_{\Phi^{D^{\mathbcal{A}}}}&\mor &
[
h(x,j)
]_{\Phi^{E}}
\\
(z,j)&\longmapsto &h(z,j)
\end{array}
\right.
\tag{E6}\label{E6}
\]

Let us unpack the equality $(\xi,u)=\mathrm{Is}_{\Sigma^{\neq 0}}(h)\circ \eta_{\mathbcal{A}}$.
\allowdisplaybreaks
\begin{align*}
(\xi,u)&=\mathrm{Is}_{\Sigma^{\neq 0}}(h)\circ \eta_{\mathbcal{A}}
\tag{Assumption}
\\&=
(
\mathrm{Sl}(h),
t(h)
)\circ \eta_{\mathbcal{A}}
\tag{Def.~$\mathrm{Is}_{\Sigma^{\neq 0}}(h)$}
\\&=
(
\mathrm{Sl}(h),
t(h)
)\circ (\eta_{\mathbcal{A}}, v_{\mathbcal{A}})
\tag{Def.~$\eta_{\mathbcal{A}}$}
\\&=
(
\mathrm{Sl}(h)\circ \eta_{\mathbcal{A}},
(
t(h)_{\alpha^{\mathbcal{A}}(i)}
\circ 
v_{\mathbcal{A},i}
)_{i\in I}
).
\tag{Composition}
\end{align*}

From this we conclude that, for every $i\in I$, the following equality holds
\[
u_{i}=t(h)_{\alpha^{\mathbcal{A}}(i)}
\circ 
v_{\mathbcal{A},i}.
\tag{E7}\label{E7}
\]

Let $(x,j)$ be an element of $\coprod_{i\in I}A_{i}$. Then the following chain of equalities holds.
\allowdisplaybreaks
\begin{align*}
h(x,j)&=h(v_{\mathbcal{A},j}(x))
\tag{Def.~$v_{\mathbcal{A},j}$}
\\&=
t(h)_{\alpha^{\mathbcal{A}}(j)}
(v_{\mathbcal{A},j}(x))
\tag{by~\ref{E6}}
\\&=
u_{j}(x)
\tag{by~\ref{E7}}
\\&=
(\xi,u)^{\sharp}(x,j).
\tag{Claim~\ref{CUniv0}}
\end{align*}

All in all, we conclude that $h=(\xi,u)^{\sharp}$.

This proves Claim~\ref{CUniv8}.

We denote by $\mbox{\upshape{P{\l}}}_{\Sigma^{\neq 0}}$ the left adjoint of $\mathrm{Is}_{\Sigma^{\neq 0}}$ and we call it \emph{the functor of P{\l}onka}. Moreover, for an inductive system of $\Sigma^{\neq 0}$-algebras $\mathbcal{A}$, we have that 
$\mbox{\upshape{P{\l}}}_{\Sigma^{\neq 0}}(\mathbcal{A})$ is the P{\l}onka $\Sigma^{\neq 0}$-algebra 
$(\mbox{\bfseries{\upshape{P{\l}}}}_{\Sigma^{\neq 0}}(\mathbcal{A}), D^{\mathbcal{A}})$; and for a morphism $(\xi,u)$ from an inductive system of $\Sigma^{\neq 0}$-algebras $\mathbcal{A}$ to another $\mathbcal{B}$, we have that $\mbox{\upshape{P{\l}}}_{\Sigma^{\neq 0}}(\xi,u)$ is $((\xi,u)\circ\eta_{\mathbcal{B}})^{\sharp}$, the unique morphism from 
$(\mbox{\bfseries{\upshape{P{\l}}}}_{\Sigma^{\neq 0}}(\mathbcal{A}), D^{\mathbcal{A}})$ to 
$(\mbox{\bfseries{\upshape{P{\l}}}}_{\Sigma^{\neq 0}}(\mathbcal{B}), D^{\mathbcal{B}})$ such that
$$
(\xi,u)\circ \eta_{\mathbcal{B}}  = ((\xi,u)\circ\eta_{\mathbcal{B}})^{\sharp}\circ \eta_{\mathbcal{A}}.
$$

This completes the proof of Proposition~\ref{LAdj}.
\end{proof}

\begin{proposition}
The functor $\mbox{\upshape{P{\l}}}_{\Sigma^{\neq 0}}$ from $\int^{\mathsf{Ssl}}\mathrm{Isys}_{\Sigma^{\neq 0}}$ 
to $\mbox{\sffamily{\upshape{P{\l}Alg}}}(\Sigma^{\neq 0})$ is essentially surjective. 
\end{proposition}

\begin{proof}
Let $(\mathbf{A},D)$ be a P{\l}onka $\Sigma^{\neq 0}$-algebra. Then every component $\varepsilon_{(\mathbf{A},D)}$ of the counit $\varepsilon$ of the adjunction 
$\mbox{\upshape{P{\l}}}_{\Sigma^{\neq 0}}\dashv \mathrm{Is}_{\Sigma^{\neq 0}}$ is an isomorphim from $\mbox{\upshape{P{\l}}}_{\Sigma^{\neq 0}}(\mathrm{Is}_{\Sigma^{\neq 0}}(\mathbf{A},D))$ to $(\mathbf{A},D)$. The details are left to the reader.
\end{proof}

\begin{figure}
\begin{center}
\begin{tikzpicture}
[ACliment/.style={-{To [angle'=45, length=5.75pt, width=4pt, round]}},scale=1,
AClimentD/.style={double equal sign distance,
-implies
}, scale=1.2]

\node[] (A) at (0,0) [] {$\int^{\mathsf{Ssl}}\mathrm{Isys_{\Sigma^{\neq 0}}}$};
\node[] (B) at (6,0) [] {$\mbox{\sffamily{\upshape{P{\l}Alg}}}(\Sigma^{\neq 0})$};
\node[] (C) at (0,-4) [] {$\mathsf{Ssl}$};
\node[] (D) at (6,-4) [] {$\mathsf{Alg}(\Sigma^{\neq 0})$};
\node[] (E) at (3,-2) [] {$\mathsf{Lnb}$};
\node[] (F) at (6,-2) [] {$\mathsf{Alg}(\Sigma^{\neq 0})\otimes\mathsf{Lnb}$};

\draw[ACliment] (B) to node [right] {$J_{\Sigma^{\neq 0}}$} (F);
\draw[ACliment] (F) to node [above] {$Q_{\Sigma^{\neq 0}}$} (E);
\draw[ACliment] (F) to node [right] {$P_{\Sigma^{\neq 0}}$} (D);
\draw[ACliment] (A) to node [midway, fill=white] {$\pi_{\mathrm{Isys}_{\Sigma^{\neq 0}}}$} (C);
\draw[ACliment, bend left] (C) to node [left, pos=.67] {$L_{\Sigma^{\neq 0}}$} (A);
\draw[ACliment, bend right] (C) to node [right, pos=.67] {$K_{\Sigma^{\neq 0}}$} (A);

\node[] () at (.3,-2.7) [] {$\dashv$};
\node[] () at (-.3,-2.7) [] {$\dashv$};

\node[] (c) at (0,-1.5) [fill=white] {$\qquad$};
\node[] (c) at (0,-1.35) [fill=white] {$\qquad$};
\node[] (c) at (0,-1.2) [fill=white] {$\qquad$};
\node[] (a) at (-.5,-1.5) [] {};
\node[] (b) at (.5,-1.5) [] {};
\draw[AClimentD] (a) to node [above] {$\gamma_{\Sigma^{\neq 0}}$} (b);

\node[] (G) at (3,-4) [rotate=180] {$\perp$};
\draw[ACliment] ($(G)+(-2.7,.2)$) to node [above, pos=.55] {$W_{\Sigma^{\neq 0}}$} ($(G)+(2.3,.2)$);
\draw[ACliment] ($(G)+(2.3,-.2)$) to node [below, pos=.45] {$M_{\Sigma^{\neq 0}}$} ($(G)+(-2.7,-.2)$);

\node[] (H) at (3,0) [] {$\perp$};
\draw[ACliment] ($(H)+(-2.1,.2)$) to node [above, pos=.53] {$\mbox{\upshape{P{\l}}}_{\Sigma^{\neq 0}}$} ($(H)+(2.05,.2)$);
\draw[ACliment] ($(H)+(2.05,-.2)$) to node [below, pos=.47] {$\mathrm{Is}_{\Sigma^{\neq 0}}$} ($(H)+(-2.1,-.2)$);

\node[] () at (1.6,-2.85) [rotate=38] {$\perp$};
\node[] (I) at (1.5,-3) [rotate=38] {};
\draw[ACliment, bend right, rotate=38] ($(I)+(1.4,.1)$) to node [above, pos=.5] {$\mathrm{Sl}\quad$} ($(I)+(-1,.1)$);
\draw[ACliment, bend right, rotate=38] ($(I)+(-1,-.1)$) to node [below, pos=.75] {$\quad\quad\qquad\mathrm{In}_{\mathsf{Ssl},\mathsf{Lnb}}$} ($(I)+(1.4,-.1)$);

\draw[ACliment, rounded corners] (D) -- ++(1.8,0) 
-- node[right, pos=0.5] {$V_{\Sigma^{\neq 0}}$}  ++(0,4) 
-- (B);

\node[] () at (7.4,-2) [] {$\vdash$};

\end{tikzpicture}
\end{center}
\caption{A general overview of the paper.}
\label{FOv}
\end{figure}

\begin{remark}
For two signatures $\Sigma$ and $\Lambda$ and a morphism $d$ from $\Sigma$ to $\Lambda$ (e.g., a derivor) we obtain the split indexed categories $(\mathsf{Ssl},\mathrm{Isys_{\Sigma}})$ and 
$(\mathsf{Ssl},\mathrm{Isys_{\Lambda}})$ and a morphism $(\mathrm{Id}_{\mathsf{Ssl}},\mathrm{Isys}_{d})$ from $(\mathsf{Ssl},\mathrm{Isys_{\Lambda}})$ to $(\mathsf{Ssl},\mathrm{Isys_{\Sigma}})$ (note that there exists a contravariant functor $\mathrm{Isys}$ from the category $\mathsf{Sig}_{\mathfrak{d}}$, of signatures and derivors between signatures, to the category $\mathsf{Cat}^{\mathsf{Ssl}^{\mathrm{op}}}$, that sends a signature $\Sigma$ to the contravariant functor $\mathrm{Isys_{\Sigma}}$ and a derivor $d$ from $\Sigma$ to $\Lambda$ to the natural transformation $\mathrm{Isys}_{d}$ from $\mathrm{Isys_{\Lambda}}$ to $\mathrm{Isys_{\Sigma}}$). Then a natural connection between the diagram in Figure~\ref{FOv} and the diagram of the corresponding figure for $\Lambda^{\neq 0}$ is obtained.
\end{remark}


\end{document}